\documentstyle[11pt]{article}

\author{Yu.N. Bespalov
\thanks{Institute for Theoretical Physics \
           Metrologichna str., 14-b \
           Kiev 143, 252143 Ukraine}
}
\title{Crossed modules and quantum groups in braided categories}
\date{}

\expandafter\chardef\csname pre amssym.def at\endcsname=\the\catcode`\@
\catcode`\@=11
\def\hexnumber@#1{\ifcase#1 0\or 1\or 2\or 3\or 4\or 5\or 6\or 7\or 8\or
 9\or A\or B\or C\or D\or E\or F\fi}
\font\tenmsa=msam10
\font\sevenmsa=msam7
\font\fivemsa=msam5
\newfam\msafam
\textfont\msafam=\tenmsa
\scriptfont\msafam=\sevenmsa
\scriptscriptfont\msafam=\fivemsa
\edef\msafam@{\hexnumber@\msafam}
\def\emptybox{\mathrel{\mathchar"0\msafam@03}}
\font\tenmsb=msbm10
\font\sevenmsb=msbm7
\font\fivemsb=msbm5
\newfam\msbfam
\textfont\msbfam=\tenmsb
\scriptfont\msbfam=\sevenmsb
\scriptscriptfont\msbfam=\fivemsb
\edef\msbfam@{\hexnumber@\msbfam}
\def\ltimes{\mathrel{\mathchar"0\msbfam@6E}}
\def\rtimes{\mathrel{\mathchar"0\msbfam@6F}}
\catcode`\@=\csname pre amssym.def at\endcsname
\newcommand{\subsec}{\subsection{}\vskip -21pt\hskip 1.5\parindent%
\setcounter{equation}{0}}
\newcommand{\sect}[1]{\section{#1}\setcounter{equation}{0}}

\newtheorem{theorem}{\rm THEOREM}[subsection]
\newtheorem{proposition}[theorem]{\rm PROPOSITION}
\newtheorem{corollary}[theorem]{\rm COROLLARY}
\newtheorem{lemma}[theorem]{\rm LEMMA}
\newtheorem{definition}{\rm DEFINITION}[subsection]

\newenvironment{proof}{\par\noindent{\it Proof.}}%
{\hfill$\emptybox$\newline\medskip}
\def\proof{\par\noindent{\it Proof. }}

\def\DY#1{{\cal DY}\left({#1}\right)}
\def\YD#1{{\cal DY}{#1}}
\def\cO#1{{{\cal O}(#1)}}
\def\RR{I\!\!R}

\def\ZZ{Z\!\!\!Z}

\def\C{{\cal C}}
\def\op{{\rm op}}
\def\id{{\rm id}}
\def\Obj{{\rm Obj}}

\def\krr{\kern -.16667em}%
\def\kr{}%
\def\krrr{\kern -.3\unitlength}%
\def\krl{}%
\def\krrrr%
{\kern -4.65\unitlength}
\newlength{\textwd}%
\def\hhstep{\kr\kr
\kern -.5\unitlength}
\def\hstep{\kr\kr
\kern .5\unitlength}
\def\step{\kr\kr
\kern \unitlength}
\def\Step{\kr\kr
\kern 2\unitlength}
\def\vvbox#1{{\offinterlineskip\vcenter{%
\def\coev{\kr
\begin{picture}(2,2)\put(1,0){\oval(2,2)[t]}\end{picture}}
\def\ev{\kr
\begin{picture}(2,2)\put(1,2){\oval(2,2)[b]}\end{picture}}
\def\hcoev{\kr
\begin{picture}(1,2)\put(.5,0){\oval(1,1)[t]}\end{picture}}
\def\hev{\kr
\begin{picture}(1,2)\put(.5,2){\oval(1,1)[b]}\end{picture}}
\def\COEV{\kr
\begin{picture}(2,2)\put(3,0){\oval(6,6)[t]}\end{picture}}
\def\EV{\kr
\begin{picture}(2,2)\put(3,2){\oval(6,6)[b]}\end{picture}}
\def\unit{\kr
\begin{picture}(0,2)
\put(0,0){\line(0,1){1}}\put(0,1.2){\circle{0.4}}
\end{picture}}
\def\counit{\kr
\begin{picture}(0,2)
\put(0,1){\line(0,1){1}}\put(0,.8){\circle{0.4}}
\end{picture}}
\def\Q##1{\kr
\begin{picture}(0,2)
\put(0,0){\line(0,1){0.4}}\put(0,1){\circle{1.2}}
\put(-0.6,0.4){\makebox(1.2,1.2)[cc]{$\scriptstyle ##1$}}
\end{picture}}
\def\O##1{\kr
\begin{picture}(0,2)
\put(0,0){\line(0,1){0.4}}\put(0,1.6){\line(0,1){0.4}}\put(0,1){\circle{1.2}}
\put(-0.6,0.4){\makebox(1.2,1.2)[cc]{$\scriptstyle ##1$}}
\end{picture}}
\def\S{\O{S}}                 \def\SS{\O{S^-}}
\def\tS{\O{\overline S}}     \def\tSS{\O{\overline S^-}}
\let\P\O
\def\dash##1{\kr
\begin{picture}(2,2)
\put(-.5,0){\dashbox{.1}(3,2){$\scriptstyle ##1$}}
\end{picture}}
\def\Dash##1{\kr
\begin{picture}(2,2)
\put(-1,0){\dashbox{.1}(4,2){$\scriptstyle ##1$}}
\end{picture}}
\def\DDash##1{\kr
\begin{picture}(2,2)
\put(-2,0){\dashbox{.1}(6,2){$\scriptstyle ##1$}}
\end{picture}}
\def\x{\kr
\begin{picture}(2,2)
\put(0,2){\line(1,-1){2}}\put(0,0){\line(1,1){.7}}\put(2,2){\line(-1,-1){.7}}
\end{picture}}
\def\xx{\kr
\begin{picture}(2,2)
\put(0,2){\line(1,-1){.7}}\put(0,0){\line(1,1){2}}\put(2,0){\line(-1,1){.7}}
\end{picture}}
\def\hx{\kr
\begin{picture}(1,2)
\put(0,2){\line(1,-2){1}}\put(0,0){\line(1,2){.35}}\put(1,2){\line(-1,-2){.35}}
\end{picture}}
\def\hxx{\kr
\begin{picture}(1,2)
\put(0,2){\line(1,-2){.35}}\put(0,0){\line(1,2){1}}\put(1,0){\line(-1,2){.35}}
\end{picture}}
\def\d{\kr
\begin{picture}(1,2)\put(0,2){\line(1,-2){1}}\end{picture}}
\def\dd{\kr
\begin{picture}(1,2)\put(0,0){\line(1,2){1}}\end{picture}}
\def\hd{\kr
\begin{picture}(1,2)
\put(0,2){\line(1,-2){.5}}
\put(.5,1){\line(0,-1){1}}
\end{picture}}
\def\hdd{\kr
\begin{picture}(1,2)
\put(1,2){\line(-1,-2){.5}}
\put(0,1){\line(0,-1){1}}
\end{picture}}
\def\ld{\kr
\begin{picture}(1,2)
\put(1,0){\oval(2,2)[lt]}\put(1,0){\line(0,1)2}
\end{picture}}
\def\Ld{\kr
\begin{picture}(2,2)
\put(2,0){\oval(4,2)[lt]}\put(2,0){\line(0,1)2}
\end{picture}}
\def\cd{\kr
\begin{picture}(2,2)
\put(1,0){\oval(2,2)[ct]}\put(1,1){\line(0,1)1}
\end{picture}}
\def\hdcd{\kr
\begin{picture}(1,2)
\put(0,2){\line(1,-2){.5}}
\put(.5,0){\oval(1,1)[ct]}\put(.5,.5){\line(0,1){.5}}
\end{picture}}
\def\hddcd{\kr
\begin{picture}(1,2)
\put(1,2){\line(-1,-2){.5}}
\put(.5,0){\oval(1,1)[ct]}\put(.5,.5){\line(0,1){.5}}
\end{picture}}
\def\hcd{\kr
\begin{picture}(1,2)
\put(.5,0){\oval(1,1)[ct]}\put(.5,.5){\line(0,1){1.5}}
\end{picture}}
\def\Cd{\kr
\begin{picture}(4,2)
\put(2,0){\oval(4,2)[ct]}\put(2,1){\line(0,1)1}
\end{picture}}
\def\rd{\kr
\begin{picture}(1,2)
\put(0,0){\oval(2,2)[rt]}\put(0,0){\line(0,1)2}
\end{picture}}
\def\Rd{\kr
\begin{picture}(2,2)
\put(0,0){\oval(4,2)[rt]}\put(0,0){\line(0,1)2}
\end{picture}}
\def\lu{\kr
\begin{picture}(1,2)
\put(1,2){\oval(2,2)[lb]}\put(1,0){\line(0,1)2}
\end{picture}}
\def\Lu{\kr
\begin{picture}(2,2)
\put(2,2){\oval(4,2)[lb]}\put(2,0){\line(0,1)2}
\end{picture}}
\def\cu{\kr
\begin{picture}(2,2)
\put(1,2){\oval(2,2)[cb]}\put(1,0){\line(0,1)1}
\end{picture}}
\def\hdcu{\kr
\begin{picture}(1,2)
\put(1,0){\line(-1,2){.5}}
\put(.5,2){\oval(1,1)[cb]}\put(.5,1){\line(0,1){.5}}
\end{picture}}
\def\hddcu{\kr
\begin{picture}(1,2)
\put(0,0){\line(1,2){.5}}
\put(.5,2){\oval(1,1)[cb]}\put(.5,1){\line(0,1){.5}}
\end{picture}}
\def\hcu{\kr
\begin{picture}(1,2)
\put(.5,2){\oval(1,1)[cb]}\put(.5,0){\line(0,1){1.5}}
\end{picture}}
\def\Cu{\kr
\begin{picture}(4,2)
\put(2,2){\oval(4,2)[cb]}\put(1,0){\line(0,1)1}
\end{picture}}
\def\ru{\kr
\begin{picture}(1,2)
\put(0,2){\oval(2,2)[rb]}\put(0,0){\line(0,1)2}
\end{picture}}
\def\Ru{\kr
\begin{picture}(2,2)
\put(0,2){\oval(4,2)[rb]}\put(0,0){\line(0,1)2}
\end{picture}}
\def\k{\kr
\begin{picture}(1,2)
\put(0,2){\oval(2,1)[rb]}
\put(0,0){\oval(2,1)[rt]}
\put(0,0){\line(0,1)2}
\end{picture}}
\def\kk{\kr
\begin{picture}(1,2)
\put(1,2){\oval(2,1)[lb]}
\put(1,0){\oval(2,1)[lt]}
\put(1,0){\line(0,1)2}
\end{picture}}
\def\ro##1{\kr
\begin{picture}(2,2)
\put(.4,0){\oval(.8,.8)[lt]}\put(1.6,0){\oval(.8,.8)[rt]}
\put(1,0.4){\circle{1.2}}
\put(0.4,-0.2){\makebox(1.2,1.2)[cc]{$\scriptstyle ##1$}}%
\end{picture}}
\def\coro##1{\kr
\begin{picture}(2,2)
\put(.4,2){\oval(.8,.8)[lb]}\put(1.6,2){\oval(.8,.8)[rb]}
\put(1,1.6){\circle{1.2}}
\put(0.4,1){\makebox(1.2,1.2)[cc]{$\scriptstyle ##1$}}%
\end{picture}}
\def\Ro##1{\kr
\begin{picture}(4,2)
\put(1.4,0){\oval(2.8,1.2)[lt]}\put(2.6,0){\oval(2.8,1.2)[rt]}
\put(2,.6){\circle{1.2}}
\put(1.4,0){\makebox(1.2,1.2)[cc]{$\scriptstyle ##1$}}%
\end{picture}}
\def\coRo##1{\kr
\begin{picture}(4,2)
\put(1.4,2){\oval(2.8,1.2)[lb]}\put(2.6,2){\oval(2.8,1.2)[rb]}
\put(2,1.4){\circle{1.2}}
\put(1.4,.8){\makebox(1.2,1.2)[cc]{$\scriptstyle ##1$}}%
\end{picture}}
\def\r{\ro{\cal R}}              \def\rr{\ro{{\cal R}^-}}
            \def\rrr{\ro{{\cal R}^{\tilde{}}}}
\def\ra{\ro{{\cal R}_A}}        \def\rra{\ro{{\cal R}^-_A}}
\def\rb{\ro{{\cal R}_B}}        \def\rrb{\ro{{\cal R}^-_B}}
\def\rh{\ro{{\cal R}_H}}
\def\R{\Ro{\cal R}}           \def\RR{\Ro{{\cal R}^-}}
\def\Ra{\Ro{{\cal R}_A}}        \def\RRa{\Ro{{\cal R}^-_A}}
\def\Rb{\Ro{{\cal R}_B}}        \def\RRb{\Ro{{\cal R}^-_B}}
\def\Rh{\Ro{{\cal R}_H}}
\def\tu##1{\kr
\begin{picture}(2,2)
\put(.4,2){\oval(.8,.8)[lb]}\put(1.6,2){\oval(.8,.8)[rb]}
\put(1,1.6){\circle{1.2}}
\put(0.4,1){\makebox(1.2,1.2)[cc]{$\scriptstyle ##1$}}%
\put(1,0){\line(0,1)1}
\end{picture}}
\def\id{\kr
\begin{picture}(0,2)\put(0,0){\line(0,1)2}\end{picture}}
\def\obj##1{\settowidth{\textwd}{$##1$}%
\raise .2\unitlength\hbox{\kern -.5\textwd $##1$ \kern -.5\textwd \krrr}}
\def\Obj##1{\settowidth{\textwd}{$##1$}%
\raise 1.1\unitlength\hbox{\kern -1\textwd $##1$}}
\def\hhbox##1{\hbox{\krrrr
\def\coev{\kr
\begin{picture}(1,1)\put(.5,0){\oval(1,1)[t]}\end{picture}}
\def\ev{\kr
\begin{picture}(1,1)\put(.5,1){\oval(1,1)[b]}\end{picture}}
\def\ld{\kr
\begin{picture}(1,1)
\put(1,0){\oval(2,2)[lt]}\put(1,0){\line(0,1)1}
\end{picture}}
\def\Ld{\kr
\begin{picture}(2,1)
\put(2,0){\oval(4,2)[lt]}\put(2,0){\line(0,1)1}
\end{picture}}
\def\rd{\kr
\begin{picture}(1,1)
\put(0,0){\oval(2,2)[rt]}\put(0,0){\line(0,1)1}
\end{picture}}
\def\Rd{\kr
\begin{picture}(2,1)
\put(0,0){\oval(4,2)[rt]}\put(0,0){\line(0,1)1}
\end{picture}}
\def\cd{\kr
\begin{picture}(1,1)
\put(.5,0){\oval(1,1)[ct]}\put(.5,.5){\line(0,1){.5}}
\end{picture}}
\def\lu{\kr
\begin{picture}(1,1)
\put(1,1){\oval(2,2)[lb]}\put(1,0){\line(0,1)1}
\end{picture}}
\def\Lu{\kr
\begin{picture}(2,1)
\put(2,1){\oval(4,2)[lb]}\put(2,0){\line(0,1)1}
\end{picture}}
\def\cu{\kr
\begin{picture}(1,1)
\put(.5,1){\oval(1,1)[cb]}\put(.5,0){\line(0,1){.5}}
\end{picture}}
\def\ru{\kr
\begin{picture}(1,1)
\put(0,1){\oval(2,2)[rb]}\put(0,0){\line(0,1)1}
\end{picture}}
\def\Ru{\kr
\begin{picture}(2,1)
\put(0,1){\oval(4,2)[rb]}\put(0,0){\line(0,1)1}
\end{picture}}
\def\hru{\kr
\begin{picture}(.5,1)
\put(0,1){\oval(1,1)[rb]}\put(0,0){\line(0,1)1}
\end{picture}}
\def\hlu{\kr
\begin{picture}(.5,1)
\put(.5,1){\oval(1,1)[lb]}\put(.5,0){\line(0,1)1}
\end{picture}}
\def\hrd{\kr
\begin{picture}(.5,1)
\put(0,0){\oval(1,1)[rt]}\put(0,0){\line(0,1)1}
\end{picture}}
\def\hld{\kr
\begin{picture}(.5,1)
\put(.5,0){\oval(1,1)[lt]}\put(.5,0){\line(0,1)1}
\end{picture}}
\def\id{\kr
\begin{picture}(0,1)\put(0,0){\line(0,1)1}\end{picture}}
\def\d{\kr
\begin{picture}(.5,1)\put(0,1){\line(1,-2){0.5}}\end{picture}}
\def\dd{\kr
\begin{picture}(.5,1)\put(0,0){\line(1,2){0.5}}\end{picture}}
##1}}#1}\normalbaselines}}
\def\object#1{\settowidth{\textwd}{$#1$}%
                        \hbox{%
                        \kern -.5\textwd $#1$ \kern -.5\textwd}}

\begin{document}
\maketitle
\begin{abstract}
Let $A$ be a Hopf algebra in a braided category $\cal C$.
Crossed modules over $A$ are introduced and studied as objects with both
module and comodule structures satisfying a compatibility condition.
The category $\DY{\cal C}^A_A$ of crossed modules is braided and is a
concrete realization of a known general construction of a double or
center of a monoidal category.  For a quantum braided group
$(A,\overline A,{\cal R})$
the corresponding braided category of modules
${\cal C}_{\cO{A,\overline A}}$ is
identified with a full subcategory in $\DY{\cal C}_A^A$.
The connection
with cross products is discussed and a suitable cross product in the
class of quantum braided groups is built.  Majid--Radford theorem, which
gives equivalent conditions for an ordinary Hopf algebra to be such a
cross product, is generalized to the braided category.  Majid's
bosonization theorem is also generalized.
\end{abstract}

\sect{Introduction}  

Crossed modules over a finite group $G$ arisen in topology \cite{Whitehead1}.
Crossed modules over groups and over Lie algebras were studied, as an
algebraic object, in several contexts. In particular, see \cite{Brown} in
connection with cohomology of groups.  A generalization for an arbitrary
Hopf algebra $A$ was noted by Yetter \cite{Yetter1}, who called the
structures "crossed bimodules".  This construction was extensively
studied by Radford and Towber \cite{RadTow} under the name
"Yetter-Drinfel'd structures". See also \cite{LR}.  A crossed module is a
vector space with both module and comodule structures over $A$ satisfying
a compatibility condition.  The category ${}_A^A\YD{}$ of crossed modules
is a convenient reformulation of the category of modules over Drinfel'd's
quantum double ${\cal D}(A)$ \cite{Drinfel'd1} and has the corresponding
braiding due to Drinfel'd.  This was explained in \cite{M1} where also a
functor embedding ${}_A{\cal M}\hookrightarrow{}_A^A\YD{}$ was introduced
in the case where $A$ is quasitriangular (a strict quantum group).  One
can obtain ${}_A^A\YD{}$ also as a 'center' or an 'inner double' of the
monoidal category ${}_A{\cal M}$ of left modules by a general
construction \cite{Majid9}.
Crossed modules appear in different contexts related to Hopf algebras
and quantum groups.
For example, a subalgebra of left or right invariant forms in a bicovariant
differential algebra over a Hopf algebra has a crossed module structure.

In this paper we introduce and study the category $\DY{\cal C}_A^A$
of crossed modules over a Hopf algebra $A$
living in an arbitrary braided monoidal category $\cal C$
({\em braided-Hopf algebra} or {\em braided group}
\cite{M2}-\cite{M4},
\cite{Lyubashenko1}-\cite{Lyubashenko4},
\cite{Majid6}-\cite{Majid10}).
This category is also braided and most of the results for ordinary Hopf
algebras hold in this situation.
We also define the category $\DY{\cal C}_{A,H}$ depending on bialgebra
pairing $\rho:\, A\otimes H\rightarrow\underline 1$, which is a fully
braided analog of the category of modules ${\cal M}_{{\cal D}(A,H,\rho)}$
over Drinfel'd double ${\cal D}(A,H,\rho)$, and discuss (in a special
case) the question when this category can be realized as a category of
modules over something.

Braided groups have been extensively studied over the last few years and
play an important role in $q$-deformed physics and mathematics
\cite{Majid8}-\cite{Majid10}.
Examples, applications and the basic theory of braided groups have been
introduced and developed by Majid
(some similar concepts arise independently
in \cite{Lyubashenko1,Lyubashenko4} inspired by results
on conformal field theory).
In particular, Majid \cite{Majid6}-\cite{Majid7} defined a quantum braided
group as a pair of a braided Hopf algebra $H$ with a quasitriangular
structure $\cal R$ satisfying axioms which are a generalization of ones
for an ordinary quantum group \cite{Drinfel'd1}, and a non-trivial
class $\cal O$ of modules over $H$.  He showed that the largest such
class ${\cal C}_{\cO{H}}$ is closed under tensor product and braided
\cite{Majid7}; we will see that it becomes a full subcategory of
$\DY{\cal C}^H_H$ in the same way as for usual quantum groups in
\cite{M1}. See also \cite{M5}.  This embedding of modules over a
quantum braided group into crossed modules is a key to applications to
certain cross products and the bosonization construction
\cite{Majid6},\cite{M6} and allows us to generalizes them to the fully
braided setting.

Let $A$ be a Hopf algebra in a braided category $\cal C$ and
$B$ a Hopf algebra in the category $\DY{\cal C}_A^A$ of crossed modules.
Then similarly to unbraided case
\cite{Majid10} the tensor product $A\otimes B$ can be equipped with a natural
cross product algebra and coalgebra $A\ltimes B$.
The Majid-Radford theorem \cite{M6} gives equivalent condition for
an ordinary Hopf algebra to be such a cross product;
we show that a braided variant of this theorem holds when idempotents
in a category $\cal C$ are split.
The last condition is not essential because
any braided category can be extended to one with split idempotents.

Similarly, let $(A,{\cal R}_A)$ be a quantum braided group in $\cal C$
and $(B,{\cal R}_B)$ a quantum braided group in ${\cal C}_{\cO{H}}$.
Then the cross product Hopf algebra $A\ltimes B$ in $\cal C$ is also
a quantum braided group with ${\cal R}_{A\ltimes B}$ built from ${\cal R}_A$
and ${\cal R}_B$.
This construction generalizes the Majid's bosonization procedure \cite{Majid6}
which was defined in the case when $A$ is an ordinary quantum group.
Equivalent conditions for a quantum braided group to be a such cross product
are also obtained.

Then we describe fully braided variant of Majid's transmutation procedure
\cite{M19,Majid7} and show that analog of Majid's result \cite{Majid6}
about relation between transmutation and bosonization is true in our
setting.

Finally, we suppose that a category $\cal C$ is balanced and define a
ribbon structure on a quantum braided group $(A,\overline A,{\cal R})$ in
$\C$ in a such way that the category ${\cal C}_{\cO{A,\overline A}}$ is
also ribbon.
We show that a ribbon structure is well-behaved under cross products and
transmutation.

Results of this paper where announced in the note \cite{B2}.
See also \cite{Bespalov1} where we work in the framework of ordinary Hopf
algebras. The theory of crossed modules which is developed here is used
to study of Hopf bimodules in \cite{BD}.
Constructions of \cite{B3} illustrate our general theory.

An outline of the paper is as follows.
In the preliminary Section 2 necessary notations connected with
braided categories and braided groups are recalled from
\cite{Majid8},\cite{Majid10}.
In Section 3 the categories of crossed modules over a bialgebra
(braided group) in a braided category $\cal C$ are introduced and studied.
Section 4 is devoted to certain cross products of braided Hopf algebras and
a braided variant of Majid-Radford theorem.
Majid's definition of quantum braided group from \cite{Majid7} are discussed
in Section 5.
Finally, in Section 6 the results of Section 4 are extended to include
cross products by quantum braided groups
(generalizing \cite{Majid6} for cross products by
braided groups with trivial ${\cal R}_A$), allowing us to prove a
generalized bosonization theorem.
Connections with transmutation and ribbon structure are discussed.
Diagrammatic proofs of theorems are moved in appendix.

        \sect{Preliminaries}                                             %

In this preliminary section we recall the basic notations and results of
the theory of braided groups from
\cite{Majid6}-\cite{Majid10}.
See also \cite{Lyubashenko1},\cite{Lyubashenko4},
and see \cite{Freyd1},\cite{JS} in connection with braided categories
themselves.
We assume that the reader is familiar with ordinary Hopf algebras
\cite{Sweedler1} and quantum groups \cite{Drinfel'd1}.

\subsec{}
Unless otherwise stated, we will suppose that
${\cal C}=({\cal C},\otimes ,\underline 1,\Psi)$
is {\em a braided (monoidal) category}
with tensor product  $\otimes$, unit object $\underline 1$
and braiding $\Psi$.
Without loss of generality by Mac Lane's coherence theorem \cite{MacLane}
we will assume that underlying monoidal category is strict, i.e. the functors
$\_\otimes (\_\otimes\_)$ and $(\_\otimes\_)\otimes\_$ coincide and
$\underline 1\otimes X=X=X\otimes\underline 1$.
A category $\cal C$ is called {\em pre-braided} if existence of $\Psi^{-1}$
is not postulated.

{\em A (braided) monoidal functor}
$F=(F,\lambda ):({\cal C},\otimes ,\underline 1,\Psi)\rightarrow
   ({\cal C}^\prime,\otimes^\prime ,\underline 1^\prime,\Psi^\prime)$
is a pair of a functor
  $F:{\cal C}\rightarrow{\cal C}^\prime$
and isomorphism of functors
 $\lambda :\,F(\_)\otimes F(\_)\rightarrow F(\_\otimes\_)$
which are compatible with braiding:
 $\lambda\circ\Psi^\prime=F(\Psi )\circ\lambda$.
We say that $F$ is {\em strict} if $\lambda$ is identity.

We actively use diagrammatic calculus in braided categories which is not a
trivial generalization of 'wiring diagrams' for usual linear algebra (as in
Penrose's spin networks and \cite{Yetter1}) developed by Majid
(\cite{Majid6},\cite{Majid8}).
Morphisms $\Psi$ and $\Psi^{-1}$ are represented by under and over crossing
and algebraic information 'flows' along braids and tangles according to
functoriality and the coherence theorem for braided categories \cite{JS}:
\begin{equation}
\Psi=\enspace
\vvbox{\hbox{\hx}}
\qquad\quad
\Psi^{-1}=\enspace
\vvbox{\hbox{\hxx}}
\qquad\qquad
\vvbox{\hbox{\O{f}\step\id}
       \hbox{\hx}}
\enspace =\enspace
\vvbox{\hbox{\hx}
       \hbox{\id\step\O{f}}}
\qquad\quad
\vvbox{\hbox{\id\step\O{f}}
       \hbox{\hx}}
\enspace =\enspace
\vvbox{\hbox{\hx}
       \hbox{\O{f}\step\id}}
\label{Psi}
\end{equation}

\begin{figure}
\begin{displaymath}
\eta_A={}\enspace
\vvbox{\hbox{\unit}}
\qquad\qquad
\mu_A={}\enspace{}
\vvbox{\hbox{\hcu}}
\qquad\qquad
\vvbox{\hbox{\unit\step\id}\hbox{\hcu}}
\enspace ={}\enspace{}
\vvbox{\hbox{\id}\hbox{\id}}
\enspace ={}\enspace{}
\vvbox{\hbox{\id\step\unit}\hbox{\hcu}}
\qquad\qquad
\vvbox{\hhbox{\cu\hstep\id}
       \hbox{\hstep\hcu}}
\enspace ={}\enspace{}
\vvbox{\hhbox{\id\hstep\cu}
       \hbox{\hcu}}
\end{displaymath}
{\scriptsize{ a) An algebra in a monoidal category}}
$$
\epsilon ={}\enspace{}\;
\vvbox{\hbox{\counit}}
\qquad\qquad
\Delta ={}\enspace{}
\vvbox{\hbox{\hcd}}
\qquad\qquad
\vvbox{\hbox{\hcd}
       \hbox{\counit\step\id}}
\enspace ={}\enspace{}
\vvbox{\hbox{\id}\hbox{\id}}
\enspace ={}\enspace{}
\vvbox{\hbox{\hcd}
       \hbox{\id\step\counit}}
\qquad\qquad
\vvbox{\hbox{\hstep\hcd}
       \hhbox{\cd\hstep\id}}
\enspace ={}\enspace{}
\vvbox{\hbox{\hcd}
       \hhbox{\id\hstep\cd}}
$$
{\scriptsize{ b) A coalgebra in a monoidal category}}
$$
\vvbox{\hhbox{\cu}
       \hhbox{\cd}}
\enspace ={}\enspace{}
\vvbox{\hhbox{\cd\step\cd}
       \hbox{\id\step\hx\step\id}
       \hhbox{\cu\step\cu}}
\qquad\qquad
\vvbox{\hbox{\hstep\unit}
       \hhbox{\cd}}
\enspace ={}\enspace{}
\vvbox{\hbox{\unit\step\unit}}
\qquad\qquad
\vvbox{\hhbox{\cu}
       \hbox{\hstep\counit}}
\enspace ={}\enspace{}
\vvbox{
       \hbox{\counit\step\counit}}
\qquad\qquad
\vvbox{\hbox{\unit}
       \hbox{\counit}}
\enspace =\enspace
{\rm id}_{\underline 1}
$$
{\scriptsize{c)Bialgebra axioms}}
$$
  \vvbox{\hhbox{\cd}
       \hbox{\S\step\id}
       \hhbox{\cu}}
\enspace ={}\enspace{}
\vvbox{\hbox{\counit}
       \hbox{\unit}}
\enspace ={}\enspace{}
\vvbox{\hhbox{\cd}
       \hbox{\id\step\S}
       \hhbox{\cu}}
\qquad\qquad\qquad
\vvbox{\hhbox{\cu}
       \hbox{\hstep\S}}
\enspace ={}\enspace{}
\vvbox{\hbox{\S\step\S}
       \hbox{\hx}
       \hhbox{\cu}}
\qquad\qquad
\vvbox{\hbox{\hstep\S}
              \hhbox{\cd}}
\enspace ={}\enspace{}
\vvbox{\hhbox{\cd}
       \hbox{\hx}
       \hbox{\S\step\S}}
$$
{\scriptsize{d) Antipode axiom and properties}}
$$
\mu_r:={}\enspace
\matrix{\object{X}\step\object{A}\cr
         \vvbox{\hbox{\ru}}\cr
         \object{X}\step \cr}
\qquad\qquad
\vvbox{\hbox{\id\step\unit}
       \hbox{\ru}}
\enspace ={}\enspace
\vvbox{\hbox{\id}
       \hbox{\id}}
\qquad\qquad
\vvbox{\hhbox{\ru\step\id}
       \hbox{\Ru}}
\enspace ={}\enspace
\vvbox{\hhbox{\id\hstep\cu}
       \hbox{\ru}}
\qquad\qquad
\matrix{\object{X}\step\object{Y}\hstep\step\object{A}\hstep\cr
	\vvbox{\hhbox{\id\step\id\step\cd}
	       \hbox{\id\step\hx\step\id}
	       \hbox{\ru\step\ru}}\cr
	\object{X}\Step\object{Y}\step}
$$
{\scriptsize e) Module axioms and
               module structure on tensor product of modules}
\caption{ The basic algebraic structures in a braided category }
\label{Fig-Main}
\end{figure}

\subsec{}
{\em An algebra} in a monoidal category $\cal C$ is an object $A$
equipped with unit $\eta=\eta_A:\, \underline 1\rightarrow A$
and multiplication $\mu=\mu_A:\, A\otimes A\rightarrow A$
obeying the axioms in Fig.\ref{Fig-Main}a.
{\em A coalgebra} is object $C$ equipped with
counit $\epsilon=\epsilon_A:\,C\rightarrow\underline 1$ and
comultiplication $\Delta=\Delta_A:\,A\rightarrow A\otimes A$
obeying the axioms of algebra turned upside-down
(Fig.\ref{Fig-Main}b).
Let $A$, $B$ be two algebras in a braided category $\cal C$.
Then $A\otimes B$ is also an algebra with the multiplication
\begin{equation}
\mu_{A\otimes B}:=
(\mu_A\otimes\mu_B)\circ (A \otimes\Psi_{B,A}\otimes B)
\end{equation}
Finally, \cite{M2},\cite{Majid7}
{\em A bialgebra $A$ in a braided category $\cal C$} is
an object in $\cal C$ equipped with algebra and coalgebra structures
obeying the compatibility axioms in Fig.\ref{Fig-Main}c
(i.e. $\Delta_A$, $\epsilon_A$ are algebra homomorphisms or, equivalently,
     $\mu_A$, $\eta_A$ are coalgebra homomorphisms).
{\em A Hopf algebra $A$ in a braided category} $\cal C$
({\em braided group\/} or {\em braided Hopf algebra}) is
a bialgebra in $\cal C$ with antipode $S:\,A\rightarrow A$ which is
convolution-inverse to identical map
(the first identity in Fig.\ref{Fig-Main}d).
\par
For example, in the category of sets with monoidal structure given by
Cartesian product any object becomes a coalgebra uniquely,
and a braided group in this category is an ordinary group.
Mainly, we do not exploit additive structure of the category $\cal C$,
and so we prefer Majid's term {\em braided group}.
However, all results of this paper hold
if we suppose categories and functors to be $k$-linear,
where $k$ is a commutative ring.

Now we recall basic properties of braided groups from
\cite{Majid6}-\cite{Majid10}.
Firstly, the last two identities in Fig.\ref{Fig-Main}d follow
from the axioms \cite{Majid7}.
Axioms for (right) module $X$ over an algebra $A$
with action $\mu_r :\,X\otimes A\rightarrow X$
in Fig.\ref{Fig-Main}e are obtained by "polarization" of the algebra axioms.
We use slightly modified graphical notation for action
to differ it from multiplication.
The category ${\cal C}_A$ of right modules over bialgebra $A$ in $\cal C$
is monoidal \cite{Majid7} with the module structure on the tensor product
$X\otimes Y$ of modules $X$ and $Y$ defined using comultiplication
(Fig.\ref{Fig-Main}e).

\subsec{}
\label{Subsec-Sym}
In diagrammatic language we have important symmetries \cite{Majid8} of the
axioms of braided group which extend to symmetries of theorems when
hypotheses are also symmetric.
For a diagram in a braided category $\cal C$ one can applies any of
the following symmetry transformations:
{\em input-output\/} or {\em up-side-down\/} which turns a structure into
dual or co- structure,
{\em left-right\/} and {\em mirror symmetry transformation};
the last is trivial in the unbraided case.
For tangles the first two are axial symmetries and the third is reflection.
Obtained diagram can be respectively attributed to
the opposite category ${\cal C}^{\rm op}$ with the braiding
$X\otimes Y \buildrel {\Psi _{Y,X}} \over \longleftarrow Y\otimes X$, or
to the category $\cal C _{\rm op}$ with reversed tensor product and braiding
\begin{displaymath}
X\otimes _{\rm op} Y=Y\otimes X  \buildrel{\Psi _{Y,X}}\over\longrightarrow
  X\otimes Y=Y\otimes _{\rm op} X,
\end{displaymath}
or to the category $\overline{\cal C}$ with the same tensor product and with
inverse braiding $\Psi^{-1}$.
In particular, note that the collection of bialgebra (braided group) axioms
for $A$ in a braided category $\cal C$ is input-output and left-right
symmetrical, i.e. $A$ is a bialgebra (braided group) in the categories
${\cal C}^{\rm op}$ and ${\cal C}_{\rm op}$ also.

For any algebra (resp. coalgebra) $A$ in $\cal C$ we will always consider
{\em the opposite algebra\/}
$(A^{\rm op},\mu_{A^{\rm op}}:=\mu_A\circ\Psi^{-1})$
(resp. {\em the opposite coalgebra\/}
$(A_{\rm op},\Delta_{A_{\rm op}}:=\Psi^{-1}\circ\Delta_A$)
as an object of the category $\overline{\cal C}$.
In particular, $(A^{\rm op})^{\rm op}=A$.
If $A$ is a bialgebra in $\cal C$ then $A^{\rm op}$ and $A_{\rm op}$ are
bialgebras in $\overline{\cal C}$ (cf. \cite{Majid8}).
Antipode $S^-$ for $A^{\rm op}$ (or, the same, for $A_{\rm op}$) is called
{\em skew antipode} and equals $S^{-1}$ if both $S$ and $S^-$
exist. The last two identities in Fig.\ref{Fig-Main}d mean exactly that
antipode $S_A$ is a bialgebra morphism
$(A^{\rm op})_{\rm op}\rightarrow A$
(or $A\rightarrow (A_{\rm op})^{\rm op}$) in $\cal C$.

\begin{figure}
$$\matrix{\object{X}\Step\cr
	\vvbox{\hbox{\id\step\hcoev}
	       \hbox{\hev\step\id}}\cr
	\Step\object{X}}
\enspace =\enspace
\matrix{\object{X}\cr
        \vvbox{\hbox{\id}\hbox{\id}}\cr
        \object{X}}
\qquad\quad
\matrix{\Step\object{X^\vee}\cr
	\vvbox{\hbox{\hcoev\step\id}
	       \hbox{\id\step\hev}}\cr
	\object{X^\vee}\Step}
\enspace =\enspace
\matrix{\object{X^\vee}\cr
        \vvbox{\hbox{\id}\hbox{\id}}\cr
        \object{X^\vee}}
\qquad\qquad
\vvbox{\hbox{\id\step\O{f^\vee}}
       \hbox{\hev}}
\enspace =\enspace
\vvbox{\hbox{\O{f}\step\id}
       \hbox{\hev}}
\quad\Leftrightarrow\quad
f^\vee =\enspace
\vvbox{\hbox{\hcoev\step\id}
       \hbox{\id\step\O{f}\step\id}
       \hbox{\id\step\hev}}
$$
{\scriptsize a) (Right) dual object and arrow}
$$
\rho^{*n}=\enspace
\matrix{\object{X\!\dots\! X}\step\Step\object{Y\!\dots\! Y}\cr
\vvbox{\hbox{\id\step\coro{\rho}\hhstep\hhstep\obj{\vdots}\step\step\id}
       \hbox{\coRo{\rho}}}}
\qquad\qquad
\matrix{\object{X^{\otimes m}}\Step\object{Y^{\otimes n}}\cr
\vvbox{\hbox{\O{f}\Step\id}
       \hbox{\coro{\rho^{*n}}}}}
\enspace=\enspace
\matrix{\object{X^{\otimes n}}\Step\object{Y^{\otimes m}}\cr
\vvbox{\hbox{\id\Step\O{g}}
       \hbox{\coro{\rho^{*m}}}}}
$$
{\scriptsize b) $\rho$-dual arrows}
$$
\rho\cdot\rho^\prime=\enspace
\vvbox{\hhbox{\step\cd\step\cd\step}
       \hbox{\dd\step\hx\step\d}
       \hbox{\coro{\rho}\step\coro{\rho^\prime}}}
\qquad\qquad
\rho\,\widetilde\cdot\,\rho^\prime=\enspace
\vvbox{\hhbox{\cd\Step\cd}
       \hbox{\id\step\coro{\rho^\prime}\step\id}
       \hbox{\coRo{\rho}}}
$$
{\scriptsize c) Two products of pairings}
\caption{Dual and pairings.}
\label{Fig-Pairing}
\end{figure}

\subsec{}
\label{subsec_dual}
An object $X^\vee$ is called {\em (right\/) dual} for an object $X$ in a
(braided) monoidal category ${\cal C}$ if pairing
$\cup :\,X\otimes X^\vee\rightarrow\underline 1$
and copairing $\cap :\,\underline 1\rightarrow X^\vee\otimes X$
obeying the first two identities in Fig.\ref{Fig-Pairing}a are given.
Dual arrow $f^\vee$ is defined by one of the two equivalent conditions
in Fig.\ref{Fig-Pairing}a.
In this way a braided monoidal functor
$(\_)^\vee:\,{\cal C}\rightarrow{\cal C}^{\rm op}_{\rm op}$
can be defined if $X^\vee$ exists for each $X\in{\rm Obj}({\cal C})$.
Without loss of generality by coherence theorem \cite{KL}
we shall assume that $(\_)^\vee$ is a strict monoidal functor:
$(X\otimes Y)^\vee=Y^\vee\otimes X^\vee$,
$(f\otimes g)^\vee=g^\vee\otimes f^\vee$.
Application of $(\_)^\vee$ to a bialgebra (braided group) $A$ gives the same
structure on $A^\vee$ (cf. \cite{Majid8}):
\begin{equation}
\mu_{A^\vee}:={\Delta_A}^\vee,
\qquad
\eta_{A^\vee}:={\epsilon_A}^\vee,
\qquad
\Delta_{A^\vee}:={\mu_A}^\vee,
\qquad
\epsilon_{A^\vee}:={\eta_A}^\vee
\end{equation}
{\em Left dual} ${}^\vee X$ is the right dual in the category
${\cal C}^{\rm op}$ or ${\cal C}_{\rm op}$.

More generally, let $\cal C$ be a monoidal category,
$X,Y\in{\rm Obj}({\cal C})$,
$\rho\in{\rm Hom}_{\cal C}({X\otimes Y},\underline 1)$
and morphisms
$\rho^{*n}:\,X^{\otimes n}\otimes Y^{\otimes n}\rightarrow\underline 1$
be defined by Fig.\ref{Fig-Pairing}b.
We say that arrows $f:\,X^{\otimes m}\rightarrow X^{\otimes n}$ and
$g:\,Y^{\otimes n}\rightarrow Y^{\otimes m}$ are $\rho$-{\em dual}
if the identity in Fig.\ref{Fig-Pairing}b is satisfied.
Let $A$ and $H$ be bialgebras in braided category $\cal C$.
Morphism $\rho:\,A\otimes H\rightarrow\underline 1$ is called
{\em a bialgebra pairing}
if algebra (resp. coalgebra) structure on $A$ and
coalgebra (resp. algebra) structure on $H$ are $\rho$-dual.
Convolution product $\cdot$ and 'the second' product $\widetilde\cdot$ for
$\rho,\rho^\prime\in{\rm Hom}_{\cal C}(X\otimes Y,\underline 1)$
are defined in Fig.\ref{Fig-Pairing}c.
We denote by $\rho^-$, $\rho^{\sim{}}$ corresponding inverse to $\rho$.
Let $\overline\rho:=\rho^-\circ\Psi^{-1}$.
If $A$ or $H$ has (skew) antipode then $\rho^{\sim{}}$ (resp. $\rho^-$)
exists and
\begin{equation}
\label{anti-pair}
\rho\circ(S_A\otimes H)=\rho^{\sim{}}=\rho\circ(A\otimes S_H)
\qquad\enspace
   \rho\circ(S_A^-\otimes H)=\rho^-=\rho\circ(A\otimes S_H^-)
\end{equation}
If $\rho^-$ or $\rho^{\sim{}}$ exists then $\rho$-duality between
multiplications and comultiplications implies $\rho$-duality between units and
counits.
If $(A,H,\rho)$ is bialgebra pairing in $\cal C$ then
$(A_{\rm op},H_{\rm op},\rho^-)$, $(A^{\rm op},H^{\rm op},\rho^{\sim{}})$,
$(H^{\rm op},A^{\rm op},\overline\rho)$ are bialgebra pairing in
$\overline{\cal C}$.
\par
{\em A bialgebra copairing} $\cal R$ in $\cal C$ is
a bialgebra pairing in ${\cal C}^{\rm op}$.
In this case we will say about $\cal R$-{\em codual morphisms} and
{\em bialgebra copairing} and use similar notations
${\cal R}^-$, ${\cal R}^{\sim{}}$, $\overline{\cal R}$.

\begin{figure}
$$
\mu^{\rm op}:=\enspace
\matrix{\object{A}\step\object{X}\cr
\vvbox{\hbox{\hxx}
       \hbox{\ru\hhstep\obj{\mu}}}}
\qquad\qquad
\matrix{\object{X}\step\hstep\object{A}\hstep\cr
\vvbox{\hhbox{\krl\id\step\cd}
       \hhbox{\krl\ru\hstep\obj{\mu}\hstep\id}
       \hbox{\Ru\hhstep\obj{\mu^-}}}}
\enspace=\enspace
\matrix{\object{X}\step\object{A}\cr
\vvbox{\hbox{\id\step\counit}}}
\enspace=\enspace
\matrix{\object{X}\step\hstep\object{A}\hstep\cr
\vvbox{\hhbox{\krl\id\step\cd}
       \hhbox{\krl\ru\hstep\obj{\mu^-}\hstep\id}
       \hbox{\Ru\hhstep\hhstep\obj{\mu}}}}
\qquad\qquad
\matrix{
\vvbox{\hbox{\id\hstep\r}
       \hhbox{\krl\hru\Step\id}}\cr
\object{X}\Step\hstep\object{H}}
$$
$$\hbox{\scriptsize a) (Co)module structures on $X$}$$
$$
\matrix{
\matrix{
\matrix{\object{X}\step\object{A}\cr
\vvbox{\hbox{\ru\obj{\mu_r}}}}
\enspace=\enspace
\vvbox{\hhbox{\krl\id\step\cd}
       \hbox{\hx\step\id}
       \hbox{\id\step\ru\obj{\mu_{\rm ad}}}
       \hbox{\obj{\mu_\ell}\lu}}
\cr\begin{picture}(1,1)\end{picture}\cr
\matrix{\object{X}\step\object{A}\cr
\vvbox{\hbox{\ru\obj{\mu_{\rm ad}}}}}
\enspace=\enspace
        \vvbox{\hhbox{\krl\id\step\cd}
               \hbox{\hx\step\id}
               \hbox{\S\step\ru\obj{\mu_r}}
               \hbox{\obj{\mu_\ell}\lu}} }
&\qquad\qquad&
\matrix{
\matrix{\object{X}\step\object{A}\step\object{X^\vee}\cr
\vvbox{\hhbox{\krl\id\step\lu}
       \hbox{\ev}}}
\;:=\;
\matrix{\object{X}\step\object{A}\step\object{X^\vee}\cr
\vvbox{\hhbox{\krl\ru\step\id}
       \hbox{\ev}}}
\cr
\matrix{\Step\object{X^\vee}\step\object{A}\cr
        \vvbox{\hbox{\hcoev\step\id\step\S}
               \hbox{\id\step\id\step\hx}
               \hbox{\id\step\ru\step\id}
               \hbox{\id\step\ev}}\cr
       \object{X^\vee}\Step\step}
\quad
\matrix{\Step\step\object{X^\vee}\cr
        \vvbox{\hbox{\hcoev\Step\id}
               \hbox{\id\step\rd\step\id}
               \hbox{\id\step\hx\step\id}
               \hbox{\id\step\SS\step\hev}}\cr
        \object{X^\vee}\step\object{A}\Step}    }\cr
\hbox{\scriptsize b) Adjoint action}&&
\hbox{\scriptsize c) (Co)module structures on $X^\vee$}}
$$
\caption{Constructions of new (co)modules}
\label{Fig-Mod}
\end{figure}
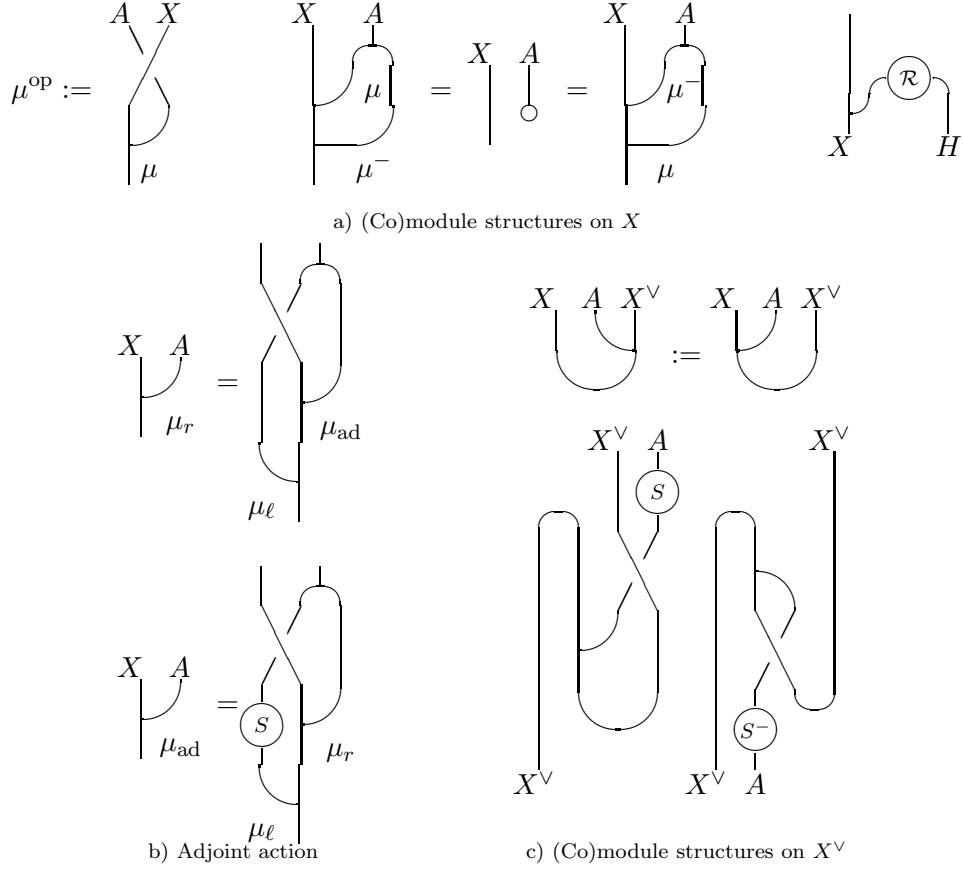

\subsec{}
\label{New-(co)mod-str}
Let $A$ be a bialgebra (braided group) in $\cal C$ and $X$ a (co)module over
$A$.
Then there exist other various (co)module structures on $X$ or $X^\vee$
over $A$ or its opposite or dual \cite{Majid8},
which can be obtained as a result of sequential application of the following
basic procedures (and their left-right, input-output and mirror reversed
forms).
One can define new (co)module structures on the underlying object of
a right $A$-module $X$ with action $\mu$ as shown in Fig.\ref{Fig-Mod}a:
\begin{itemize}
\item
left $A^{\rm op}$-module with {\em opposite\/} action $\mu^{\rm op}$
(this construction extends to isomorphism of monoidal categories
${\cal C}_A$ and ${}_{A^{\rm op}}\overline{\cal C}$);
\item
right $A^{\rm op}$-module with {\em inverse\/} action $\mu^-$
if $A$ is a braided group (and, therefore, $\mu^-$ exists and equals
$\mu\circ (X\otimes S_A)$)
(in this case one has a strict monoidal functor from ${\cal C}_A$ to
${\cal C}_{(A^{\rm op})_{\rm op}}$);
\item
right $H$-comodule via bialgebra copairing
${\cal R}:\,\underline 1\rightarrow A\otimes H$
(this correspondence extends to strict monoidal functor
${\cal C}_A\rightarrow\overline{\cal C}^{H^{\rm op}}$).
\end{itemize}
To complete the picture we note that identity functor together with
the natural transformation $\Psi$ defines equivalence of
the monoidal categories
${\cal C}_A$ and $(\overline{\cal C}_{A_{\rm op}})_{\rm op}$.
\par
Let $A$ be a bialgebra and $(X,\mu_\ell,\mu_r)$ a bimodule over $A$.
An action $\mu_{\rm ad}:\,X\otimes A\rightarrow A$ is called (right) adjoint
if the first identity in Fig.\ref{Fig-Mod}b is satisfied.
If $A$ is a braided group adjoint action always exists and is defined
by the second identity in Fig.\ref{Fig-Mod}b.
See \cite{M7} about properties of adjoint action.
Adjoint coaction is defined by input-output reversed diagrams.
\par
Let $X$ be a right module over braided group $A$.
The identity in Fig.\ref{Fig-Mod}b defines left $A$-module structure
on dual object $X^\vee$.
Then one can use procedures described above to convert
left action into right action. The last two diagrams in Fig.\ref{Fig-Mod}c
describe (co)module structure on $X^\vee$ which turns it into dual in the
category ${\cal C}_H$ (resp. ${\cal C}^H$).

\subsec{}
Quantum braided groups in a braided category were introduced in \cite{Majid7}
and basic theory was developed there.
The following are definitions from \cite{Majid7} in a slightly modified form
suitable for our use.

\begin{figure}
$$
\overline\Delta^{\rm op}\cdot{\cal R}
\quad:=\quad
\vvbox{\hbox{\cd\hstep\Obj{\overline\Delta^{\rm op}}\hstep\r}
       \hbox{\id\Step\hx\Step\id}
       \hbox{\cu\step\cu}}
\quad ={}\quad
\vvbox{\hbox{\r\step\cd}
       \hbox{\id\Step\hx\Step\id}
       \hbox{\cu\step\cu}}
\quad=:\quad
{\cal R}\cdot\Delta
$$
{\scriptsize a) Comultiplications $\Delta$ and $\overline\Delta^{\rm op}$
		 are adjoint}
$$
\matrix{
\matrix{\object{X}\step\hstep\object{A}\hstep\cr
\vvbox{\hhbox{\krl\id\step\cd\hstep\obj{\Delta}}
       \hbox{\hx\step\id}
       \hbox{\id\step\ru}}\cr
	\object{A}\step\object{X}\step}
\enspace=\enspace
\matrix{\object{X}\step\hstep\object{A}\hstep\cr
\vvbox{\hhbox{\krl\id\step\cd\hstep\obj{\overline\Delta^{\rm op}}}
       \hbox{\hxx\step\id}
       \hbox{\id\step\ru}}\cr
	\object{A}\step\object{X}\step}
&\enspace&
\Psi=\enspace
\vvbox{\hbox{\id\step\id\hstep\ro{\cal R}}
       \hhbox{\krl\id\step\hru\step\hstep\dd}
       \hbox{\hx\step\dd}
       \hbox{\id\step\ru}}
\qquad
\Psi^{-1}=\enspace
\vvbox{\hbox{\id\step\id\hstep\ro{\overline{\cal R}}}
       \hhbox{\krl\id\step\hru\step\hstep\dd}
       \hbox{\hxx\step\dd}
       \hbox{\id\step\ru}}
\cr
\hbox{\scriptsize b) The condition on modules from
                     ${\cal C}_\cO{A,\overline A}$}
&&
\hbox{\scriptsize c) Braiding in ${\cal C}_\cO{A,\overline A}$}
}
$$
$$
\matrix{\object{X}\Step\Step\step\cr
	\vvbox{\hbox{\id\step\r}
               \hbox{\hx\Step\id}
               \hbox{\id\step\Ru\r}
               \hbox{\id\step\Ru\Step\id}
               \hbox{\id\step\id\step\r\step\id}
               \hbox{\id\step\hxx\Step\hcu}
               \hbox{\hcu\step\id\Step\hstep\id}}}
\quad ={}\quad
\matrix{\object{X}\Step\Step\step\cr
	\vvbox{\hbox{\id\step\R}
               \hbox{\hx\step\r\step\id}
               \hbox{\id\step\ru\Step\hcu}
               \hbox{\id\step\id\step\r\hstep\id}
               \hbox{\id\step\hx\Step\id\hstep\id}
               \hbox{\hcu\step\Ru\hstep\id}}}
$$
{\scriptsize d) A relative form of the Yang-Baxter equation.}
\caption{The axiom for a quantum braided group $(A,{\cal R}$).
         The braided category ${\cal C}_\cO{A,\overline A}$ of modules. }
\label{Fig-QBG}
\end{figure}

{\em A quasitriangular bialgebra} in a braided category $\cal C$ is a pair of
bialgebras $A$ in $\cal C$ and $\overline A$ in $\overline{\cal C}$ with the
same underlying algebra
($\Delta$ and $\overline\Delta$ are comultiplications
 in $A$ and $\overline A$ respectively),
and convolution invertible bialgebra copairing
({\em quasitriangular structure})
${\cal R}:\,\underline 1\rightarrow\overline A_{\rm op}\otimes A$,
satisfying the condition in Fig.\ref{Fig-QBG}a.
(It follows directly from the definition that counits for $A$ and
for $\overline A$ are the same.)
{\em A quantum braided group} or {\em a quasitriangular Hopf algebra} in
$\cal C$ is a quasitriangular bialgebra such that $A$ and $\overline A$
have antipodes $S$ and $\overline S$ respectively.
(In this case ${\cal R}^-=(\overline S\otimes A)\circ{\cal R}$ and
 ${\cal R}^{\sim{}}=(A\otimes S)\circ{\cal R}$.)
\par
In particular,  for any bialgebra (braided group) $A$ the pair
$(A,A_{\rm op})$  is a quasitriangular bialgebra (quantum braided group)
with the trivial quasitriangular structure ${\cal R}=\eta\otimes\eta$.

Category ${\cal C}_\cO{A,\overline A}$ is
a full subcategory of ${\cal C}_A$ with objects $X$ satisfying
the identity in Fig.\ref{Fig-QBG}b.
${\cal C}_\cO{A,\overline A}$ is a monoidal subcategory of ${\cal C}_A$
and braided with $\Psi$ and $\Psi^{-1}$ shown in Fig.\ref{Fig-QBG}c
(see.  \cite{Majid7} where corresponding categories of left modules are
introduced and studied).
We use a brief notation ${\cal C}_\cO{A}$ for
${\cal C}_\cO{A,A_{\rm op}}$.

As Majid showed, the basic formulas for ordinary quantum groups
\cite{Drinfel'd2} have analogues in this more general context but some of
them exist only in a relative form: actions on arbitrary modules from
${\cal C}_{\cO{A,\overline A}}$ take part in them.
One can obtain the standard formulas for ordinary quantum groups if
one considers the action on the unit element of the regular module.
As an example a relative form of the Yang-Baxter equation for
$X\in{\rm Obj}({\cal C}_{\cO{A,\overline A}})$ is shown on
Fig.\ref{Fig-QBG}d.

%
                    \sect{Crossed modules.}
%

\subsec{}
\label{cross-intro}
Here we introduce categories of crossed modules over
a bialgebra (braided group) $A$ in a braided category $\cal C$.

\begin{definition}
{\em A right} (resp. {\em left-right}) {\em crossed module} over
a bialgebra $A$ in a braided category $\cal C$ is an object $X$ with right
(resp. left) $A$-module and right $A$-comodule structures
obeying the first (resp. the second) identity in Fig.\ref{Fig-DY}a.
$\DY{\cal C}_A^A$ (resp. ${}_A\DY{\cal C}^A$) is the
{\em category of right} (resp. {\em left-right}) {\em crossed modules}
with morphisms which are both module and comodule maps.
Objects of categories ${}^A_A\DY{\cal C}$ and ${}^A\DY{\cal C}_A$ are
described by left-right reversed axioms.
\par
Note, that the category $\DY{\cal C}_A^A$ is defined when $\cal C$ is only
pre-braided, whereas definition of ${}_A\DY{\cal C}^A$ uses both $\Psi$ and
$\Psi^{-1}$.
\end{definition}

Let $H$ be also a bialgebra in $\cal C$ and
$\rho : A\otimes H\rightarrow\underline 1$ a bialgebra pairing.

\begin{definition}
Objects of category ${}_A\DY{\cal C}_H$ (resp. $\DY{\cal C}_{A,H}$) are
both $A$- and $H$- modules satisfying the first (resp. the second) identity
in Fig.\ref{Fig-DY}b.
\par
If the pairing is trivial ($\rho=\epsilon_A\otimes\epsilon_H$) then objects of
${}_A\DY{\cal C}_H$ are $A$-$H$-bimodules.
\end{definition}

\begin{figure}
$$
\matrix{\object{X}\hstep\step\object{A}\hstep\cr
        \vvbox{\hhbox{\krl\id\step\cd}
               \hbox{\hx\step\id}
               \hbox{\id\step\k}
               \hbox{\hx\step\id}
               \hhbox{\krl\id\step\cu}}\cr
       \object{X}\hstep\step\object{A}\hstep}
\enspace ={}\enspace
\matrix{\object{X}\Step\hstep\object{A}\hstep\cr
        \vvbox{\hbox{\rd\step\hcd}
               \hbox{\id\step\hx\step\id}
	       \hbox{\ru\step\hcu}}\cr
       \object{X}\Step\hstep\object{A}\hstep}
\qquad\qquad
\matrix{\hstep\object{A}\step\object{X}\hstep\cr
\vvbox{\hhbox{\krl\cd\hstep\id}
       \hhbox{\krl\d\hstep\hlu}
       \hbox{\hstep\hxx}
       \hhbox{\krl\hstep\hrd\hstep\d}
       \hhbox{\krl\hstep\id\hstep\cu}}\cr
\hstep\object{X}\step\object{A}\hstep}
\enspace ={}\enspace
\matrix{\hstep\object{A}\step\hstep\object{X}\step\cr
\vvbox{\hbox{\hcd\step\rd}
       \hbox{\id\step\hx\step\id}
       \hbox{\lu\step\hcu}}\cr
\step\object{X}\step\hstep\object{A}\hstep}
$$
{\scriptsize a) axioms for objects of
              $\DY{\cal C}_A^A$ and ${}_A\DY{\cal C}^A$}
$$
\matrix{\hstep\object{A}\hstep\step\object{X}\hstep\step\object{H}\hstep\cr
\vvbox{\hhbox{\krl\cd\step\id\step\hstep\id}
       \hbox{\id\step\lu\step\hcd}
       \hbox{\id\Step\hx\step\id}
       \hbox{\coro{\rho}\step\ru}}}
\enspace ={}\enspace
\matrix{\hstep\object{A}\hstep\step\object{X}\hstep\step\object{H}\hstep\cr
\vvbox{\hhbox{\krl\hstep\id\step\hstep\id\step\cd}
       \hbox{\hcd\step\ru\step\id}
       \hbox{\id\step\hx\Step\id}
       \hbox{\lu\step\coro{\rho}}}}
\qquad\quad
\matrix{\object{X}\step\hstep\object{A}\Step\step\object{H}\hstep\cr
\vvbox{\hhbox{\krl\id\step\cd\Step\cd}
       \hbox{\d\d\coro{\rho}\dd}
       \hbox{\step\id\step\xx}
       \hbox{\step\ru\step\dd}
       \hbox{\step\Ru}}}
\enspace ={}\enspace
\matrix{\object{X}\step\hstep\object{A}\Step\object{H}\cr
\vvbox{\hhbox{\krl\id\step\cd\step\cd}
       \hbox{\id\step\hxx\step\hxx}
       \hbox{\ru\dd\step\id\step\id}
       \hbox{\hx\Step\id\step\id}
       \hbox{\hx\Step\id\step\id}
       \hbox{\d\coro{\rho}\dd}
       \hbox{\step\Ru}}}
$$
{\scriptsize b) axioms for objects of
        ${}_A\DY{\cal C}_H$ and $\DY{\cal C}_{A,H}$}
\caption{Drinfel'd-Yetter compatibility conditions.}
\label{Fig-DY}
\end{figure}

Connections between the categories
${}_A^A\YD{}$, ${}_A\YD{}^A$, ${}^A\YD{}_A$, $\YD{}^A_A$
of crossed modules over an ordinary bialgebra (Hopf algebra)
were studied in \cite{RadTow}.
Categories $\DY{\cal C}_A^A$ and ${}_A^A\DY{\cal C}$ in
the fully braided setting were introduced and studied in
the preprint version of the present paper (May 1994).
Left-right crossed modules over a braided group were defined independently in
\cite{CY} where an interesting connection of these objects with topological
quantum field theories was discovered.
The category $\DY{\cal C}_{A,H}$ turns into the category of modules over
the corresponding Drinfel'd's double in the case of ordinary bialgebras
$A$ and $H$.

\subsec{}
We will use abbreviated form
${\rm L}_{\DY{\cal C}_A^A}^X$ (resp. ${\rm R}_{\DY{\cal C}_A^A}^X$)
for the left hand side (resp. for the right hand side) of right crossed module
axiom for $X$ over a bialgebra $A$ in $\cal C$ and similar notations for
other variants of crossed modules.

\begin{lemma}
\label{Lemma-R-LR}
Let $A$ be a bialgebra in $\cal C$ and $(X,\mu_r,\Delta_r)$ right crossed
module over $A$.
Then $X$ with opposite action $\mu_r\circ\Psi^{-1}$ is left-right crossed
module over $A^{\rm op}$.
\par
This construction defines categorical isomorphism
$\DY{\cal C}^A_A\simeq{}_{A^{\rm op}}\DY{\overline{\cal C}}^{A^{\rm op}}$.
\end{lemma}

\begin{proof}
$$
{\rm L}_{{}_{A^{\rm op}}\DY{\cal C}^{A^{\rm op}}}^X=
\Psi^{-1}_{X,A}\circ{\rm L}_{\DY{\cal C}^A_A}^X=
\Psi^{-1}_{X,A}\circ{\rm R}_{\DY{\cal C}^A_A}^X=
{\rm R}_{{}_{A^{\rm op}}\DY{\cal C}^{A^{\rm op}}}^X
$$
\end{proof}

\subsec{}
Let us describe a monoidal structure on the categories of crossed modules.

\begin{lemma}
\label{Lemma-R-times}
If $X$ and $Y$ are right crossed modules over a bialgebra $A$ then
$X\otimes Y$ also is, with module (resp. comodule) structure the braided
tensor product one from $X$ and $Y$ defined by
the last diagram in Fig.\ref{Fig-Main}e.
(resp. by the input-output reversed diagram).
This turns $\DY{\cal C}^A_A$ into a monoidal category.
\end{lemma}

\begin{proof}
See Fig.\ref{Proof-DY-otimes}.
\end{proof}

The following lemma can be verified directly or obtained as a corollary of
lemmas \ref{Lemma-R-LR} and \ref{Lemma-R-times}.

\begin{lemma}
If $X$ and $Y$ are left-right crossed modules over a bialgebra $A$ then
$X\otimes Y$ also is, with underlying module (resp. comodule) the tensor
product one in ${}_{A^{\rm op}}\overline{\cal C}$ (resp. in ${\cal C}^A$).
\end{lemma}

We will use modified notation ${}_{A_{\rm op}}\DY{\cal C}^A$ for the category
of left-right crossed modules with monoidal structure described in previous
lemma. An obvious symmetry implies existence of the second monoidal structure
${}_A\DY{\cal C}^{A^{\rm op}}$ on this category.

One can also verify that the categories ${}_A\DY{\cal C}_H$,
$\DY{\cal C}_{A,H^{\rm op}}$ are monoidal
(where 'op' emphasize that in the second category $H$-module structure on
$X\otimes Y$ is defined by tensor product
in $\overline{\cal C}_{H^{\rm op}}$).

\subsec{}
(Large) category $\bf MonCat$ of all (small) monoidal categories is
monoidal with usual Cartesian product of categories and functors.
Let $\Psi$ be braiding on monoidal category $\C$.
Then tensor product functor $\otimes: \C\times\C\rightarrow\C$ together
with a natural transformation
\begin{displaymath}
X\otimes\Psi_{X^\prime,Y}\otimes Y^\prime:\;
(X\otimes X^\prime)\otimes (Y\otimes Y^\prime)\rightarrow
(X\otimes Y)\otimes (X^\prime\otimes Y^\prime)
\end{displaymath}
define a monoidal functor.

\hfill
More generally, let
\hfill
$\C,\,\C_1,\,\C_2$
\hfill
be monoidal categories and
\hfill
$F_1:\C_1\rightarrow\C\,,$
\hfill\hfill
\newline
$F_2:\C_2\rightarrow\C\,$  monoidal functors
(for simplicity, we suppose that $\C,\,\C_1,\,\C_2$ are strict monoidal
categories and $F_1,F_2$ are strict monoidal functors).
And let $\Psi=\{\Psi_{X,Y}\vert{X\in{\rm Obj}(\C_1),Y\in{\rm Obj}(\C_2)}\}$
be a natural transformation of functors
\begin{displaymath}
\otimes\circ(F_1\times F_2)\buildrel\Psi\over\Longrightarrow
\otimes_{\rm op}\circ(F_1\times F_2):\;
\C_1\times\C_2\rightarrow\C\,,
\end{displaymath}
where $\otimes$ (resp. $\otimes_{\rm op}$) is tensor product (resp.
opposite tensor product) in $\C$.
We consider a pair of functor $F:=\otimes\circ(F_2\times F_1)$ and
natural transformation
$\lambda:F(\_)\otimes F(\_)\rightarrow F(\_\otimes\_)$,
where
$
\lambda_{Y\times X,Y^\prime\times X^\prime}:=
\id_{F_2(Y)}\otimes\Psi_{X,Y^\prime}\otimes\id_{F_1(X^\prime)}$.

\begin{lemma}
$(F,\lambda):\C_1\times\C_2\rightarrow\C$ is a monoidal functor iff
$\Psi$ satisfies the following 'hexagon identities'
\begin{eqnarray}
\Psi_{X\otimes X^\prime,Y}&=&
(\Psi_{X,Y}\otimes\id_{F_1(X^\prime)})\circ
(\id_{F_1(X)}\otimes\Psi_{X^\prime,Y})\,,\nonumber\\
\label{equ-hex}
\Psi_{X,Y\otimes Y^\prime}&=&
(\id_{F_2(Y)}\otimes\Psi_{X,Y^\prime})\circ
(\Psi_{X,Y}\otimes\id_{F_2(Y^\prime)})\,,
\end{eqnarray}
for $X,X^\prime\in{\rm Obj}(\C_1)$
and $Y,Y^\prime\in{\rm Obj}(\C_2)$.
\par
In this case we will say that $\Psi$ is {\em a generalized braiding} for
$\C_1\buildrel{F_1}\over\rightarrow\C\buildrel{F_2}\over\leftarrow\C_2$.
\end{lemma}

\begin{figure}
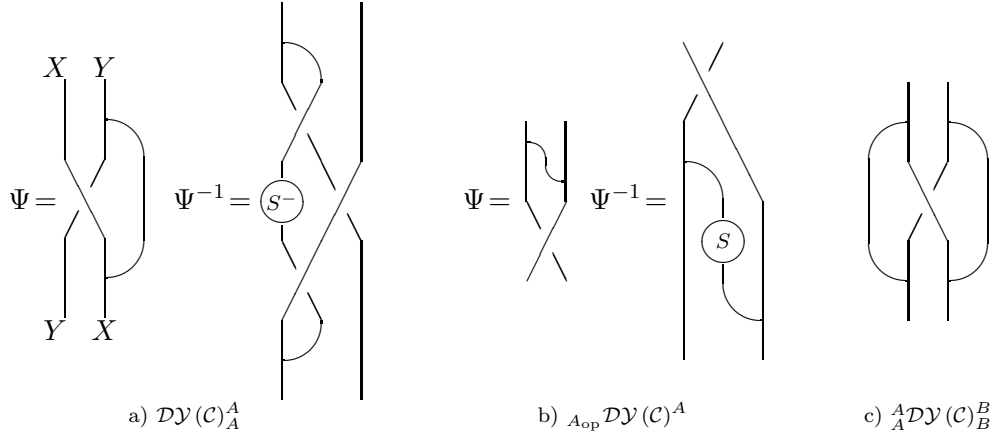

$$
\matrix{
\Psi\!=\!\!\!
\matrix{\object{X}\step\object{Y}\step\cr
\vvbox{\hbox{\id\step\rd}
       \hbox{\hx\step\id}
       \hbox{\id\step\ru}}\cr
\object{Y}\step\object{X}\step}
\enspace
\Psi^{-1}\!=\enspace\;
\vvbox{\hbox{\rd\step\id}
       \hbox{\hxx\step\id}
       \hbox{\SS\step\hxx}
       \hbox{\hxx\step\id}
       \hbox{\ru\step\id}}
&\quad\enspace&
\Psi\!=\;
\vvbox{\hhbox{\krl\hrd\hstep\id}
       \hhbox{\krl\id\hstep\hlu}
       \hbox{\hxx}}
\enspace\;
\Psi^{-1}\!=\;
\vvbox{\hbox{\hx}
       \hbox{\rd\d}
       \hbox{\id\step\S\step\id}
       \hbox{\id\step\lu}}
&\quad\enspace&
\vvbox{\hbox{\ld\step\rd}
       \hbox{\id\step\hx\step\id}
       \hbox{\lu\step\ru}}
\cr \hbox{\scriptsize a) $\DY{\cal C}_A^A$} &&
 \hbox{\scriptsize b) ${}_{A_{\rm op}}\DY{\cal C}^A$}  &&
 \hbox{\scriptsize c) ${}_A^A\DY{\cal C}^B_B$}}
$$
\caption{Braidings in the categories of crossed modules}
\label{Fig-DY-Psi}
\end{figure}

The following is an example a generalized braiding.
Let $A$ be a bialgebra in $\cal C$.
For right $A$-module $(X,\mu^X_r)$ and right $A$-comodule $(Y,\Delta^Y_r)$
we define morphism
\begin{equation}
\Psi^A_{X,Y}:=(Y\otimes\mu^X_r)\circ(\Psi_{X,Y}\otimes A)\circ
(X\otimes\Delta^Y_r)
\end{equation}
or, equivalently, by the first diagram in Fig.\ref{Fig-DY-Psi}.

\begin{lemma}
The collection $\Psi^A:=\{\Psi^A_{X,Y}\}$ is a generalized braiding for
pair of underlying functors
$\C_A\rightarrow\C\leftarrow\C^A$.
\end{lemma}

\begin{figure}
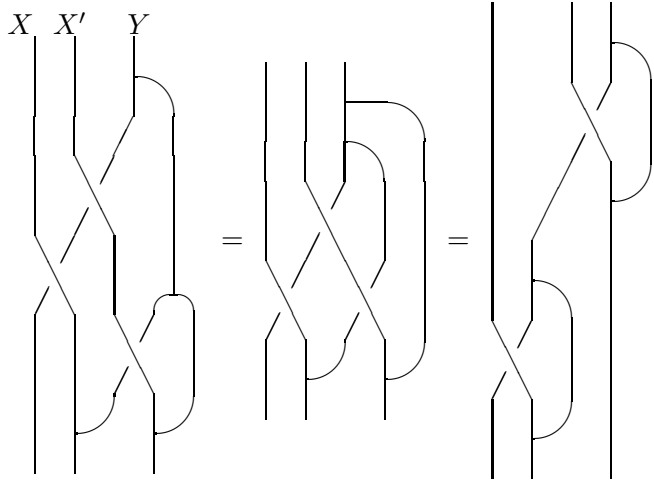

\begin{displaymath}
\matrix{\object{X}\step\object{X^\prime}\step\hstep\object{Y}\step\hstep\cr
        \vvbox{\hbox{\id\step\id\step\hstep\rd}
               \hhbox{\krl\id\step\id\step\dd\step\id}
               \hbox{\id\step\hx\step\hstep\id}
               \hbox{\hx\step\id\step\hcd}
               \hbox{\id\step\id\step\hx\step\id}
               \hbox{\id\step\ru\step\ru}}}
=\enspace
\vvbox{\hbox{\id\step\id\step\Rd}
       \hhbox{\krl\id\step\id\step\rd\step\id}
       \hbox{\id\step\hx\step\id\step\id}
       \hbox{\hx\step\hx\step\id}
       \hbox{\id\step\ru\step\ru}}
\enspace =\enspace
\vvbox{\hbox{\id\Step\id\step\rd}
       \hbox{\id\Step\hx\step\id}
       \hbox{\id\step\dd\step\ru}
       \hbox{\id\step\rd\step\id}
       \hbox{\hx\step\id\step\id}
       \hbox{\id\step\ru\step\id}}
\end{displaymath}
\caption{Proof of the hexagon identities for $\Psi^A$}
\label{Proof-Hex}
\end{figure}

\begin{proof}
Fig.\ref{Proof-Hex} proves the first identity from (\ref{equ-hex}).
A proof of the second one uses input-output reversed diagrams.
\end{proof}

\begin{theorem}
\hfill
Let $A$ be bialgebra in $\cal C$.
Then the categories
\hfill
$\DY{\cal C}^A_A$
\hfill
and
\newline
${}_{A_{\rm op}}\DY{\cal C}^A$
are pre-braided with $\Psi$ shown in Fig.\ref{Fig-DY-Psi}.
If $A$ (resp. $A^{\rm op}$) is a braided group
then ${}_{A_{\rm op}}\DY{\cal C}^A$ (resp. $\DY{\cal C}^A_A$) is braided.
The lemma \ref{Lemma-R-LR} defines isomorphism of (pre-)braided categories
$\DY{\cal C}^A_A$ and
${}_{(A^{\rm op})_{\rm op}}\DY{\overline{\cal C}}^{A_{\rm op}}$.
\end{theorem}

\begin{figure}
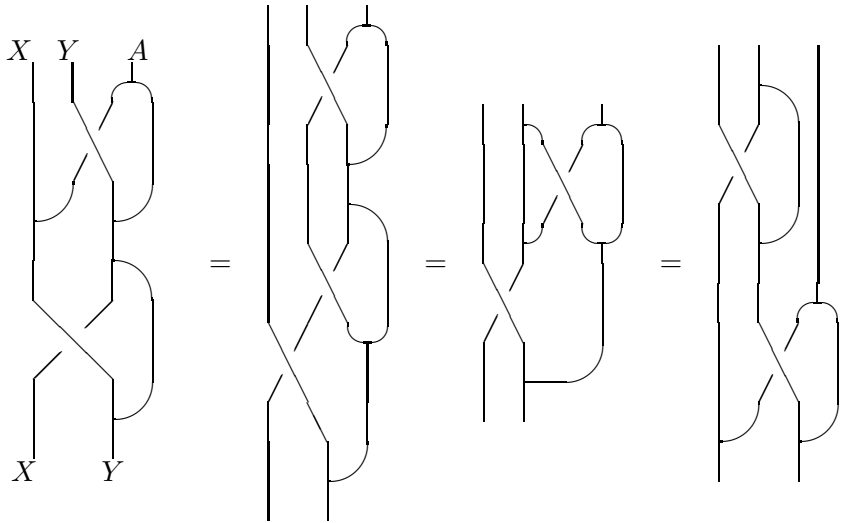

\begin{displaymath}
\matrix{\object{X}\step\object{Y}\step\hstep\object{A}\hstep\cr
\vvbox{\hhbox{\krl\id\step\id\step\cd}
       \hbox{\id\step\hx\step\id}
       \hbox{\ru\step\ru}
       \hhbox{\krl\id\Step\rd}
       \hbox{\x\step\id}
       \hbox{\id\Step\ru}}\cr
	\object{X}\Step\object{Y}\step}
\quad ={}\quad
\vvbox{\hhbox{\krl\id\step\id\step\cd}
       \hbox{\id\step\hx\step\id}
       \hhbox{\krl\id\step\id\step\ru}
       \hbox{\id\step\id\step\rd}
       \hbox{\id\step\hx\step\id}
       \hbox{\hx\step\hcu}
       \hhbox{\krl\id\step\d\step\id}
       \hbox{\id\step\hstep\ru}}
\quad ={}\quad
\vvbox{\hhbox{\krl\id\step\hrd\step\cd}
       \hbox{\id\step\id\hstep\hx\step\id}
       \hhbox{\krl\id\step\hru\step\cu}
       \hbox{\hx\Step\id}
       \hbox{\id\step\Ru}}
\quad ={}\quad
\vvbox{\hbox{\id\step\rd\hstep\id}
       \hbox{\hx\step\id\hstep\id}
       \hbox{\id\step\ru\hstep\id}
       \hhbox{\krl\id\step\id\step\cd}
       \hbox{\id\step\hx\step\id}
       \hbox{\ru\step\ru}}
\end{displaymath}
\caption{Braiding in $\DY{\cal C}_A^A$ is a module map:
\protect\newline
The first, the second and the third equalities are the module,
crossed module and comodule axioms respectively.}
\label{Proof-Mod-Map}
\end{figure}

\begin{proof}
Hexagon identities for $\DY{\cal C}^A_A$ are proven in the previous lemma.
Diagrams in Fig.\ref{Proof-Mod-Map} (resp. the input-output reversed diagrams)
show that $\Psi$ in Fig.\ref{Fig-DY-Psi}a is module (resp. comodule) map.
It is obvious that $\Psi$ and $\Psi^{-1}$ in Fig.\ref{Fig-DY-Psi}b are
inverse each to other.
\end{proof}

(Pre-)braided structures on
${}_A\DY{\cal C}^{A^{\rm op}}$, ${}_A^A\DY{\cal C}$,
${}^A\DY{\cal C}_{A_{\rm op}}$, ${}^{A^{\rm op}}\DY{\cal C}_A$
are defined in similar way using the symmetries described in
\ref{Subsec-Sym}.

\subsec{}
\label{Subsec-New-Cross}
Let $A$ be a bialgebra (braided group) in $\cal C$.
The lemma \ref{Lemma-R-LR} and the following lemmas \ref{Lemma-Inverse},
\ref{Lemma-Dual}, \ref{Lemma-Pairing} describe analogues
of constructions from \ref{New-(co)mod-str} which produce new crossed modules
from given one.
As a corollary we obtain isomorphisms between categories of crossed modules.

\begin{lemma}
\label{Lemma-Inverse}
Let $A$ be a braided group in $\cal C$ and $(X,\mu_\ell,\Delta_r)$ left-right
crossed module over $A$.
Then $X$ with inverse action $\mu_\ell^-:=\mu_\ell\circ(S\otimes X)$ is
a left-right crossed module over $A^{\rm op}$.
\par
This construction extends to braided functor
${}_A\DY{\cal C}^{A^{\rm op}}\rightarrow
 {}_{(A^{\rm op})_{\rm op}}\DY{\overline{\cal C}}^{A^{\rm op}}$
(isomorphism of braided categories if antipode $S$ is invertible).
\end{lemma}

\begin{proof}
Crossed module axiom for $X^-:=(X,\mu_\ell^-,\Delta_r)$ is verified on
Fig.\ref{Proof-new-DY}a.
\end{proof}

\begin{lemma}
\label{Lemma-Dual}
Let $A$ be a bialgebra in $\cal C$ and $X\in{\rm Obj}({}^A\DY{\cal C}_A)$,
an object $X^\vee$ dual in $\cal C$ exists and equipped with left
$A$-module and right $A$-comodule structures defined by the identity
in Fig.\ref{Fig-Mod}c and by the input-output reversed identity.
Then $X^\vee\in{\rm Obj}({}_A\DY{\cal C}^A)$.
\par
If $\cal C$ has (right) duals this construction defines a braided functor
from ${}^{A^{\rm op}}\DY{\cal C}_A$
to $({}_{A_{\rm op}}\DY{\cal C}^A)_{\rm op}^{\rm op}$.
\end{lemma}

\begin{proof}
Crossed module compatibility condition follows from identities on
Fig.\ref{Proof-new-DY}b.
\end{proof}

\begin{corollary}
Right dual for object $X$ in the category $\DY{\cal C}_A^A$ of
crossed modules
over braided group $A$ with invertible antipode is right dual $X^\vee$ for
underlying object in $\cal C$ (if this latter exists) with the action and
coaction defined in Fig.\ref{Fig-Mod}c.
\end{corollary}

\begin{proof}
Right crossed module structure on $X^\vee$ is the result of sequential
application of the constructions from lemma \ref{Lemma-Dual},
lemmas \ref{Lemma-R-LR}, \ref{Lemma-Inverse} and their input-output reversed
forms.
Module (resp. comodule) structure on $X^\vee$ is the same as in the category
${\cal C}_A$ (resp. ${\cal C}^A$).
Hence, pairing and copairing are module and comodule maps.
\end{proof}

\begin{lemma}
\label{Lemma-Two-LR}
Let $A$ be bialgebra in $\cal C$. Then
${}_A\DY{\cal C}^{A^{\rm op}}
 \buildrel{({\rm Id},\Psi)}\over\longrightarrow
 ({}_{A_{\rm op}}\DY{\cal C}^A)_{\rm op}$
is isomorphism of pre-braided categories.
\end{lemma}

Taking into attention lemmas \ref{Lemma-R-LR}, \ref{Lemma-Inverse},
\ref{Lemma-Two-LR} and the fact that the functor
${\cal C} \buildrel{({\rm Id},\Psi )}\over\longrightarrow
 {\cal C}_{\rm op}$
defines isomorphism of braided categories, we obtain the following
result about isomorphisms between categories of crossed modules as in
the case of ordinary Hopf algebras \cite{RadTow}.

\begin{corollary}
Let $A$ be a braided group in $\cal C$ with invertible antipode.
Then the following braided categories of crossed modules are isomorphic:
$$\DY{\cal C}^A_A\simeq
{}_{A_{\rm op}}\DY{\cal C}^A\simeq{}_A\DY{\cal C}^{A^{\rm op}}\simeq
{}_A^A\DY{\cal C}\simeq
{}^A\DY{\cal C}_{A_{\rm op}}\simeq{}^{A^{\rm op}}\DY{\cal C}_A.$$
\end{corollary}

In particular, there exist two braided functors from $\DY{\cal C}^A_A$
to ${}_A^A\DY{\cal C}$ defined on object in the following way:
\begin{equation}
\matrix{
(X,\mu_r,\Delta_r)\mapsto
(X,\,\mu_r\circ\Psi^{-1}\circ (S^-\otimes X),\,
                   (S\otimes X)\circ\Psi\circ\Delta_r),
\cr
(X,\mu_r,\Delta_r)\mapsto
(X,\,\mu_r\circ\Psi\circ (S\otimes X),\,
                   (S^-\otimes X)\circ\Psi^{-1}\circ\Delta_r).}
\end{equation}
In \ref{Subsec-Antipode} we will describe isomorphism between these functors.

\begin{lemma}
\label{Lemma-Pairing}
Let, $A$ and $H$ be bialgebras in $\cal C$ and
$\rho:A\otimes H\rightarrow\underline 1$ a bialgebra pairing.
Then any left-right crossed module over $A$ becomes an object of
${}_A\DY{\cal C}_H$ with $H$-module structure defined via $\rho$
(by input-output reversed form of the last diagram from Fig.\ref{Fig-Mod}a).
\par
This construction defines a monoidal functor
from ${}_A\DY{\cal C}^{A^{\rm op}}$ to ${}_A\DY{\cal C}_H$.
\end{lemma}

\proof
$$
{\rm  L}^X_{{}_A\DY{\cal C}_H}=
 (X\otimes\rho)\circ({\rm L}^X_{{}_A\DY{\cal C}^A}\otimes H)=
 (X\otimes\rho)\circ({\rm R}^X_{{}_A\DY{\cal C}^A}\otimes H)=
{\rm  R}^X_{{}_A\DY{\cal C}_H}.\qquad\emptybox
$$

The image of the functor from the previous lemma is a pre-braided subcategory
in ${}_A\DY{\cal C}_H$.
In similar way one can define a monoidal functor
from $\DY{\cal C}^A_A$ to $\DY{\cal C}_{A,H^{\rm op}}$.

\begin{corollary}
Let $A$ be bialgebra in $\cal C$ and $A^\vee$ exists.
Then the pre-braided categories ${}_A\DY{\cal C}^{A^{\rm op}}$ and
${}^{(A^\vee)^{\rm op}}\DY{\cal C}_{A^\vee}$ are isomorphic.
\end{corollary}

Finally, we note the following obvious isomorphisms of (pre-)braided
categories:
\begin{equation}
\DY{{\cal C}^{\rm op}}^A_A\simeq (\DY{\cal C}_A^A)^{\rm op},\qquad\qquad
\DY{{\cal C}_{\rm op}}^A_A\simeq ({}^A_A\DY{\cal C})_{\rm op}.
\end{equation}

\subsec{}
Braided category ${}_A^A\YD{}$ of crossed modules over an ordinary
Hopf algebra
$A$ can be obtained directly from monoidal category ${}_A{\cal M}$ as
a 'center' or 'inner double'. This construction is due to Drinfel'd
(unpublished) and a formal proof is in \cite{Majid9,M5}.
Here we show that the same with slight modification
is true in the braided case.

{\em A center}  $\cal Z(V)$ of monoidal category $\cal V$ is
a special case of {\em Pontryagin dual monoidal category} \cite{Majid9,M5}.
Objects of $\cal Z(V)$ are pairs $(V,\lambda_V)$ where $V$ is an object of
$\cal C$ and $\lambda_V$ is a natural isomorphism in
$\hbox{Nat}(V\otimes\hbox{id},\hbox{id}\otimes V)$ such that
\begin{displaymath}
\lambda_{V,\underline 1}={\rm id}\qquad
  (\hbox{id}\otimes\lambda_{V,Z})(\lambda_{V,W}\otimes\hbox{id})=
  \lambda_{V,W\otimes Z}
\end{displaymath}
and morphisms are $\phi :\; V\rightarrow W$ such that the objects are
intertwined in the form
\begin{displaymath}
(\hbox{id}\otimes\phi )\lambda_{V,Z}=
  \lambda_{W,Z}(\phi\otimes\hbox{id})\,,\qquad
  \forall Z \in{\rm Obj}({\cal C})
\end{displaymath}
The tensor product and braiding are the following:
\begin{displaymath}
 (V,\lambda_V)\otimes (W,\lambda_W)=
  (V\otimes W,\lambda_{V\otimes W})\qquad
  \lambda_{V\otimes W,Z}=
  (\lambda_{V,Z}\otimes\hbox{id})(\hbox{id}\otimes\lambda_{W,Z})\,,
\end{displaymath}
\begin{displaymath}
\Psi_{(V,\lambda_V),(W,\lambda_W)}=\lambda_{V,W}
\end{displaymath}
Existence of braiding in this case was pointed out by Drinfel'd.

Let $A$ be a braided group with invertible antipode in $\cal C$.
For a special case ${\cal V}={}_A\C$ we denote by
${\cal Z}_{\cal C}({}_A{\cal C})$ the full subcategory of
${\cal Z}({}_A{\cal C})$
with the following condition on $\lambda_V$:
for any object $W$ in $\cal C$ with trivial action (through counit)
$\lambda_{V,W}$ coincides with braiding $\Psi_{V,W}$ in $\cal C$.
(In the tensor category $\bf Vect$ of vector spaces the unit object
$\underline 1=k$ is a generator of the category and then
${\cal Z}_{\rm Vect}({}_A{\cal M})={\cal Z}({}_A{\cal M})$.)
The following proposition is analog of corresponding result from
\cite{Majid9,M5}.

\begin{proposition}
Braided monoidal categories
${}_A^A\DY{\cal C}$ and ${\cal Z}_{\cal C}({}_A{\cal C})$ are isomorphic.
\end{proposition}

\begin{proof}
Proof is like in the corresponding proposition from \cite{Majid9,M5}.
One can identify a crossed module $X$ with a pair of the underlying module
of $X$ and the braiding ${}^{({}^A_A\DY{\cal C})}\Psi_{X,\_}$.
Conversely, coaction on $X$ can be reconstructed from $(X,\lambda_X)$ as
the composition morphism
$X\buildrel{{\rm id}\otimes\eta}\over\longrightarrow X\otimes A
  \buildrel{\lambda_X}\over\longrightarrow A\otimes X$.
(We consider $A$ with left regular module structure.)
\end{proof}

Similarly, one can identify $\DY{\cal C}^A_A$ and
${\cal Z}_{\cal C}({\cal C}^A)$.
So 'centers' of categories of modules and comodules are the same and
coincide
with the category of crossed modules just as in the standard case over
$\bf Vect$.

\subsec{}
In \ref{cross-intro} we defined the category $\DY{\cal C}_{A,H}$
depending on bialgebra
pairing $\rho:\, A\otimes H\rightarrow\underline 1$, which is a fully
braided analog of the category of modules ${\cal M}_{{\cal D}(A,H,\rho)}$
over Drinfel'd double ${\cal D}(A,H,\rho)$ of ordinary Hopf algebras.
Here we want to discuss when this category can be realized as a category of
modules over something.

\begin{figure}
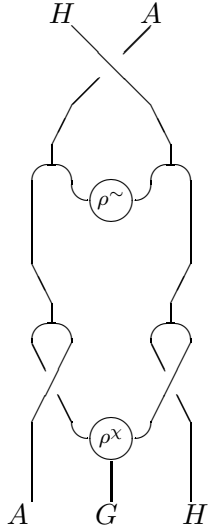

$$
\matrix{\step\object{H}\Step\object{A}\step\cr
\vvbox{\hbox{\step\x}
       \hbox{\hddcd\Step\hdcd}
       \hbox{\id\step\coro{\rho^\sim}\step\id}
       \hbox{\hdcd\Step\hddcd}
       \hbox{\hxx\Step\hxx}
       \hbox{\id\step\tu{\rho^\chi}\step\id}}\cr
\object{A}\Step\object{G}\Step\object{H}}
$$
\caption{Multiplication in $A\otimes G\otimes H$.}
\label{Fig-Prod}
\end{figure}

We will consider a special case. Let $k$ be a field, $(\Gamma,\cdot)$
an Abelian group, $k\Gamma$ a group algebra
equipped with a coquasitriangular structure
given by a function $\chi:\Gamma\times\Gamma\to k$ obeying the bicharacter
conditions \cite{Majid10}.
And let $\C$ be the category ${\cal M}^{k\Gamma}$  of
right $k\Gamma$-comodules $(X,\Delta^X_r)$ (or, equivalently, the category
of $\Gamma$-graded spaces $X=\oplus X_{\alpha\in\Gamma}$ where
$X_\alpha:=\{x\in X\vert\Delta^X_r(x)=x\otimes\alpha\}$) with braiding
defined by this coquasitriangular structure:
$$
\Psi: x\otimes y\mapsto\chi(\alpha,\beta)\cdot y\otimes x\,,\quad
 x\in X_\alpha\,,\enspace y\in Y_\beta\,.$$
Let $A=\oplus_{\alpha\in\Gamma}A_\alpha$,
$H=\oplus_{\alpha\in\Gamma}H_\alpha$ be braided groups and
$\rho:A\otimes H\to k$ a bialgebra pairing in $\C$
(where $\otimes$ means a tensor product over $k$).
We denote by $G$ a copy of $k\Gamma$ considered as an object of $\C$ with
the trivial $k\Gamma$-comodule structure (really, $G$ is a Hopf algebra in
$\C$).
$k\Gamma$-comodule property of $\rho$ means that
$\rho(a\otimes h)\neq 0$ for $a\in A_\alpha, h\in H_\beta$
only if $\alpha=\beta^{-1}$.
The formula
$\rho^\chi(a\otimes h):=\rho(a\otimes h)\cdot\alpha$ for such $a$ and $h$
define '$G$-valued pairing' $\rho^\chi:A\otimes H\rightarrow G$ in $\C$.
(Note that the underlying algebra of $G$ is commutative and there are no
problems to transform definition of bialgebra pairing from \ref{subsec_dual}
to this case.)
Let us also define right $G$-module structure on $X$ in $\C$ by the
formula
\begin{equation}
\label{mod-str-dub}
X_\alpha\otimes G\ni x_\alpha\otimes\beta\mapsto\chi(\alpha,\beta)\cdot
\chi(\beta,\alpha)\cdot x_\alpha\in X_\alpha\,.
\end{equation}
There exists
an algebra structure on $A\otimes G\otimes H$ which is uniquely determined
by the following conditions that $A$, $G$, $H$ are subalgebras of
$A\otimes G\otimes H$ with embeddings $A\otimes\eta_G\otimes\eta_H$,
$\eta_A\otimes G\otimes\eta_H$ and $\eta_A\otimes\eta_G\otimes H$
respectively, $G$ is in the center of $A\otimes G\otimes H$,
the product $a\cdot\alpha\cdot h$ for $a\in A$, $\alpha\in G$,
$h\in H$ equals to $a\otimes\alpha\otimes h$,
and the product
$h\cdot a:=(1_A\otimes 1_G\otimes h)\cdot (a\otimes 1_G\otimes 1_H)$
is defined by the diagram in Fig.\ref{Fig-Prod},
where crossings mean inverse to braiding in $\C$.

\begin{proposition}
The category $\DY{\cal C}_{A,H}$ is identified with a full subcategory
of $\C_{(A\otimes G\otimes H)}$ of $(A\otimes G\otimes H)$-modules such
that    underlying $G$-module structure coincides with
(\ref{mod-str-dub}).
\end{proposition}

Let us consider an example.
Let $\C$ be a category of $\ZZ$-graded vector spaces
$X=\oplus_{n\in\ZZ}X_n$ with braiding
\begin{equation}
\label{anyone-br}
x\otimes y\mapsto q^{mn}\cdot y\otimes x\qquad\hbox{for}\quad
x\in X_m,\enspace y\in Y_n,
\end{equation}
where $q\in k^*$ is a fixed invertible element such that $q^n\neq 1$ for
$n=1,2,\dots$.
Majid's braided line is the simplest example of a properly braided group.
This is a free $k$-algebra $A:=k[x]$ with one generator $x$ of degree $1$.
Counit, comultiplications and antipode are the following:
\begin{equation}
\label{br-line}
\epsilon(x^n)=\delta_{n,0}\,,\qquad
\Delta(x^n)=\Sigma_{k=0}^n\left[ n\atop k \right]_qx^k\otimes x^{n-k}\,,
\qquad
S(x^n)=(-1)^n\cdot q^{\frac{n(n-1)}{2}}\cdot x^n\,,
\end{equation}
where
$$
\left[ n\atop k\right]_q:=\frac{[n]_q!}{[k]_q!\cdot[n-k]_q!}\,,\qquad
[n]_q!:=\Pi_{k=1}^n[k]_q\,,\qquad
[n]:=\frac{q^n-1}{q-1}\,.
$$
We also consider a similar Hopf algebra $H$ in $\C$ with one generator
$y$ of degree $-1$.
Majid noted that there exists the following non-degenerated
bialgebra pairing between $A$
and $H$:
\begin{equation}
\rho:A\otimes H\rightarrow k\,,\qquad
x^m\otimes y^n\mapsto\delta_{m,n}\cdot [n]_q!\,.
\end{equation}
We want to describe the corresponding category $\DY{\C}_{A,H}$.
Let $G:=k[t,t^{-1}]$ be a group algebra of free group with one generator
$t$.
We consider $G$ as an object of $\C$ and suppose that degree of $t$ is $0$.
The formula in Fig.\ref{Fig-Prod} for multiplication in
$A\otimes G\otimes H$ takes the form

\begin{equation}
y^n\cdot x^m=
\sum_{{k,\ell,m^\prime,n^\prime\geq 0}\atop
      {{m^\prime+k+\ell=m}\atop {n^\prime+k+\ell=n} }}
 \,
(-1)^k\cdot q^{mn+[{k\atop 2}]-\ell(m^\prime+n^\prime)}\cdot
\frac{[m]_q!\cdot[n]_q!}
     {[m^\prime]_q!\cdot[n^\prime]_q!\cdot[k]_q!\cdot[\ell]_q!}\cdot
x^{m^\prime}\otimes t^\ell\otimes y^{n^\prime}\,.
\end{equation}
In particular, $yx=q(xy+t-1)$.

Let us consider a quiver with points labeled by integers and arrows
$x_n$, $y_n$, $n\in\ZZ$:
$$
  \cdots\;\circ_{n-1}\;
  \mathop{\longrightarrow\atop\longleftarrow}^{x_{n-1}}_{y_n}
  \;\circ_n\;
  \mathop{\longrightarrow\atop\longleftarrow}^{x_n}_{y_{n+1}}
  \;\circ_{n+1}\;\cdots
$$
And let $\cal A$ be the category generated by this quiver with relations
$y_{n+1}\cdot x_n=q\cdot x_{n-1}\cdot y_n+q(q^{2n}-1)\cdot 1_n$.

\begin{proposition}
The category $\DY{\C}_{A,H}$ is identified with the category of
representations of $\cal A$ over $\bf Vect$.
\par
Representation $\pi:{\cal A}\rightarrow {\bf Vect}$
corresponds to crossed module from
$\DY{\C}^A_A\hookrightarrow\DY{\C}_{A,H}$
(where embedding is described by the lemma \ref{Lemma-Pairing})
iff there exists $n_0\in\ZZ$ such that $\pi(n)$ is a zero vector space
for all $n<n_0$.
In particular, representations of $\cal A$ with vacuum state correspond to
simple $A$-crossed modules.
\end{proposition}

One can also consider the category of $\ZZ_n$-graded (anionic) spaces
with braiding defined by the formula (\ref{anyone-br}) with
$q^n=1$ and algebras $A={k[x]}/{(x^n)}$, $H={k[y]}/{(y^n)}$ in this
category with braided group structure defined by (\ref{br-line}).
A similar approach allows us to describe simple crossed modules over $A$.
This give us examples of anyonic $R$-matrix, in particular, those
described in \cite{Majid12}.

More detailedly these results will be published anywhere.


\subsec{}
For a pair of braided groups (bialgebras) $A$ and $B$ in $\cal C$ one can
consider a category $_A^A\DY{\cal C}^B_B$ whose objects $X$ are both left
crossed modules over $A$ with action $\mu_\ell$ and coaction $\Delta_\ell$
and right crossed modules over $B$ with action $\mu_r$ and coaction
$\Delta_r$ such that $(X,\mu_\ell,\mu_r)$ is $A$-$B$-bimodule,
$(X,\Delta_\ell,\Delta_r)$ is $A$-$B$-bicomodule and
$$
\Delta_r\circ\mu_\ell=(\mu_\ell\otimes B)\circ (A\otimes\Delta_r),\qquad
\Delta_\ell\circ\mu_r=(A\otimes\mu_r)\circ (\Delta_\ell\otimes B).
$$
$_A^A\DY{\cal C}^B_B$ is a monoidal category with standard tensor product
for underlying (co)modules and can be equipped with (pre-)braiding defined
in Fig.\ref{Fig-DY-Psi}c.
This and more complicated structures will be studied in forthcoming papers.

\subsec{}
A Hopf bimodules appeared (under the name 'bicovariant bimodules')
as the basic notion in Woronowicz approach to differential calculus
on quantum groups \cite{Wor}.
The main theorem of the paper \cite{Sch} can be considered as
a coordinate free version of the Woronowicz result about structure of Hopf
bimodules and states the equivalence between (pre-)braided categories of
Hopf bimodules and crossed modules over a Hopf algebra
in a symmetric monoidal category which has (co)equalizers.
In \cite{BD} the same result is proven in the context of braided groups
without assumption about existence of (co)equalizers.
Consideration of Hopf bimodules simplifies proofs of certain results about
crossed bimodules. Here we describe necessary facts from \cite{BD}.
\par
Firstly, we need an assumption that idempotents are split in our braided
category $\cal C$, i.e. for any idempotent (projection)
$e=e^2:\,X\rightarrow X$ there exist object $X_e$ and morphisms
$\displaystyle X_e\mathop{\longrightarrow\atop\longleftarrow}^{i_e}_{p_e} X$
such that $e=i_e\circ p_e$ and $p_e\circ i_e={\rm id}_{X_e}$.
This condition can be always satisfied:

\begin{lemma}
For any braided category $\cal C$ there exists braided category
$\widetilde{\cal C}$ and full embedding of braided
categories ${\cal C}\hookrightarrow\widetilde{\cal C}$
with idempotents in $\widetilde{\cal C}$ split.
\end{lemma}

\begin{proof}
We need only to prove that standard construction of idempotent splitting
is compatible with braided structure. Objects in the category
$\widetilde{\cal C}$ are pairs $X_e=(X,e)$, where $X$ is object in $\cal C$
and $e:\,X\rightarrow X$ is idempotent: $e^2=e$. Morphisms are the following
$$\widetilde{\cal C}(X_e,Y_f):=\{t\in{\cal C}(X,Y)\vert fte=t\}$$
with ordinary composition.
Tensor product and braiding are:
$$X_e\otimes Y_f:=(X\otimes Y)_{e\otimes f}$$
$$\Psi_{X_e,Y_f}:=
  (f\otimes e)\circ\Psi_{X,Y}\circ (e\otimes f)=
  (f\otimes e)\circ\Psi_{X,Y}=
  \Psi_{X,Y}\circ (e\otimes f)$$
The axioms of braided category are easily verified.
One can identify $\cal C$ with the full subcategory of $\widetilde{\cal C}$
whose objects are $(X,{\rm id}_X)$.
\end{proof}

We note that if idempotents are split in $\cal C$ then the same is true
in the categories ${\cal C}_A$, ${\cal C}^A$ and $\DY{\cal C}_A^A$.

\begin{figure}
$$
\matrix{
\vvbox{\hbox{\lu}
       \hhbox{\krl\ld}}
\enspace =\enspace
\vvbox{\hbox{\hcd\step\ld}
       \hbox{\id\step\hx\step\id}
       \hbox{\hcu\step\lu}}
\qquad\quad
\vvbox{\hbox{\step\ru}
       \hhbox{\krl\ld}}
\enspace=\enspace
\vvbox{\hbox{\ld\step\hcd}
       \hbox{\id\step\hx\step\id}
       \hbox{\hcu\step\ru}}
&\qquad\qquad&
{}_X\Pi:=\enspace
\vvbox{\hbox{\ld}
       \hbox{\S\step\id}
       \hbox{\lu}}
\cr \hbox{\scriptsize a) Hopf module axioms}&&
    \hbox{\scriptsize b) An idempotent ${}_X\Pi$} }
$$
\caption{}
\label{Fig-Hopf}
\end{figure}

\begin{definition}
{\em A left Hopf module} (resp. {\em right-left Hopf module})
$X$ over a bialgebra $A$ is a left (resp. right) $A$-module and left
$A$-comodule with the compatibility condition (Fig.\ref{Fig-Hopf}a)
which is a "polarized" version of the bialgebra axiom
(Fig.\ref{Fig-Main}c).
\par
{\em A Hopf bimodule}
$ X =( X , \mu_\ell, \mu_r,\Delta_\ell, \Delta_r)$
over a bialgebra $A$ is an object $X$ which is $A$-$A$-bimodule
and $A$-$A$-bicomodule and
such that any pair of action and coaction defines on $X$ structure of
a Hopf module.
\par
Hopf bimodules together with the $A$-bimodule-$A$-bicomodule morphisms
form the category which will be denoted by ${}_A^A{\cal C}^A_A$
\end{definition}

\begin{lemma}
Let $X$ be a Hopf bimodule over a braided group $A$.
Then the following are right crossed module structures on its underlying
object:
$X_{\rm ad}$ with adjoint action and regular coaction,
$X^{\rm ad}$ with regular action and adjoint coaction.
\par
Morphism ${}_X\Pi$ defined in Fig.\ref{Fig-Hopf}b is an idempotent and
${}_X\Pi:\,X^{\rm ad}\rightarrow X_{\rm ad}$ is a crossed module map.
\end{lemma}

One can define a crossed module ${}_AX$ ({\em the object of left invariants})
and crossed module morphisms
${}_Xp:\, X^{\rm ad}\rightarrow{}_AX$ and
${}_Xi:\,{}_AX\rightarrow X_{\rm ad}$ which split ${}_X\Pi$, i.e.
${}_X\Pi={}_Xi\circ{}_Xp$.

\begin{lemma}
Let $X=(X,\mu_\ell,\Delta_\ell,\mu_r,\Delta_r)$ be a Hopf bimodule and $Y$
a right crossed module over a bialgebra $A$.
Then $X\otimes Y$ equipped with right (co)module structure tensor product
one from $X$ and $Y$ and left action $\mu_\ell\otimes Y$ and coaction
$\Delta_\ell\otimes Y$ is a Hopf bimodule.
\end{lemma}

In the case when $X=A$ is a regular Hopf bimodule we use notation
$A\ltimes Y$ for tensor product equipped with a Hopf module structure
described in the previous lemma.

\begin{proposition}
Let $A$ be a braided group in $\cal C$.
Then the constructions descri\-bed above extend to functors
\begin{equation}
\label{Functor-Hopf-Cross}
{}_A^A{\cal C}^A_A
\mathop{\longrightarrow\atop\longleftarrow}^{{}_A(\_)}_{A\ltimes (\_)}
\DY{\cal C}^A_A
\end{equation}
\end{proposition}

We define the following tensor products on the category
${}_A^A{\cal C}_A^A$.
For any Hopf bimodules $(X,\mu_\ell,\Delta_\ell,\mu_r,\Delta_r)$ and
$(X^\prime,\mu_\ell^\prime,\Delta_\ell^\prime,\mu_r^\prime,\Delta_r^\prime)$
we denote by $\displaystyle X\mathop{\otimes}^\circ X^\prime$ (resp. by
$\displaystyle X\mathop{\otimes}_\circ X^\prime$) tensor product of
their underlying left and right modules (resp. comodules) equipped with
left and right comodule structures $\Delta_\ell\otimes X^\prime$,
$X\otimes\Delta_r^\prime$ (resp. module structures
$\mu_\ell\otimes X^\prime$, $X\otimes\mu_r^\prime$)
and put $\displaystyle f\mathop{\otimes}^\circ f^\prime:=f\otimes f^\prime=:
f\mathop{\otimes}_\circ f^\prime$ for Hopf bimodule maps $f$ and $f^\prime$.
One can directly verify that $\displaystyle\mathop{\otimes}^\circ$ and
$\displaystyle\mathop{\otimes}_\circ$ are bifunctors
${}^A_A{\cal C}_A^A\times{}^A_A{\cal C}_A^A\rightarrow{}^A_A{\cal C}_A^A$
satisfying the associativity condition.
But there are no units for these tensor products.

\begin{theorem}
Let $A$ be a braided group in $\cal C$.
Then there exist pre-braided structure
$\displaystyle ({}_A^A{\cal C}^A_A,\,\mathop{\otimes}_A,\,A,
 \Psi={}^{({}^A_A{\cal C}^A_A)}\Psi)$  (braided if antipode $S_A$ is
invertible)
and natural transformation of functors
$\displaystyle(\_)\mathop{\otimes}_\circ(\_)
 \buildrel{\phi^{\otimes_A}_{\_,\_}}\over\longrightarrow
 (\_)\mathop{\otimes}_A(\_)
 \buildrel{\phi^{\emptybox_A}_{\_,\_}}\over\longrightarrow
(\_)\mathop{\otimes}^\circ(\_)$
such that
\begin{itemize}
\item
any $\phi^\otimes_{X,X^\prime}$ (resp. $\phi_{X,X^\prime}^{\emptybox}$)
is a (co)tensor product over $A$, i.e. coequalizer of
$\displaystyle  X\otimes A\otimes  X^\prime
 \mathop{\longrightarrow\atop\longrightarrow}^{\mu_r\otimes X^\prime}%
_{X\otimes\mu_\ell^\prime}  X\otimes X^\prime$
(resp. equalizer of
$\displaystyle  X\otimes X^\prime
 \mathop{\longrightarrow\atop\longrightarrow}^{\Delta_r\otimes X^\prime}
                                                _{X\otimes\Delta_\ell^\prime}
  X\otimes A\otimes  X^\prime$);
\item
functors (\ref{Functor-Hopf-Cross}) define equivalence of (pre-)braided
categories.
\end{itemize}
\end{theorem}

\subsec{}
Important examples are crossed modules $A_{\rm ad}$, $A^{\rm ad}$ related
with the regular Hopf bimodule over a braided group $A$.
New solutions of Yang--Baxter equation arise as braiding between these
objects in the category of crossed modules \cite{Majid12}.
Connection of $A_{\rm ad}$ with bicovariant differential calculi on $A$
was pointed out (in coordinate form) by Woronowicz \cite{Wor} (see also
\cite{BD} for fully braided setting).
\par
A new concept of {\em a braided Lie algebra} was introduced and studied by
Majid \cite{M7}.
He showed that $A_{\rm ad}$ satisfy the axioms of braided Lie algebra if
underlying module of $A_{\rm ad}$ is an object of ${\cal C}_{\cO{A}}$.
He also built enveloping bialgebra $U({\cal L})$ for any braided Lie algebra
$\cal L$ and defined adjoint action of $U({\cal L})$ on $\cal L$ and
on itself.
We would like to note here that $U({\cal L})_{\rm ad}$ is a crossed module
over $U({\cal L})$ with crossed submodule $\cal L$.

\subsec{}
\label{Subsec-Antipode}
As was noted in \cite{Bespalov1}, one can define analog of antipode
(square of antipode) for any Hopf bimodule (resp. crossed module).
See \cite{BD} for more details about relative antipode.

\begin{definition}
Let $A$ be a braided group, $X$ a Hopf bimodule and $Y$ a right crossed
module over $A$.
{\em (Relative) antipode} $S=S_{X/A}$ for $X$ (resp. {\em square antipode}
$\sigma=\sigma_{Y/A}$ for $Y$) is defined by the first formula on
Fig.\ref{Fig-Antipodes}a (resp. in Fig.\ref{Fig-Antipodes}b),
where $3$-vertices denote compositions of left and right (co)actions.
\end{definition}

\begin{figure}
$$
\matrix{
S:=\enspace
\matrix{\step\object{X}\step\cr
\vvbox{\hbox{\ld\rd}
       \hbox{\S\step\id\step\S}
       \hbox{\lu\ru}}\cr
\step\object{X}\step}
\qquad\enspace
S^-:=\enspace
\vvbox{\hbox{\ld\rd}
       \hbox{\SS\step\id\step\SS}
       \hbox{\id\step\hxx}
       \hbox{\hxx\step\id}
       \hbox{\id\step\hxx}
       \hbox{\lu\ru}}
&\qquad\qquad&
\sigma:=\enspace
\matrix{\object{Y}\step\cr
\vvbox{\hbox{\rd}
       \hbox{\id\step\S}
       \hbox{\ru}}\cr
\object{Y}\step}
\qquad\enspace
\sigma^-:=\quad
\vvbox{\hbox{\rd}
       \hbox{\hxx}
       \hbox{\O{S^{-2}}\step\id}
       \hbox{\hxx}
       \hbox{\ru}}
\cr
\hbox{\scriptsize{ a) Relative antipode}}
&&
\hbox{\scriptsize{ b) Square antipode}}}
$$
\caption{}
\label{Fig-Antipodes}
\end{figure}

For the regular Hopf bimodule relative antipode coincide with the usual
Hopf algebra antipode.
The polarized forms of the two last identities in Fig.\ref{Fig-Main}d are
true:
\begin{equation}
\label{Rel-antip-prop}
S_{X/A}\circ\mu_r=\mu_\ell\circ\Psi\circ (S_{X/A}\otimes S_A)\qquad
\hbox{(and three similar identities).}
\end{equation}
If antipode of a braided group is invertible then
relative antipode for a Hopf bimodule also has inverse
given by the second formula in Fig.\ref{Fig-Antipodes}.
This follows directly from (\ref{Rel-antip-prop}).

Idempotent $\Pi_X$ for a Hopf bimodule $X$ over a braided group $A$
is defined by left-right reversed form of the diagram in Fig.\ref{Fig-Hopf}.
Relative antipode transforms idempotents $\Pi_ X$ and ${}_ X\Pi$
one to other:
\begin{equation}
S_{X/A}\circ\Pi_ X ={}_X\Pi\circ\Pi_X={}_ X\Pi\circ S_{X/A},
\qquad
S_{X/A}\circ{}_ X\Pi =\Pi_X\circ{}_X\Pi=\Pi_ X\circ S_{X/A}
\end{equation}

\begin{corollary} For a Hopf bimodule $X$ and for the corresponding right
crossed module ${}_AX$ the following identities are true:
\begin{equation}
S_{X/A}^2\circ{}_X\Pi={}_X\Pi\circ S_{X/A}^2\circ{}_X\Pi=
                                     {}_X\Pi\circ S_{X/A}^2\,,
\end{equation}
\begin{equation}
\label{Rel-to-Square}
{}_Xi\circ\sigma_{{{}_AX}/A}=S_{X/A}^2\circ{}_Xi\,,\qquad
\sigma_{{{}_AX}/A}\circ{}_Xp={}_Xp\circ S_{X/A}^2\,.
\end{equation}
\end{corollary}

\begin{corollary}
\hfill\hfill
Let $A$ be a braided group with invertible antipode and
\newline
$(Y,\mu^Y_r,\Delta^Y_r)$ a right crossed
module over $A$. Then
\begin{displaymath}
(Y,\,\mu^Y_r\circ\Psi^{-1}\circ (S^-\otimes Y),\,
                   (S\otimes Y)\circ\Psi\circ\Delta^Y_r)
\buildrel{\sigma_{Y/A}}\over\longrightarrow
(Y,\,\mu^Y_r\circ\Psi\circ (S\otimes Y),\,
                   (S^-\otimes Y)\circ\Psi^{-1}\circ\Delta^Y_r)
\end{displaymath}
is a morphism of left crossed modules.
Moreover a collection of $\sigma_{Y/A}$ describe isomorphism between two
corresponding braided functors from $\DY{\cal C}_A^A$ to ${}^A_A\DY{\cal C}$
defined in \ref{Subsec-New-Cross}.
\par
Morphism $\sigma_{Y/A}^-$ defined in Fig.\ref{Fig-Antipodes}b is inverse to
$\sigma_{Y/A}$.
\end{corollary}

\proof
Let us put $X=A\ltimes Y$. Then one can identify $Y$ with ${}_AX$.
Formulas (\ref{Rel-antip-prop}) imply that
\begin{displaymath}
S^2_{X/A}\circ\mu^X_r\circ\Psi^{-1}\circ(S^-_A\otimes{}_AX)=
\mu^X_r\circ\Psi\circ(S_A\otimes S^2_{X/A})\,.
\end{displaymath}
Composition with ${}_Xi$ and ${}_Xp$ and then application of
(\ref{Rel-to-Square}) prove that $\sigma_{{{}_AX}/A}$ is a modules morphism.
One can also find a direct proof of the (co-)module property of
$\sigma_{X/A}$ in the preprint version of this paper.
\par
\hfill
Module or comodule property of
$\sigma_{Y/A}$
imply that both
\hfill
$\sigma_{Y/A}\circ\sigma_{Y/A}^-$
\hfill
and
\newline
$\sigma_{Y/A}^-\circ\sigma_{Y/A}$ are equal to
\begin{displaymath}
\mu^Y_r\circ(\mu^Y_r\otimes A)\circ(Y\otimes S_A\otimes A)
\circ(\Delta^Y_r\otimes A)\circ\Delta^Y_r={\rm id}_Y\,.
\qquad\emptybox
\end{displaymath}
\medskip

For right crossed modules $A^{\rm ad}$ and $A_{\rm ad}$ we have
$\sigma_{{A^{\rm ad}}/A}=\sigma_{{A_{\rm ad}}/A}=S_A^2$
as in the case of ordinary Hopf algebra \cite{M8}.

\begin{proposition}
The following identity is true for any crossed modules $X$ and $Y$ from
$\DY{\cal C}_A^A$:
\begin{equation}
{}^{(\DY{\cal C}^A_A)}\Psi_{Y,X}\circ\sigma_{{(X\otimes Y)}/A}\circ
                                                {}^{(\DY{\cal C}^A_A)}
\Psi_{X,Y}=
{}^{\cal C}\Psi_{Y,X}\circ (\sigma_{Y/A}\otimes\sigma_{X/A})\circ
                                                 {}^{\cal C}\Psi_{X,Y}
\label{antipodetensor}
\end{equation}
\end{proposition}

\begin{figure}
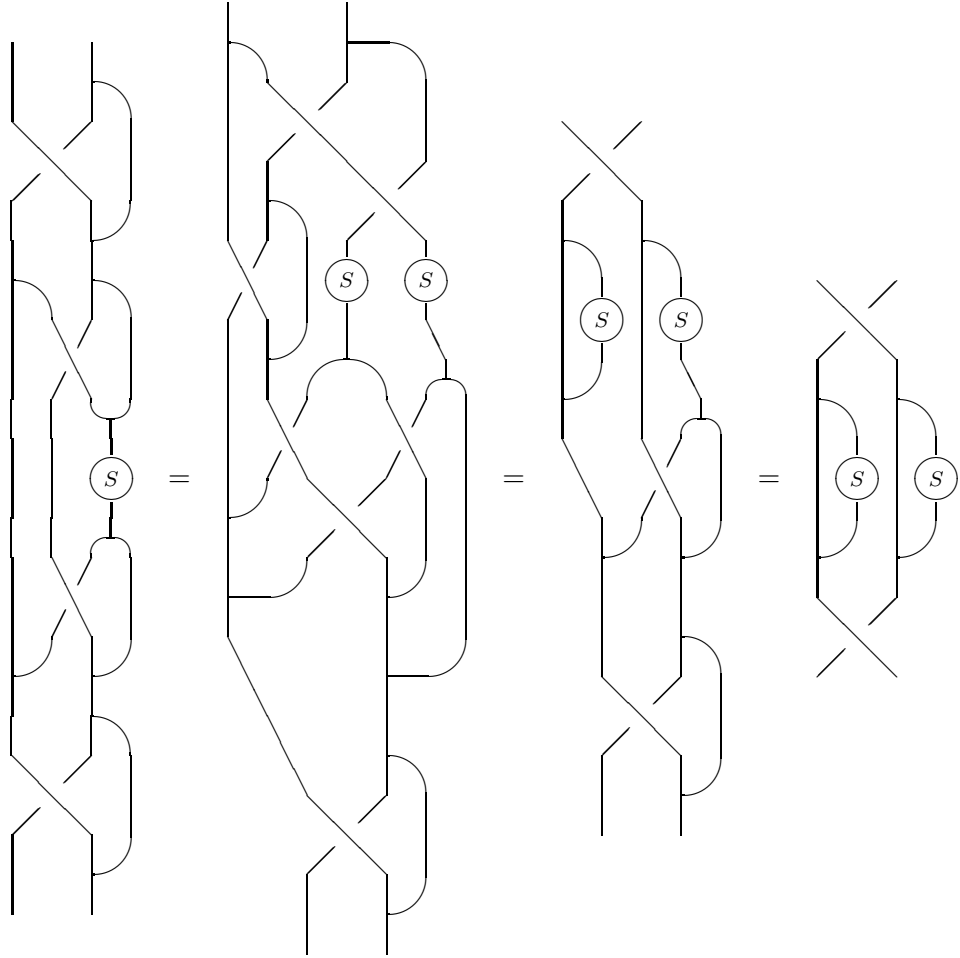

\begin{displaymath}
\vvbox{\hbox{\id\Step\rd}
       \hbox{\x\step\id}
       \hhbox{\id\Step\ru}
       \hbox{\rd\step\rd}
       \hbox{\id\step\hx\step\id}
       \hhbox{\id\step\id\step\cu}
       \hbox{\id\step\id\step\hstep\S}
       \hhbox{\id\step\id\step\cd}
       \hbox{\id\step\hx\step\id}
       \hbox{\ru\step\ru}
       \hhbox{\id\Step\rd}
       \hbox{\x\step\id}
       \hbox{\id\Step\ru}}
\quad=\quad
\vvbox{\hbox{\rd\Step\Rd}
       \hbox{\id\step\x\Step\id}
       \hbox{\id\step\rd\step\x}
       \hbox{\hx\step\id\step\S\Step\S}
       \hbox{\id\step\ru\cd\step\hdcd}
       \hbox{\id\step\hx\Step\hx\step\id}
       \hbox{\ru\step\x\step\id\step\id}
       \hbox{\Ru\Step\ru\step\id}
       \hbox{\d\Step\step\Ru}
       \hbox{\step\d\Step\rd}
       \hbox{\Step\x\step\id}
       \hbox{\Step\id\Step\ru}}
\quad=\quad
\vvbox{\hbox{\x}
       \hbox{\rd\step\rd}
       \hbox{\id\step\S\step\id\step\S}
       \hbox{\ru\step\id\step\hdcd}
       \hbox{\d\step\hx\step\id}
       \hbox{\step\ru\step\ru}
       \hbox{\step\id\Step\rd}
       \hbox{\step\x\step\id}
       \hbox{\step\id\Step\ru}}
\quad=\quad
\vvbox{\hbox{\x}
       \hbox{\rd\step\rd}
       \hbox{\id\step\S\step\id\step\S}
       \hbox{\ru\step\ru}
       \hbox{\x}}
\end{displaymath}
\caption{Proof of the identity
    ${}^{(\DY{\cal C}_H^H)}\Psi\circ S^2\circ{}^{(\DY{\cal C}_H^H)}\Psi=
          {}^{\cal C}\Psi\circ (S^2\otimes S^2)\circ{}^{\cal C}\Psi$:
\protect\newline
The first equality uses the facts that ${}^{(\DY{\cal C}_H^H)}\Psi$ is
a comodule morphism, antipode is anti-algebra morphism, the bialgebra and
module axioms.
The second uses the fact that antipode is anti-coalgebra morphism and
the antipode axiom.
The third uses the facts that ${}^{(\DY{\cal C}_H^H)}\Psi$ is a module
morphism,
antipode is anti-coalgebra morphism and the antipode axiom.}
\label{proofantipodetensor}
\end{figure}

\begin{proof}
See Fig.\ref{proofantipodetensor}.
\end{proof}

Square antipode $\sigma$ naturally appears in the formula for rank of crossed
module. Let $X\in{\rm Obj}({\cal C})$ and $X^\vee$ exists.
{\em A rank} or {\em categorical dimension} \cite{Majid7} of $X$
is the element of ${\rm End}(\underline 1)$ defined on
Fig.\ref{Fig-Rank}a.

\begin{proposition}
Let $X$ be a crossed module with dual $X^\vee$.
Then the identities in Fig.\ref{Fig-Rank}b are true.
\end{proposition}

\begin{figure}
$$
\matrix{
{\rm dim}_{\cal C}(X)
:=
\vvbox{\hbox{\obj{X}\hcoev}
       \hbox{\hx}
       \hbox{\hev}}
&\qquad&
\matrix{\object{X}\step\object{X^\vee}\cr
	  \vvbox{\hbox{\O{\sigma}\step\id}
		 \hbox{\hev}}}
=
\matrix{\object{X}\step\object{X^\vee}\cr
	  \vvbox{\hbox{\id\step\O{\sigma}}
		 \hbox{\hev}}}
\qquad
{\rm dim}_{\DY{\cal C}^A_A}(X):=\enspace
\vvbox{\hbox{\hcoev}
       \hbox{\O{\sigma}\step\id}
       \hbox{\hx}
       \hbox{\hev}}
\;=\;
\vvbox{\hbox{\hcoev}
       \hbox{\id\step\O{\sigma}}
	      \hbox{\hx}
	      \hbox{\hev}}
\cr
\hbox{\scriptsize a) A rank of $X\in{\rm Obj}({\cal C})$}&&
\hbox{\scriptsize b) A rank of $X\in{\rm Obj}(\DY{\cal C}^A_A)$} }
$$
\caption{}
\label{Fig-Rank}
\end{figure}

%
                          \sect{Cross products.}
%
\subsec{}
A structure of cross product algebra can be realized as a tensor product
algebra in a suitable braided category.

\begin{proposition}
\label{adjointalgebra}
(Cf. \cite{Majid12})
Let $A$ be a Hopf algebra in a braided category $\cal C$.
Then $A_{\rm ad}$ (resp. $A^{\rm ad}$) become commutative algebra
(resp. cocommutative coalgebra) in the category $\DY{\cal C}^A_A$ with
multiplication (resp. comultiplication) inherited from $A$.
\end{proposition}

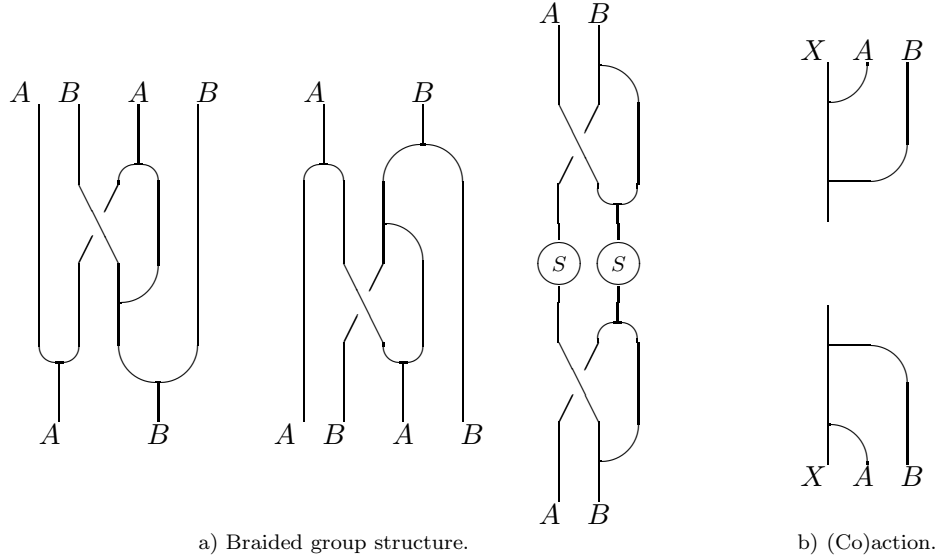
\begin{figure}
$$
\matrix{
\matrix{\object{A}\step\object{B}\step\hstep\object{A}%
\step\hstep\object{B}\cr
	\vvbox{\hbox{\id\step\id\step\hcd\step\id}
	       \hbox{\id\step\hx\step\id\step\id}
	       \hbox{\id\step\id\step\ru\step\id}
	       \hbox{\hcu\step\cu}}\cr
	\hstep\object{A}\hstep\Step\object{B}\step}
\qquad
\matrix{\hstep\object{A}\hstep\Step\object{B}\step\cr
	\vvbox{\hbox{\hcd\step\cd}
	       \hbox{\id\step\id\step\rd\step\id}
	       \hbox{\id\step\hx\step\id\step\id}
               \hbox{\id\step\id\step\hcu\step\id}}\cr
        \object{A}\step\object{B}\step\hstep\object{A}\step\hstep\object{B}}
\qquad
\matrix{\object{A}\step\object{B}\step\cr
	\vvbox{\hbox{\id\step\rd}
	       \hbox{\hx\step\id}
	       \hhbox{\id\step\cu}
	       \hbox{\S\step\hstep\S}
	       \hhbox{\id\step\cd}
	       \hbox{\hx\step\id}
	       \hbox{\id\step\ru}}\cr
        \object{A}\step\object{B}\step}
&\qquad\quad&
\matrix{\object{X}\step\object{A}\step\object{B}\cr
	\vvbox{\hbox{\ru\step\id}
               \hbox{\Ru}}
        \cr\begin{picture}(1,2)\end{picture}\cr
        \vvbox{\hbox{\Rd}
               \hbox{\rd\step\id}}\cr
        \object{X}\step\object{A}\step\object{B}}
\cr
\hbox{\scriptsize a) Braided group structure.}&&
\hbox{\scriptsize b) (Co)action.}}
$$
\caption{Cross product $A\ltimes B$.}
\label{Fig-Cross-Prod}
\end{figure}

Let $B$ be a Hopf algebra in a category $\DY{\cal C}_A^A$.
Denote by $A\ltimes B$ object $A\otimes B$ equipped with
algebra structure $A_{\rm ad}\otimes B$,
coalgebra structure $A^{\rm ad}\otimes B$
(tensor product algebra and coalgebra in $\DY{\cal C}_A^A$) and antipode
$$S_{A\ltimes B}:=
  {}^{(\DY{\cal C}_A^A)}\Psi_{B,A_{\rm ad}}\circ (S_B\otimes S_A)\circ
                         {}^{(\DY{\cal C}_A^A)}\Psi_{A^{\rm ad},B}.$$
(It's easy to see that there exists ${S_{A\ltimes B}}^-$ if
$S_A^-$ and $S_B^-$ exist because all factors in the last formula are
invertible.)
Formulas for multiplication, comultiplication and antipode
in terms of the category $\cal C$ are given in Fig.\ref{Fig-Cross-Prod}a.

\begin{theorem}
$A\ltimes B$ is a Hopf algebra in a category $\cal C$.
\label{CrossProduct}
\end{theorem}

\begin{proof}
Underlying algebra of $A\ltimes B$ is exactly Majid's cross product by
braided group in \cite{Majid6} so associativity is clear.
The coalgebra is precisely dual construction by input-output symmetry.
Nontrivial part is a verification of the bialgebra axiom
(see Fig.\ref{Proof-Bialg}).
\end{proof}

This is the analog of construction of such cross products and cross
coproducts for ordinary Hopf algebras in \cite{Radford1}, \cite{M1}.
The analog of the converse, which is the Majid-Radford theorem
\cite{Radford1}, \cite{M6} is also true in our braided situation.

Let $A$ and $H$ be bialgebras in $\cal C$.
A bialgebra projection is a pair of bialgebra morphisms
\begin{equation}
\label{Projection}
A\buildrel{i_A}\over\rightarrow
 H\buildrel{p_A}\over\rightarrow A
\end{equation}
such that the composition homomorphism
$p_A\circ i_A$ equals ${\rm id}_A$.
A direct proof of the following theorem was given in the preprint
version of this paper.

\begin{theorem}\label{Radford}
Let $\cal C$ be a braided category with split idempotents, $A$ and $H$
braided groups in $\cal C$, and bialgebra projection
(\ref{Projection}) be given.
Then there exists a braided group $B$ living in the category
$\DY{\cal C}_A^A$ such that $H\simeq A\ltimes B$.
\end{theorem}

\subsec{}
For given a bialgebra projection (\ref{Projection}) one can turn $H$
into an object $\underline H$ of ${}_A^A{\cal C}_A^A$ taking composition
of the regular actions (resp. coactions) with $i_A$ (resp. $p_A$).
Both theorems
\ref{CrossProduct} and \ref{Radford} are corollaries of the following
theorem and the result about equivalence of the categories ${}^A_A{\cal
C}^A_A$ and $\DY{\cal C}^A_A$.

\begin{figure}
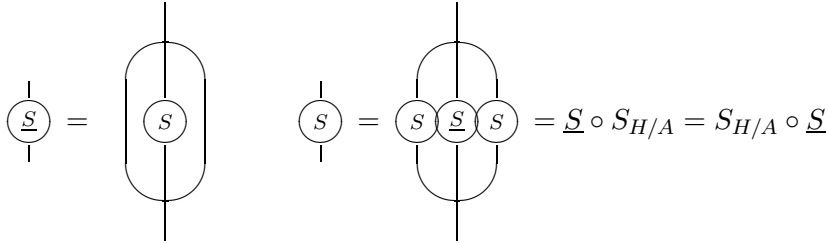

$$
\vvbox{\hbox{\O{\underline S}}}\quad =\quad
\vvbox{\hbox{\ld\rd}
       \hbox{\id\step\S\step\id}
       \hbox{\lu\ru}}
\qquad\qquad
\vvbox{\hbox{\S}}\quad =\quad
\vvbox{\hbox{\ld\rd}
       \hbox{\S\step\O{\underline S}\step\S}
       \hbox{\lu\ru}}
\quad =\underline S\circ S_{H/A}=S_{H/A}\circ\underline S
$$
\caption{Relations between antipodes on $H$ and $\underline H$. }
\label{Antipodes-BMR}
\end{figure}

\begin{theorem}(cf. \cite{BD})
Let $(A,\mu_A,\eta_A,\Delta_A,\epsilon_A)$ be a braided group in $\cal C$.
For any bialgebra $(H,\mu_H,\eta_H,\Delta_H,\epsilon_H)$ and bialgebra
projection (\ref{Projection}) morphism $\mu$
(resp. $\Delta$) is uniquely factorized through (co-)tensor
product over $A$:
\begin{equation}
\mu_H=\underline\mu\circ\phi^{\otimes_A}_{\underline H,\underline H},
\qquad
\qquad
\Delta_H =\phi_{\underline H,\underline H}^{\emptybox_A}\circ%
\underline\Delta\,,
\label{BMR}
\end{equation}
And
$\underline H:=
 (H,\underline\mu,\underline\eta,\underline\Delta,\underline\epsilon )$
is a bialgebra in ${}_A^A{\cal C}^A_A$, where $\underline\eta:=i_A$ and
$\underline\epsilon:=p_A$.
It this way one to one correspondence between bialgebra projections
(\ref{Projection}) and bialgebras in ${}_A^A{\cal C}^A_A$ is defined.
\par
There also exists one to one correspondence between
braided group structures on $H$ and $\underline H$. The relations
between corresponding antipodes $S$ and $\underline S$ are shown
in Fig.\ref{Antipodes-BMR}.
\end{theorem}

The idempotent ${}_H\Pi$ for $H\in{\rm Obj}({}^A_A{\cal C}^A_A)$ from
the previous theorem takes the form
\begin{equation}
\label{Idemp}
{}_H\Pi=\mu_H\circ(i_A\circ S_A\circ p_A\otimes H)\circ\Delta_H\,.
\end{equation}
Let
$B\displaystyle\mathop{\longrightarrow\atop\longleftarrow}^{i_B}_{p_B} H$
split this idempotent.

\begin{corollary}
\label{Cor-Radford}
$B$ is a bialgebra in the category $\DY{\cal C}^A_A$ with
multiplication $p_B\circ\mu_H\circ(i_B\otimes i_B)$,
comultiplication $(p_B\otimes p_B)\circ\Delta_H\circ i_B$,
right $A$-module structure $p_B\circ\mu_H\circ(i_B\otimes i_A)$ and
right $A$-comodule structure $(p_B\otimes p_A)\circ\Delta_H\circ i_B$.
This bialgebra satisfy the theorem \ref{Radford}.
\end{corollary}

Let $A$ be a braided group in a category $\cal C$, $B$ a Hopf algebra in
$\DY{\cal C}^A_A$ and $X$ be a right module over $B$ in $\DY{\cal C}^A_A$.
It is easy  to verify that the first diagram in Fig.\ref{Fig-Cross-Prod}b
defines ($A\!\ltimes\! B$)-module structure on $X$.
Moreover, in this way one can construct full embedding of categories
$(\DY{\cal C}^A_A)_B\hookrightarrow{\cal C}_{A\ltimes B}$

\begin{proposition}
Let $A$ be a braided group in $\cal C$ and $B$ a braided group in
$\DY{\cal C}^A_A$.
Then the braided categories
$(\DY{\cal C}^A_A)^B_B$ and $\DY{\cal C}^{A\ltimes B}_{A\ltimes B}$ are
isomorphic.
\end{proposition}

\begin{figure}
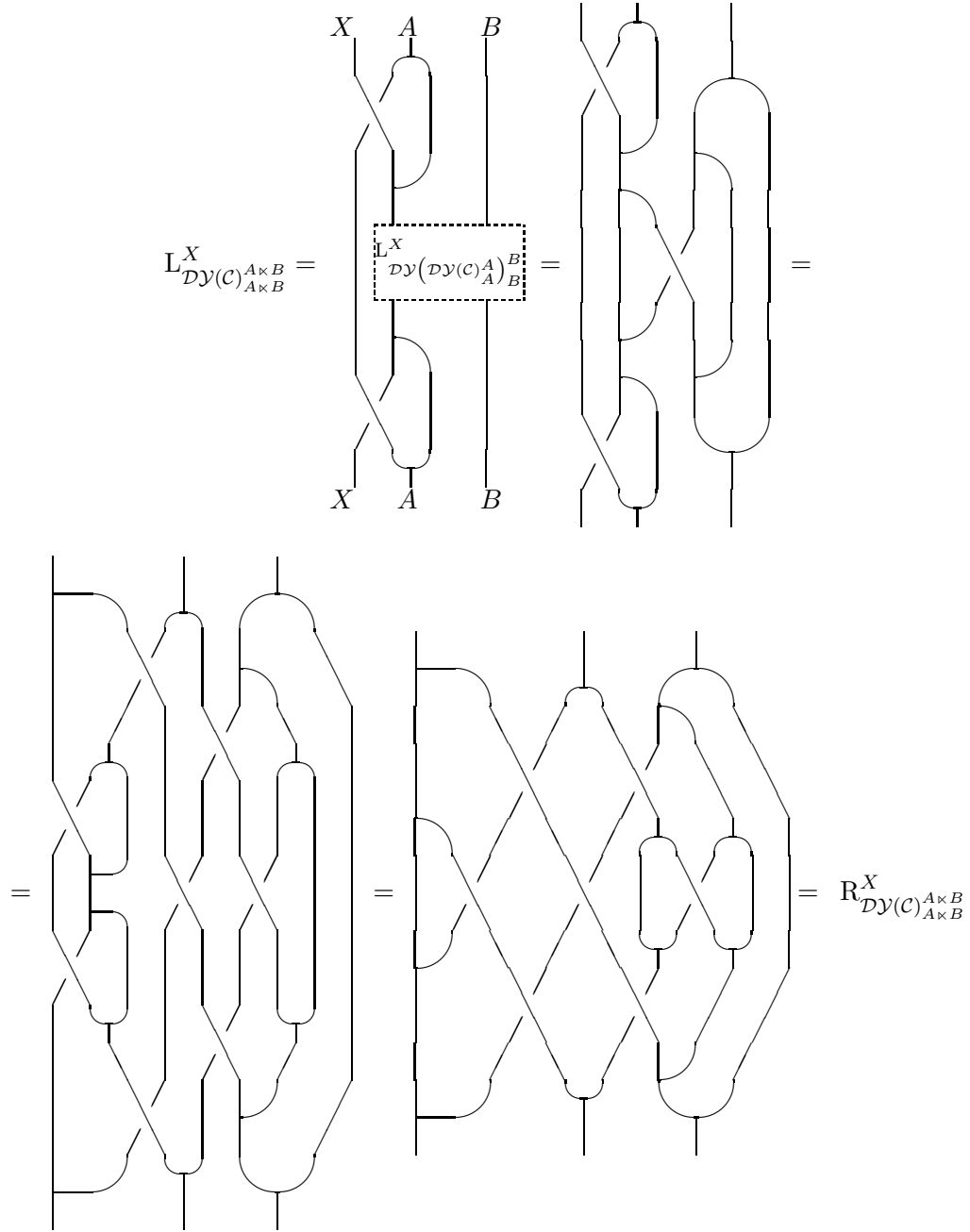

\begin{displaymath}
{\rm L}^X_{\DY{\cal C}^{A\ltimes B}_{A\ltimes B}}=\enspace
\matrix{\object{X}\step\hstep\object{A}\Step\object{B}\cr
        \vvbox{\hhbox{\krl\id\step\cd\hstep\step\id}
               \hbox{\hx\step\id\hstep\step\id}
               \hbox{\id\step\ru\hstep\step\id}
               \hbox{\id\step\hstep%
\Dash{{\rm L}^X_{\DY{\DY{\cal C}^A_A}^B_B}}}
               \hbox{\id\step\rd\hstep\step\id}
               \hbox{\hx\step\id\hstep\step\id}
               \hhbox{\krl\id\step\cu\hstep\step\id}}\cr
        \object{X}\step\hstep\object{A}\Step\object{B}}
\quad=\enspace
\vvbox{\hhbox{\id\step\cd\Step\id}
       \hbox{\hx\step\id\step\cd}
       \hbox{\id\step\ru\step\rd\step\id}
       \hhbox{\id\step\rd\step\id\step\id\step\id}
       \hbox{\id\step\id\step\hx\step\id\step\id}
       \hhbox{\id\step\ru\step\id\step\id\step\id}
       \hbox{\id\step\rd\step\ru\step\id}
       \hbox{\hx\step\id\step\cu}
       \hhbox{\id\step\cu\Step\id}}
\enspace=
\end{displaymath}
\begin{displaymath}
=\enspace
\vvbox{\hbox{\Rd\step\hcd\step\cd}
      \hbox{\id\Step\hx\step\id\step\rd\step\d}
      \hbox{\id\step\hddcd\step\id\step\hx\step\hdcd\step\id}
      \hbox{\hx\step\id\step\id\step\id\step\id\step\id\step\id\step\id}
      \hbox{\id\step\k\step\hx\step\hx\step\id\step\id}
      \hbox{\hx\step\id\step\id\step\id\step\id\step\id\step\id\step\id}
      \hbox{\id\step\hdcu\step\id\step\hx\step\hddcu\step\id}
      \hbox{\id\Step\hx\step\id\step\ru\step\dd}
      \hbox{\Ru\step\hcu\step\cu}}
\enspace=\enspace
\vvbox{\hbox{\Rd\Step\hcd\step\hstep\cd}
       \hhbox{\id\Step\d\step\dd\step\d\step\rd\step\d}
       \hbox{\id\Step\hstep\hx\Step\hx\step\d\hstep\d}
       \hhbox{\rd\step\dd\step\d\step\dd\hstep\cd\step\cd\step\id}
       \hbox{\id\step\hx\Step\hx\step\id\step\hx\step\id\step\id}
       \hhbox{\ru\step\d\step\dd\step\d\hstep\cu\step\cu\step\id}
       \hbox{\id\Step\hstep\hx\Step\hx\step\dd\hstep\dd}
       \hhbox{\id\Step\dd\step\d\step\dd\step\ru\step\dd}
       \hbox{\Ru\Step\hcu\step\hstep\cu}}
=\enspace{\rm R}^X_{\DY{\cal C}^{A\ltimes B}_{A\ltimes B}}
\end{displaymath}
\caption{Proof of the crossed module axiom over cross product $A\ltimes B$:
\protect\newline
The second equality is the crossed module axiom for $X$ over $B$.
The third uses the fact that $B$-action is $A$-module and
$B$-coaction is $A$-comodule morphisms respectively.
The 4th is the crossed module axiom for $X$ over $A$.
The 5th uses the bialgebra axiom for $A$.}
\label{Proof-over-Cross}
\end{figure}

\begin{proof}
($A\!\ltimes\! B$)-(co-)module structure on object $X$ of
$\DY{\DY{\cal C}_A^A}^B_B$
is defined in Fig.\ref{Fig-Cross-Prod}b.
Nontrivial part is a verification of crossed module axiom for $X$
in Fig.\ref{Proof-over-Cross}.
Conversely, $A$-, ($B$-) crossed module structure on objects of
$\DY{\cal C}^{A\ltimes B}_{A\ltimes B}$ is defined by means of $i_A$, $p_A$
($i_B$, $p_B$).
\end{proof}

So for object $C$ in $\cal C$ to be a braided group in
$\DY{\cal C}^{A\ltimes B}_{A\ltimes B}$ or in $\DY{\DY{\cal C}^A_A}^B_B$
are equivalent. In this case on $A\otimes B\otimes C$ there exist two
braided group structures $(A\ltimes B)\ltimes C$ and
$A\ltimes (B\ltimes C)$.

\begin{proposition} {\em (Transitivity of cross product)}
Braided groups $(A\ltimes B)\ltimes C$ and $A\ltimes (B\ltimes C)$ coincide.
\label{transitivity}
\end{proposition}

\subsec{}
The following results about cross products involve a square antipode
$\sigma_{X/A}$.

\begin{lemma}
\label{cross-sigma}
Let $H=A\ltimes B$ be a cross product of braided groups $A$ in $\cal C$
and $B$ in $\DY{\cal C}_A^A$.
Then for any object $X$ of $\DY{\cal C}^H_H=\DY{\DY{\cal C}^A_A}^B_B$
\begin{equation}
\sigma_{X/H}=\sigma_{X/A}\circ\sigma_{X/B}=\sigma_{X/B}\circ\sigma_{X/A}.
\end{equation}
\end{lemma}

\begin{lemma}
\label{cross-square}
Square of antipode in a cross product $A\ltimes B$ of braided groups
has the form
\begin{equation}
{S_{A\ltimes B}}^2=(S_A^2\otimes S_B^2)\circ\Psi^2\circ
(A\otimes\sigma_{B/A}).
\end{equation}
\end{lemma}

%
              \sect{Remarks on Quantum braided groups.}
%

\subsec{}
\label{element-u}
We note that
in definition of quasitriangular bialgebra $(A,\overline A,{\cal R})$
coassociativity and bialgebra axiom for $\overline A$ follow from
other axioms.
\par
Really, the bialgebra axiom is versified directly
(taking into attention that
${\cal R}^-$ is a copairing between ${{\overline A}_{\rm op}}^{\rm op}$
and $A^{\rm op}$):
$$\overline\Delta^{\rm op}\circ\mu_A=
{\cal R}\cdot(\Delta\circ\mu_A)\cdot{\cal R}^-
={\cal R}\cdot(\mu_{A\otimes A}\circ(\Delta\otimes\Delta))\cdot{\cal R}^-=
$$
$$ =\mu_{A\otimes A}\circ
     ({\cal R}\cdot\Delta\cdot{\cal R}^-\otimes
    {\cal R}\cdot\Delta\cdot{\cal R}^-)=
   \mu_{A\otimes A}\circ(\overline\Delta^{\rm op}\otimes
                                       \overline\Delta^{rm op}) $$
Then one can use axioms of quasitriangular structure, coassociativity of
$\Delta$, bialgebra axiom for $A$ and $\overline A$ to prove
coassociativity of $\overline\Delta^{\rm op}$:
$$(\overline\Delta^{\rm op}\otimes A)\circ\overline\Delta^{\rm op}=
 {\cal R}_{23}\cdot{\cal R}_{13}\cdot{\cal R}_{12}\cdot
 ((\Delta\otimes{\rm id})\circ\Delta)\cdot
 {\cal R}_{12}^-\cdot{\cal R}_{13}^-\cdot{\cal R}_{23}^-=
 (A\otimes\overline\Delta^{\rm op})\circ\overline\Delta^{\rm op}\,.$$

The following obvious lemma is important for our further considerations.

\begin{lemma}
If $(A,\overline A,{\cal R})$ is a quantum
braided group (quasitriangular bialgebra) in $\cal C$ then
$(\overline A,A,\overline{\cal R})$ is a quantum braided group
(quasitriangular bialgebra) in $\overline {\cal C}$.
(Pre-)braided categories $\overline{\cal C}_{\cO{\overline A,A}}$
and $\overline{{\cal C}_{\cO{A,\overline A}}}$ are naturally identified.
\end{lemma}

\subsec{}
Here we describe relations between antipodes in quantum braided group
and distinguished 'element' $u$.
These are analogues of the results obtained in \cite{Drinfel'd2}, \cite{RT},
\cite{R1} for ordinary quantum groups.

Let $(A,\overline A,{\cal R})$ be a quantum braided group and
$X\in{\rm Obj}({\cal C}_{\cO{A,\overline A}})$.
Then $\sigma_{{X^{\cal R}}/A}$ is a result of action $\triangleleft u$,
where $u=\mu_A\circ{\cal R}^\sim$ (see Fig.\ref{Fig-Element-u}).
Denote by $\overline u$ corresponding elements for
quantum braided group $(\overline A,A,\overline {\cal R})$.

\begin{figure}
$$u=\enspace
\vvbox{\hbox{\ro{{\cal R}^\sim}}
       \hbox{\id\Step\id}
       \hbox{\cu}}
\enspace=\enspace
\vvbox{\hbox{\r}
       \hbox{\id\Step\S}
       \hbox{\cu}}
\quad=\quad
\vvbox{\hbox{\r}
       \hbox{\tSS\Step\id}
       \hbox{\cu}}
\qquad\quad
u^-=\quad
\vvbox{\hbox{\r}
       \hbox{\S\Step\SS}
       \hbox{\cu}}
\quad=\quad
\vvbox{\hbox{\r}
       \hbox{\tS\Step\tSS}
       \hbox{\cu}}
\qquad\quad
$$
\caption{}
\label{Fig-Element-u}
\end{figure}

\begin{proposition}
'Elements' $u,\overline u$  have inverse $u^-,(\overline u)^-$
(where $u^-$ is defined in Fig.\ref{Fig-Element-u}).
Antipodes in quantum braided group are invertible.
The following relations are true:
$$
\overline S^-=u\cdot S\cdot u^-,\qquad
S^-=\overline u\cdot\overline S\cdot(\overline u)^-\,,
$$
$$
\overline u=S^-\circ u^-=\overline S\circ u^-\,,
\qquad
\overline u^-=S^-\circ u=\overline S\circ u\,,
$$
where $u\cdot S$ is abbreviated notation for $\mu\circ (u\otimes S)$,
etc.
\par
The following (braided variant of the formula for $\Delta u$) is
true for any modules $X$ and $Y$ from ${\cal C}_{\cO{A,\overline A}}$:
\begin{equation}
{}^{({\cal C}_{\cO{H}})}\Psi_{Y,X}\circ(\triangleleft u)
                                  \circ{}^{({\cal C}_{\cO{H}})}\Psi_{X,Y}=
{}^{\cal C}\Psi_{Y,X}\circ ((\triangleleft u)\otimes(\triangleleft u))
                                  \circ{}^{\cal C}\Psi_{X,Y}
\end{equation}
\end{proposition}

\begin{proof}
The axiom of quasitriangular structure in Fig.\ref{Fig-QBG}a
implies directly that
\begin{equation}
\mu\circ(A\otimes u\cdot S)\circ\overline\Delta^{\rm op}=u\circ\epsilon\,,
\qquad
\mu\circ(S\cdot u^-\otimes A)\circ\overline\Delta^{\rm op}=u^-\circ\epsilon
\label{equ-u1}
\end{equation}
with $u^-:=\mu\circ(S\circ\overline S\otimes A)\circ{\cal R}$
Application of $\overline S$ to the first identity from (\ref{equ-u1})
and then convolution product
with ${\rm id}_A$ give:
$u\cdot(S\circ\overline S)={\rm id}_A\cdot u.$
The later identity implies that
$u^-:=\mu\circ(S\circ\overline S\otimes A)\circ{\cal R}$ is right inverse
to $u$.
Then it follows from (\ref{equ-u1}) that $u\cdot S\cdot u^-$ is a
skew antipode for $\overline A$.
This and the same considerations for $(\overline A, A, \overline{\cal R})$
prove the first part of the theorem.
\par
The second part is a special case of
(\ref{antipodetensor}).
\end{proof}

For the ordinary quantum group our $(\overline u)^-$ coincide with $u$ from
\cite{Drinfel'd2}, \cite{RT}.

\subsec{}
Let $(A,\overline A,{\cal R})$ be a quasitriangular bialgebra in $\cal C$.
It is convenient to describe the category
${\cal C}_{\cO{A,\overline A}}$ in terms of crossed modules.
For right $A$-module $X$ denote by $X^{\cal R}$ object with additional
comodule structure defined by the last diagram in Fig.\ref{Fig-Mod}a.
An idea to turn a module into a comodule by $\cal R$ is due to Majid
\cite{M1}.
A proof of the following lemma is obtained immediately in a similar way
from axioms of quasitriangular structure.

\begin{lemma}
\hfill
Let $Y$ be a right module over $A$.
Then
\hfill
$Y\in{\rm Obj}({\cal C}_{\cO{A,\overline A}})$
\hfill
iff
\newline
$Y^{\cal R}\in{\rm Obj}\left(\DY{\cal C}^A_A\right).$
\par
In this case
$(X\otimes Y)^{\cal R}=X^{\cal R}\otimes Y^{\cal R}$
for any module $X$.
\end{lemma}

So one can identify ${\cal C}_{\cO{A,\overline A}}$ with a full
braided monoidal subcategory of $\DY{\cal C}_H^H$ and braiding coincides
with that defined by Majid (Fig.\ref{Fig-QBG}c).

\begin{lemma}
Let $(A,\overline A,{\cal R})$ be a quantum braided group in $\cal C$
and $X\in{\rm Obj}({\cal C}_{\cO{A,\overline A}})$.
Then the identities in Fig.\ref{Fig-Change-Antipode} are true.
(In the second case we suppose that there exist right dual $X^\vee$ in
$\cal C$.)
\end{lemma}

\begin{figure}
$$
\matrix{\object{A}\step\object{X}\cr
\vvbox{\hbox{\tSS\step\id}
       \hbox{\hx}
       \hbox{\ru}}}
\quad ={}\quad
\matrix{\object{A}\step\object{X}\cr
\vvbox{\hbox{\SS\step\id}
       \hbox{\hxx}
       \hbox{\ru}}}
\qquad\qquad
\matrix{\object{X}\step\object{X^\vee}\step\object{A}\cr
\vvbox{\hbox{\id\step\id\step\S}
       \hbox{\id\step\hx}
       \hbox{\ru\step\id}
       \hbox{\ev}}}
\quad ={}\quad
\matrix{\object{X}\step\object{X^\vee}\step\object{A}\cr
\vvbox{\hbox{\id\step\id\step\tS}
       \hbox{\id\step\hxx}
       \hbox{\ru\step\id}
       \hbox{\ev}}}
$$
\caption{}
\label{Fig-Change-Antipode}
\end{figure}

\begin{figure}
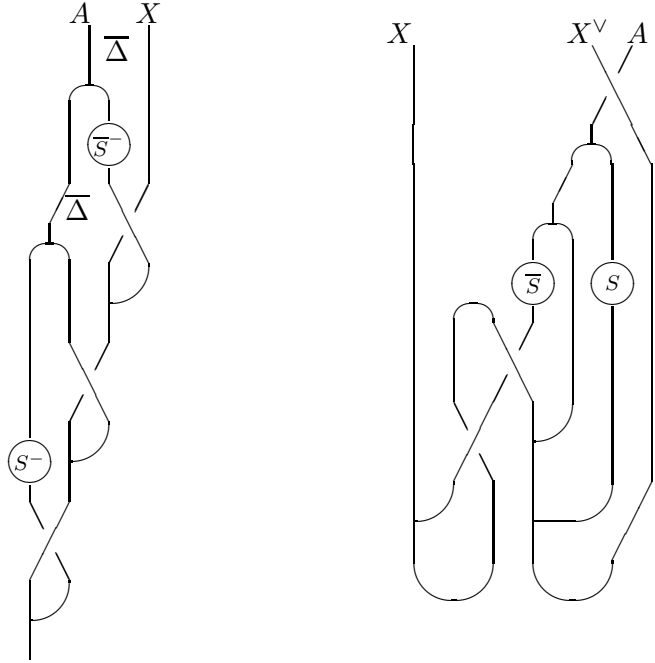

$$\matrix{\step\hstep\object{A}\hstep\step\object{X}\cr
\vvbox{\hbox{\step\hcd\hstep\Obj{\overline\Delta}\hstep\id}
       \hbox{\step\id\step\tSS\step\id}
       \hbox{\hddcd\hstep\Obj{\overline\Delta}\hstep\hx}
       \hbox{\id\step\id\step\ru}
       \hbox{\id\step\hx}
       \hbox{\SS\step\ru}
       \hbox{\hxx}
       \hbox{\ru}}}
\qquad\qquad
\qquad\qquad
\matrix{\object{X}\Step\Step\hstep\object{X^\vee}\step\object{A}\hstep\cr
\vvbox{\hbox{\id\Step\Step\hstep\hx}
       \hhbox{\id\Step\Step\cd\hstep\d}
       \hbox{\id\Step\step\hddcd\step\id\step\id}
       \hbox{\id\step\hcoev\step\tS\step\id\step\S\step\id}
       \hbox{\id\step\id\step\hx\step\id\step\id\step\id}
       \hbox{\id\step\hxx\step\ru\step\id\step\id}
       \hbox{\ru\step\id\step\Ru\dd}
       \hbox{\ev\step\ev}}}
$$
\caption{Proof of the identities in
Fig.\protect\ref{Fig-Change-Antipode}.}
\label{Proof-Change-Antipode}
\end{figure}

\begin{proof}
Both parts of the first (resp. the second) identity on
Fig.\ref{Fig-Change-Antipode} equals the first (resp. the second)
diagram in Fig.\ref{Proof-Change-Antipode} .
\end{proof}

\begin{corollary}
Let $X$ be a right module from ${\cal C}_{\cO{A,\overline A}}$ and there
exists
left dual ${}^\vee X$ (resp. right dual $X^\vee$) then
$({}^\vee X)^{\cal R}={}^\vee (X^{\cal R})$
(resp. $(X^\vee )^{\cal R}=(X^{\cal R})^\vee$) $\quad{\rm and}\quad
{}^\vee X$ (resp. $X^\vee$) belongs to ${\cal C}_{\cO{A,\overline A}}$.
\end{corollary}

\begin{figure}
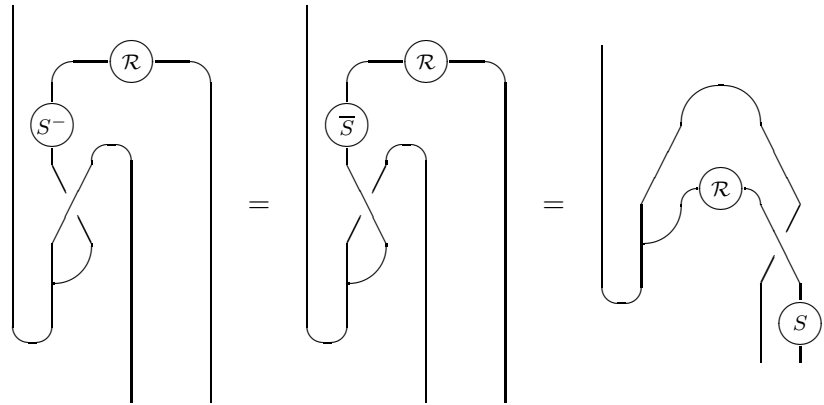

\begin{displaymath}
\vvbox{\hbox{\id\step\R}
       \hbox{\id\step\SS\step\hcoev\Step\id}
       \hbox{\id\step\hxx\step\id\Step\id}
       \hbox{\id\step\ru\step\id\Step\id}
       \hbox{\hev\Step\id\Step\id}}
\quad ={}\quad
\vvbox{\hbox{\id\step\R}
       \hbox{\id\step\tS\step\hcoev\Step\id}
       \hbox{\id\step\hx\step\id\Step\id}
       \hbox{\id\step\ru\step\id\Step\id}
       \hbox{\hev\Step\id\Step\id}}
\quad ={}\quad
\vvbox{\hbox{\id\Step\coev}
       \hbox{\id\step\dd\r\d}
       \hbox{\id\step\ru\Step\hx}
       \hbox{\hev\Step\step\id\step\S}}
\end{displaymath}
\caption{Comodule structures on $({}^\vee X)^{\cal R}$ and
         ${}^\vee (X^{\cal R})$ are the same:
\protect\newline
The first identity follows from the first identity in
Fig.\ref{Fig-Change-Antipode}.}
\label{proofld}
\end{figure}

\begin{proof}
Proposition is a direct corollary of the previous lemma and the dual
form of identities (\ref{anti-pair})
$(A\otimes S^{\pm 1})\circ{\cal R}=
 (\overline S^{\mp 1}\otimes A)\circ{\cal R}$.
For example, proof for the left dual is in Fig.\ref{proofld}.
\end{proof}

\subsec{}
Obviously,
one can consider left-right reversed form of quantum braided group
(quantum braided group in ${\cal C}_{\rm op}$),
corresponding subcategory of left modules and
embedding of this category into the category of left crossed modules.
Let $(A,\overline A,{\cal R})$ be a quantum braided group in $\cal C$.
Then $(A_{\rm op},\overline A_{\rm op},\overline R)$ is a quantum braided
group in ${\cal C}_{\rm op}$.
Then, it follows from \ref{Subsec-New-Cross} that we can turn any object
$\left(X,\mu_\ell,\Delta_\ell:=
(A\otimes\mu_\ell)\circ(\overline{\cal R}\otimes X)\right)$ of
$\overline{\cal C}_\cO{A,\overline A}\subset
{}^{A_{\rm op}}_{A_{\rm op}}\DY{\overline {\cal C}}$
into an object of ${}^A_A\DY{\cal C}$ with modified coaction
$\Psi\circ (X\otimes S)\circ\Psi\circ\Delta_\ell$,
and this construction is extended to isomorphism of braided categories.
We show in Fig.\ref{Fig-Alt}
comodule structure on an object of
${}_\cO{\overline A_{\rm op},A_{\rm op}}\overline{\cal C}$
and braiding on
${}_\cO{\overline A_{\rm op},A_{\rm op}}\overline{\cal C}$
\cite{Majid6}
obtained in this way.

\begin{figure}
$$
\matrix{
\vvbox{\hbox{\r\step\id}
       \hbox{\id\Step\hx}
       \hbox{\Lu\step\id}
       \hbox{\Step\hx}}
&\qquad\qquad&
\matrix{\Step\step\object{X}\step\object{Y}\cr
        \vvbox{\hbox{\r\step\id\step\id}
	       \hbox{\id\Step\hx\step\id}
	       \hbox{\Lu\step\lu}
	       \hbox{\Step\x}}}
\cr
}
$$
\caption{}
\label{Fig-Alt}
\end{figure}
All our further results about bosonization and transmutation have
analogues for left modules from
${}_\cO{\overline A_{\rm op},A_{\rm op}}\overline{\cal C}$.
But formulas in terms of initial quantum braided group
$(A,\overline A,{\cal R})$ for the last case are more complicated
(corresponding diagrams have more crossings).
So we prefer to work with right modules.

\sect{Generalized bosonization and transmutation.}

\subsec{}
Here we give cross product construction for quantum
groups in braided categories which generalize Majid's bosonization
construction $A\ltimes B$ to the case where $A$ is a quantum braided group
rather than a quantum group as in \cite{Majid6}.
The usual bosonization assigns this ordinary Hopf algebra to any
braided-Hopf algebra in the category of $A$-modules. Majid also showed in
\cite{Majid6} that when $B$ is braided quasitriangular then $A\ltimes B$
has a 'relative quasitriangular structure' in the sense of a second
coproduct and an ${\cal R}_{A\ltimes B}$.
This construction extends straightforwardly to our setting where $A$ is
itself a quantum braided group.
We also extend the braided version of Majid-Radford theorem in section $4$
to the case where the braided groups are equipped with quasitriangular
structures respected by the projections.

\begin{figure}
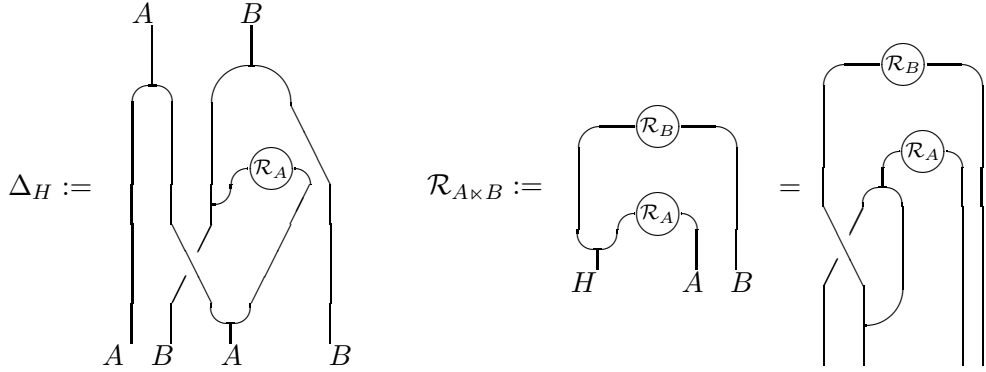

$$
\Delta_H:=\;
\matrix{\hstep\object{A}\Step\hstep\object{B}\Step\cr
        \vvbox{\hbox{\hcd\step\cd}
               \hbox{\id\step\id\step\id\hstep\ra\hhstep\d}
               \hhbox{\id\step\id\step\hru\step\hstep\dd\hstep\id}
               \hbox{\id\step\hx\step\dd\step\id}
               \hhbox{\id\step\id\step\cu\Step\id}}\cr
	\object{A}\step\object{B}\step\hstep\object{A}\Step\hstep\object{B}}
\qquad\enspace
{\cal R}_{A\ltimes B}:=\;
\matrix{
\vvbox{\hbox{\Rb}
       \hbox{\id\step\ra\step\id}
       \hhbox{\krl\cu\Step\id\step\id}}\cr
\hstep\object{H}\hstep\Step\object{A}\step\object{B}}
\enspace=\enspace
\vvbox{\hbox{\Rb}
       \hbox{\id\step\hstep\ra\hstep\id}
       \hhbox{\id\step\cd\step\hstep\id\hstep\id}
       \hbox{\hx\step\id\step\hstep\id\hstep\id}
       \hbox{\id\step\ru\step\hstep\id\hstep\id}}
$$
\caption{Quantum braided group structure on cross product $A\ltimes B$}
\label{Fig-QBG-Cross-Prod}
\end{figure}

Let $(A,\overline A,{\cal R}_A)$ be a quantum braided group in $\cal C$
and $(B,\overline B,{\cal R}_B)$ a quantum braided group in
${\cal C}_{\cO{A,\overline A}}$.
Then one can use embeddings
${\cal C}_{\cO{A,\overline A}}\rightarrow\DY{\cal C}_A^A$
to construct cross-product braided group $H:=A\ltimes B$ in $\cal C$.
Combining formulas for comultiplication from Fig.\ref{Fig-Cross-Prod}a and
for comodule structure from Fig.\ref{Fig-Mod}a one can obtain an expression
for comultiplication in $A\ltimes B$ (Fig.\ref{Fig-QBG-Cross-Prod}).
Braided group
$\overline H:=\overline A\ltimes\overline B$ in $\overline{\cal C}$
is built in the same way.
Let us define a quasitriangular structure ${\cal R}_H$
as a result of multiplication applied to ${\cal R}_A$ and ${\cal R}_B$ as
shown in Fig.\ref{Fig-QBG-Cross-Prod} (we consider $A$ and $B$ as subalgebras
of $H$).
Written explicitly this is exactly the form \cite{Majid8} of a cross
coproduct by the induced coaction as for the usual bosonization theorem.

\begin{theorem}
$(A\ltimes B,\overline A\ltimes\overline B,{\cal R}_{A\ltimes B})$
is a quantum braided group in $\cal C$.
\par
Braided categories
${\cal C}_{\cO{A\ltimes B,\overline A\ltimes\overline B}}$
and
$({\cal C}_{\cO{A,\overline A}})_{\cO{B,\overline B}}$
are isomorphic.
\label{quantumcross}
\end{theorem}

\begin{proof}
'Elements' ${\cal R}_{A\ltimes B}$ and
$\overline{\cal R}_{A\ltimes B}^{\rm op}:=
\Psi\circ{\cal R}_{\overline A\ltimes\overline B}$
are inverse to each other:
\begin{displaymath}
\overline{\cal R}_{A\ltimes B}^{\rm op}
\cdot {\cal R}_{A\ltimes B}=
\dots\hbox{ \ See Fig.\ref{randinverse} \ }\dots=
\eta_{A\ltimes B}\otimes\eta_{A\ltimes B}\,.
\end{displaymath}
It is shown in Fig.\ref{proofpairing}a) that ${\cal R}_H$ is
an algebra-coalgebra copairing.
After series of transformations one can verify that both
$\overline\Delta^{\rm op}_{A\ltimes B}\cdot{\cal R}_{A\ltimes B}$
and ${\cal R}_{A\ltimes B}\cdot\Delta_{A\ltimes B}$ are equal to
morphism given by the diagram in Fig.\ref{proofpairing}c)
(everywhere dot means multiplication in braided tensor product algebra
in corresponding category $\C$ or ${\cal C}_{\cO{A,\overline A}}$).
\par
For any bialgebra $A$ in $\C$ and right $A$-module $(X,\Delta^X_r)$
let us temporarily denote by $\psi^\C_{X,A}$ the composition morphism
$(A\otimes\Delta^X_r)\circ{}^\C\Psi_{X,A}\circ(X\otimes\Delta_A)$
To proof the second part of the theorem we note that in our case
for any object $X$ of the category ${\cal C}_{A\ltimes B}=({\cal C}_A)_B$
the following identities are true:
\begin{eqnarray}
\psi^\C_{X,A\ltimes B}&=&(A\otimes\psi^{\C_\cO{A,\overline A}}_{X,B})\circ
(\psi^\C_{X,A}\otimes B)\,,\nonumber\\
\label{mod-rest}
\psi^{\overline\C}_{X,\overline A\ltimes\overline B}&=&
(A\otimes\psi^{\overline\C_\cO{\overline A,A}}_{X,B})\circ
(\psi^{\overline \C}_{X,\overline A}\otimes B)\,.
\end{eqnarray}
But condition that $X\in{\rm Obj}({\cal C}_{\cO{H,\overline H}})$ (resp.
$X\in{\rm Obj}(({\cal C}_{\cO{A,\overline A}})_{\cO{B,\overline B}})\;)$
is exactly coincidence of the left hand sides (resp. the right hand sides)
of (\ref{mod-rest}).
\end{proof}

As in the case of braided groups (theorem \ref{transitivity})
cross product of quantum braided groups is transitive:

\begin{theorem}
Let $(A,\overline A,{\cal R}_A)$ be a quantum group in $\cal C$,
$(B,\overline B,{\cal R}_B)$ a quantum group in
${\cal C}_{\cO{A,\overline A}}$,
$(C,\overline C,{\cal R}_C)$ a quantum group in
${\cal C}_{\cO{A\ltimes B,\overline A\ltimes\overline B}}=
({\cal C}_{\cO{A,\overline A}})_{\cO{B,\overline B}}$.
Then quantum groups
$\left((A\ltimes B)\ltimes C,
 (\overline A\ltimes\overline B)\ltimes\overline C,
 {\cal R}_{(A\ltimes B)\ltimes C}\right)$
and
\newline
 $\left(A\ltimes (B\ltimes C),
  \overline A\ltimes (\overline B\ltimes\overline C),
  {\cal R}_{A\ltimes (B\ltimes C)}\right)$
coincide.
\end{theorem}

\begin{proof}
It is a corollary of the theorem \ref{transitivity}.
One needs only to verify that
${\cal R}_{(A\ltimes B)\ltimes C}={\cal R}_{A\ltimes (B\ltimes C)}$
using associativity of multiplication in $A\ltimes B\ltimes C$.
\end{proof}

\begin{definition}
Let $(A,\overline A,{\cal R}_A)$ and $(H,\overline H,{\cal R}_H)$
be quantum groups in $\cal C$.
Pair of morphisms
$\displaystyle A\mathop{\longrightarrow\atop\longleftarrow}^{i_A}_{p_A} H$
is called {\em a quantum group projection} if:
\begin{itemize}
\item
both
$\displaystyle A\mathop{\longrightarrow\atop\longleftarrow}^{i_A}_{p_A} H$
and
$\displaystyle\overline A
 \mathop{\longrightarrow\atop\longleftarrow}^{i_A}_{p_A}\overline H$
are bialgebra projections in $\cal C$ and $\overline{\cal C}$ respectively.
\item
\hfill
${}_H\Pi={}_{\overline H}\Pi$\,,
\hfill
\newline
where the first (resp. the second) is
an idempotent for
$H\in\Obj({}^A_A\C^A_A)$
defined by (\ref{Idemp})
(resp. corresponding idempotent for
$\overline H\in\Obj
({}^{\overline A}_{\overline A}{\overline \C}^{\overline A}_{\overline A})$)
(For ordinary quantum groups $A$ and $H$
this condition is satisfied automatically.);
\item and
\begin{equation}
\label{R-A-H}
\!\!\!\!
(H\otimes p_A)\circ{\cal R}_H=(i_A\otimes A)\circ{\cal R}_A \,,
\quad\enspace
(p_A\otimes H)\circ{\cal R}_H=(A\otimes i_A)\circ{\cal R}_A \,.
\end{equation}
\end{itemize}
\end{definition}

\begin{theorem}
\hfill
Let $\cal C$ be a braided category with split idempotents and
\newline
$\displaystyle (A,\overline A,{\cal R}_A)
  \mathop{\longrightarrow\atop\longleftarrow}^{i_A}_{p_A}
  (H,\overline H,{\cal R}_H)$
a quantum group projection in $\cal C$.
Then there exists quantum group $(B,\overline B,{\cal R}_B)$
in the category
${\cal C}_{\cO{A,\overline A}}$ such that
\begin{equation}
(H,\overline H,{\cal R}_H)\simeq
(A\ltimes B,\overline A\ltimes\overline B,{\cal R}_{A\ltimes B})\,.
\end{equation}
\end{theorem}

\begin{proof}
Let
$\displaystyle H \mathop{\longrightarrow\atop\longleftarrow}^{p_B}_{i_B} B$
split idempotent ${}_H\Pi$ in $\cal C$.
According to theorem \ref{Radford} and corollary \ref{Cor-Radford}
one can equip $B$ with a structure of Hopf algebra in $\DY{\C}^A_A$
and with a structure of Hopf algebra $\overline B$
in $\DY{\overline\C}^{\overline A}_{\overline A}$ with the same underline
algebra.
Canonical Hopf algebra isomorphisms $H\simeq A\ltimes B$ and
$\overline H\simeq\overline A\ltimes\overline B$ are defined.
And we will identify $H$ with $A\ltimes B$.
\par
The axiom of quasitriangular structure (Fig.\ref{Fig-QBG}a) for
${\cal R}_H$ implies that
\begin{equation}
\label{B-in-Mod}
(p_B\otimes p_A)\circ ({\cal R}_H\cdot\Delta_H)\circ i_B=
(p_B\otimes p_A)\circ (\overline\Delta_H^\op\cdot{\cal R}_H)\circ i_B
\end{equation}
One can use sequentially the fact that $p_A:H\rightarrow A$ is an algebra
morphism, identities (\ref{R-A-H})
and
$p_B\circ\mu_H\circ(i_H\otimes H)=\epsilon_A\otimes p_B$
(or, respectively the dual form
$(p_A\otimes H)\circ\overline\Delta_H\circ i_B=\eta_A\otimes i_B$)
to show that the left hand side of (\ref{B-in-Mod}) equals to
$(p_B\otimes p_B)\circ\Delta_H\circ i_B$
(that is a right $A$-comodule structure on $B$)
and, respectively, the right hand
side of (\ref{B-in-Mod}) equals to
$(p_B\circ\mu_H\circ(i_B\otimes i_A)\otimes A)\circ(B\otimes {\cal R}_A)$
(that is the composition of a right $A$-module structure on $B$ with the
quasitriangular structure ${\cal R}_A$).
And so $B$ is an object of the subcategory
$\C_\cO{A,\overline A}\subset\DY{\C}^A_A$.
\par
Let us define
\begin{equation}
{\cal R}_B:=(p_B\otimes p_B)\circ{\cal R}_H\,.
\end{equation}
One can see that
\begin{eqnarray}
(i_B\otimes B)\circ{\cal R}_B&=&(\Pi_B\otimes p_B)\circ{\cal R}_H
\nonumber\\
&=&(H\otimes p_B)\circ
 \left(((i_A\circ\overline S_A\circ p_A\otimes H)\circ{\cal R}_H)\cdot
                                              {\cal R}_H\right)
\nonumber\\
&=&
(H\otimes p_B)\circ
 \left(((i_A\circ\overline S_A\otimes i_A)\circ{\cal R}_A)\cdot
                                              {\cal R}_H\right)
\nonumber\\
&=&(H\otimes p_B)\circ{\cal R}_H\,.
\end{eqnarray}
And taking into account this identity
we obtain the expression for ${\cal R}_H$
via ${\cal R}_A$ and ${\cal R}_B$ as in Fig.\ref{Fig-QBG-Cross-Prod}:
\begin{displaymath}
{\cal R}_H=(H\otimes(p_A\otimes p_B)\circ\Delta_H)\circ{\cal R}_H=
(H\otimes p_A\otimes p_B)\circ(({\cal R}_H)_{13}\cdot({\cal R}_H)_{23})=
({\cal R}_B)_{13}\cdot({\cal R}_A)_{23}
\end{displaymath}
The identity
\begin{equation}
(p_B\otimes p_B)\circ({\cal R}_H\cdot\Delta_H)\circ i_A=
(p_B\otimes p_B)\circ(\overline\Delta_H^{\rm op}\cdot{\cal R}_H)\circ i_A
\end{equation}
imply that ${\cal R}_B$ is $A$-module morphism.
Compositions with $p_B^{\otimes 2}$ or $p_B^{\otimes 3}$ turn axioms
of quasitriangular structure for ${\cal R}_H$ into corresponding axioms
for ${\cal R}_B$.
\end{proof}

\subsec{}
A notion of a transmutation, a basic theory and
non-trivial examples of braided groups obtained via transmutation
belong to Majid \cite{M19,Majid7,Majid6}.
We describe here a braided version of Majid's transmutation
and show that analogues of Majid results
(in particular, transitivity of transmutation) keep in this situation.
Transmutation of a crossed module is defined and compatibility
of transmutation with cross product of braided groups is proven.
In particular, connection between transmutation and bosonization noted in
\cite{Majid6} is generalized for the fully braided setting.

\begin{definition}
Let $(A,\overline A,{\cal R}_A)$ be a quantum braided group in $\cal C$,
$(Y,\mu_\ell^Y,\mu_r^Y)$ bimodule over its underlying algebra.
We say that $\mu_{\rm ad}^Y: Y\otimes A\rightarrow Y$ is {\em a (right)
adjoint action of quantum braided group $(A,\overline A,{\cal R}_A)$
on bimodule $Y$} if $\mu_{\rm ad}^Y$ is adjoint action on $Y$ for both
bialgebras $A$ and $\overline A$.
\par
We denote by ${}_A{\cal C}_{\cO{A,\overline A}}$ a full subcategory of the
category ${}_A{\cal C}_A$ of bimodules $Y$ over $A$ such that there
exists adjoint action $\mu_{\rm ad}^Y$ of $(A,\overline A,{\cal R}_A)$ on
$Y$ and $(Y,\mu_{\rm ad}^Y)\in{\rm Obj}({\cal C}_{\cO{A,\overline A}})$.
Let, moreover, $(X,\mu^X_r)\in{\rm Obj}(\C_\cO{A,\overline A})$.
We define morphism
\begin{equation}
\tau_{X,Y}:=
(\mu^X_r\otimes\mu^Y_r)\circ(X\otimes{\cal R}_A^\sim\otimes Y):\;
X\otimes Y\rightarrow X\otimes Y\,.
\end{equation}
\end{definition}

\begin{figure}
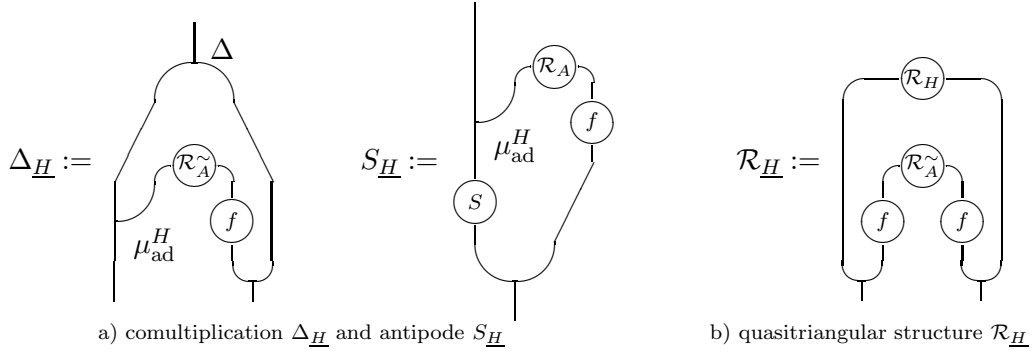

$$
\matrix{
\Delta_{\underline H}:=\enspace
\vvbox{\hbox{\step\cd\Obj{\Delta}}
       \hbox{\dd\ro{{\cal R}_A^\sim}\d}
       \hbox{\ru\obj{\mu^H_{\rm ad}}\Step\O{f}\step\id}
       \hhbox{\krl\id\Step\step\cu}}
\qquad\quad
S_{\underline H}:=\quad
\vvbox{\hbox{\id\step\ro{{\cal R}_A}}
       \hbox{\ru\obj{\mu^H_{\rm ad}}\Step\O{f}}
       \hbox{\S\Step\dd}
       \hbox{\cu}}
&\qquad&
{\cal R}_{\underline H}:=\enspace
\vvbox{\hbox{\Ro{{\cal R}_H}}
       \hbox{\id\step\ro{{\cal R}_A^\sim}\step\id}
       \hbox{\id\step\O{f}\Step\O{f}\step\id}
       \hhbox{\krl\cu\Step\cu}}%
\cr
\hbox{\scriptsize a) comultiplication $\Delta_{\underline H}$
                      and antipode $S_{\underline H}$}
&&
\hbox{\scriptsize b) quasitriangular structure ${\cal R}_{\underline H}$}
}
$$
\caption{Transmutation $\underline H$ of (quantum) braided group $H$.}
\label{trunsdef}
\end{figure}

\begin{definition}
\label{transmutation}
Let $(A,\overline A,{\cal R}_A)$ be a quantum braided group and
$(H,\mu_H,\Delta_H)$
a bialgebra (braided group) in $\cal C$ and $f:\, A\rightarrow H$
a bialgebra morphism.
Then $H$ becomes a bimodule over $A$ with actions
$\mu^H_\ell:=\mu_H\circ (f\otimes H),$ and
$\mu^H_r:=\mu_H\circ (H\otimes f)$.
We say that
$\left((A,\overline A,{\cal R}_A),f,H\right)$ are
{\em transmutation data for a bialgebra (braided group) $H$ in $\C$} if
$(H,\mu^H_\ell,\mu^H_r)$ is an object of
${}_A{\cal C}_{\cO{A,\overline A}}$.
Let us denote by $\mu^H_{\rm ad}$ adjoint action of $A$ on $H$.
{\em Transmutation $\underline H$ of bialgebra (braided group) $H$}
is the underlying algebra of $H$ with new comultiplication
$\Delta_{\underline H}:=\tau_{H,H}\circ\Delta_H$
(and antipode $S_{\underline H}$) defined on
Fig.\ref{trunsdef}a.
\par
\hfill
We denote by
\hfill
$\DY{\C}^H_{H,\cO{A,\overline A}}$
\hfill
full subcategory of
\hfill
$\,\DY{\C}^H_H$
\hfill
with objects
\newline
$(X,\mu^X_r,\Delta^X_r)$ such that
$X$ with $A$-module structure $\mu^X_r\circ(X\otimes f)$ belongs to
$\C_\cO{A,\overline A}$.
{\em A transmutation} $\underline X$ of a such crossed module
is the underline module of $X$ with a new 'coaction'
$\Delta^{\underline X}_r:=\tau_{X,H}\circ\Delta^X_r$.
\end{definition}

\begin{theorem}
Transmutation $\underline H$ of a bialgebra (resp. braided group) $H$
in $\C$ is a bialgebra (resp. braided group) in $\C_\cO{A,\overline A}$.
\par
A full embedding of (pre-)braided categories
$\DY{\cal C}^H_{H,\cO{A,\overline A}}\rightarrow
 \DY{{\cal C}_\cO{A,\overline A}}^{\underline H}_{\underline H}$
is defined by assignment $X\mapsto\underline X$ on objects and identity
on morphisms.
\end{theorem}

\proof
It follows directly from the definition of adjoint action that
$(H,\mu^H_{\rm ad})$ becomes
an algebra and $(X,\mu^X_r)$ a right a $H$-module in
$\C_\cO{A,\overline A}$.
\par
The next important step is to show that $\Delta^{\underline X}_r$ is
$A$-module morphism. A proof use the definition of adjoint action of quantum
braided group $(A,\overline A,{\cal R}_A)$ on $H$ and the fact that $X$ is
an object of both categories $\DY{\C}^H_H$ and $\C_\cO{A,\overline A}$.
Application of this result to crossed module $H_{\rm ad}$ implies that
$\Delta_{\underline H}$ is also $A$-module morphism.
\par
$A$-module property of $\Delta_{\underline H}$
(resp. of $\Delta^{\underline X}_r$), bialgebra axiom for $H$ and the fact
that ${\cal R}_A$ is a bialgebra copairing imply the coassociativity of
$\Delta_{\underline H}$ (resp. $\underline H$-comodule axiom for
$\underline X$).
\par
The bialgebra axiom for $\underline H$ is a corollary of the following more
general result (in the unbraided form noted by Majid):
Let $B$ be a right $H$-module algebra such that its underlying object with
right $A$-module structure defined via $f$ is an object of
$\C_\cO{A,\overline A}$.
Then $\tau_{B,H}:\,B\otimes H\rightarrow\underline B\otimes\underline H$
is an algebra isomorphism, where the first (resp. the second) is the
tensor product algebra in $\C$ (resp. in $\C_{A,\overline A}$).
\par
The antipode axiom
$\mu_H\circ(S_{\underline H}\otimes H)\circ\Delta_{\underline H}=
\eta_H\circ\epsilon_H$ follows directly from the definition.
To prove that
$\mu_H\circ(H\otimes S_{\underline H})\circ\Delta_{\underline H}=
\eta_H\circ\epsilon_H$ one needs to rewrite the left hand side of this
identity in the form
$\mu^H_r\circ(\dots\circ\mu^H_r\otimes A)\circ(H\otimes{\cal R}_A)$
using $A$-module property of $\Delta_{\underline H}$ and the fact that
${\cal R}_A$ is a coalgebra-algebra copairing.
The antipode axiom for $S_{\underline H}$ and
$A$-module properties of $\Delta_{\underline H}$ and $\mu_H$ imply that
$S_{\underline H}$ is $A$-module morphism also.
\par
The following identities prove the crossed module axiom for $\underline X$:
\begin{displaymath}
{\rm L}^{\underline X}%
_{\DY{\C_\cO{A,\overline A}}^{\underline H}_{\underline H}} =
\tau_{X,H}\circ{\rm L}^X_{\DY{\C}_H^H}=
\tau_{X,H}\circ{\rm R}^X_{\DY{\C}_H^H}=
{\rm R}^{\underline X}%
_{\DY{\C_\cO{A,\overline A}}^{\underline H}_{\underline H}}\,.
\quad\emptybox
\end{displaymath}

\begin{definition}
Let $(A,\overline A,{\cal R}_A)$ be a quantum braided group and
$(H,\overline H,{\cal R}_H)$ a quasitriangular bialgebra
(resp. a quantum braided group) in $\cal C$. We say that a morphism
$f:\, A\rightarrow H$ defines {\em  transmutation data for
a quasitriangular bialgebra
(resp. for a quantum braided group) $(H,\overline H,{\cal R}_H)$}
if
$\left((A,\overline A,{\cal R}_A),f,H\right)$ are
transmutation data for a bialgebra (braided group) $H$ in $\C$
and
$\left((\overline A,A,\overline{{\cal R}_A}),f,\overline H\right)$ are
transmutation data for a bialgebra (braided group) $\overline H$ in
$\overline\C$.
In this case we define {\em a transmutation of}
$(H,\overline H,{\cal R}_H)$
as a triple which consists of $\underline H$ (transmutation of $H$),
$\underline{\overline H}$ (transmutation of $\overline H$) and
${\cal R}_{\underline H}$
(transmutation of quasitriangular structure defined in Fig.\ref{trunsdef}b).
\end{definition}

\begin{theorem}
$\left(\underline H,\underline{\overline H},
{\cal R}_{\underline H}\right)$
is a quasitriangular bialgebra (resp. a quantum braided group) in
${\cal C}_\cO{A,\overline A}$.
\end{theorem}

We have the following commutative diagram of full embeddings of
(pre-)braided categories
\begin{equation}
\matrix{
\C_\cO{H,\overline H}&\longrightarrow&
\left(\C_\cO{A,\overline A}\right)
_\cO{\underline H,\overline{\underline H}}\cr
\downarrow&&\downarrow\cr
\DY{\C}^H_{H,\cO{A,\overline A}}&\longrightarrow&
\DY{\C_\cO{A,\overline A}}^{\underline H}_{\underline H}
}
\end{equation}

\begin{theorem}
(Transitivity of transmutation)\newline
Let $(A_i,\overline A_i,{\cal R}_i)\enspace i=0,1,2$ be quantum braided
groups in $\cal C$ and
morphisms $f:\,A_0\rightarrow A_1$, $g:\,A_1\rightarrow A_2$ be such that
the following triples are transmutation data:
\begin{displaymath}
\left((A_0,\overline A_0,{\cal R}_0),\,f,\,
 (A_1,\overline A_1,{\cal R}_1)\right)\,,\qquad
\left((A_0,\overline A_0,{\cal R}_0),\,g\circ f,\,
 (A_2,\overline A_2,{\cal R}_2)\right)
\end{displaymath}
and let
$(\underline{A_i},\underline{\overline A_i},\underline{{\cal R}_i})
\enspace i=1,2$
be corresponding transmutations
which are quantum braided groups in
$\C_\cO{A_0,\overline A_0}$.
Then
$\left((\underline{A_1},\underline{\overline A_1},\underline{{\cal R}_1}),\,
g,\,
(\underline{A_2},\underline{\overline A_2},\underline{{\cal R}_2})\right)$
are transmutation data iff
$\left(({A_1},{\overline A_1},{{\cal R}_1}),\,g,\,
({A_2},{\overline A_2},{{\cal R}_2})\right)$
are the same.
In this case corresponding transmutations are quantum braided groups in
$({\cal  C}_\cO{A_0,\overline A_0})
_\cO{\underline{A_1},\underline{\overline A_1}}$ and
in ${\cal  C}_\cO{A_1,\overline A_1}$ respectively and are the same if we
identify ${\cal  C}_\cO{A_1,\overline A_1}$ with a full subcategory of
$({\cal  C}_\cO{A_0,\overline A_0})
_\cO{\underline{A_1},\underline{\overline A_1}}$.
\end{theorem}

A transmutation is compatible with a braided group cross product
and generalized bosonization:

\begin{theorem}
Let $\left((A,\overline A,{\cal R}_A),f,H\right)$ be transmutation data
for bialgebra (braided group) $H$ in $\C$, and $B$ a bialgebra
(braided group) in $\DY{\cal C}^H_{H,\cO{A,\overline A}}$ and
$H\ltimes B$ cross product braided group with natural embedding
$i:\,H\hookrightarrow H\ltimes B$.
Then $\left((A,\overline A,{\cal R}_A),i\circ f,H\ltimes B\right)$
are transmutation data and
$\underline{H\ltimes B}=\underline H\ltimes B$ where the left hand side
is a transmutations of cross product and the right hand side is a cross
product of transmutation
$\underline H\in{\rm Obj}(\C_\cO{A,\overline A})$ with
$B$ considered as an object of
$\DY{\C_{\cO{A,\overline A}}}^{\underline H}_{\underline H}$.
\end{theorem}

\begin{theorem}
Let $\left((A,\overline A,{\cal R}_A),f,(H,\overline H,{\cal R}_H)\right)$
be transmutation data for a quasitriangular bialgebra
(quantum braided group)
$(H,\underline H,{\cal R}_H)$ in $\C$, and $(B,\underline B,{\cal R}_B)$
a quasitriangular bialgebra (quantum braided group)
in
\hfill
$\C_\cO{H,\overline H}$.
\hfill
And let
\newline
$(H\ltimes B,\underline H\ltimes\underline B,{\cal R}_{H\ltimes B})$
be generalized bosonization with natural embedding
$i:\,H\hookrightarrow H\ltimes B$.
Then
$\left((A,\overline A,{\cal R}_A),i\circ f,
(H\ltimes B,\underline H\ltimes\underline B,{\cal R}_{H\ltimes B})\right)$
are transmutation data and
\begin{displaymath}
(\underline{H\ltimes B},\underline{\overline H\ltimes\overline B},
{\cal R}_{\underline{H\ltimes B}}))=
(\underline H\ltimes B,\underline{\overline H}\ltimes\overline B,
{\cal R}_{\underline H\ltimes B})\,,
\end{displaymath}
where the left hand side (resp. the right hand side) is the result of the
composition of bosonization and transmutation (resp. transmutation and
bosonization).
\end{theorem}

\subsec{}
In the rest of the section we define a ribbon structure on quantum braided
group in a such way that the corresponding category of modules becomes
balanced and show that generalized bosonization and
transmutation are compatible with ribbon structures.

A braided category $\cal C$ is called {\em balanced} \cite{FY2}
if there exists an automorphism $\theta$ of the identical functor
${\rm Id}:\,{\cal C}\rightarrow{\cal C}$ such that
\begin{equation}
\theta_{X\otimes Y}=
\Psi_{Y,X}\circ\Psi_{X,Y}\circ (\theta_X\otimes\theta_Y)\,.
\end{equation}
(Existence of duals and compatibility between balancing and duality are
not essential for our considerations).

\begin{definition}
A quantum braided group $(H,\overline H,{\cal R}_H)$
in a balanced category $\cal C$ is {\em ribbon} if there exists
a group-like element
$\gamma:\,\underline 1\rightarrow H$
(i.e. $\Delta\circ\gamma=\gamma\otimes\gamma$),
which satisfies the identity
\begin{equation}
\label{spherical}
S_H^2\cdot\gamma=\gamma\cdot\theta_H
\end{equation}
\end{definition}

This is a natural generalization of an (ordinary) ribbon Hopf algebra
\cite{RT}.
And, on the other hand, one can consider the identity (\ref{spherical})
as a generalization of the axiom of spherical Hopf algebra \cite{BW}:\
$S^2(a)=\gamma\cdot a\cdot\gamma^{-1}$\ for any $a\in H$.

\begin{proposition}
If $(H,\overline H,{\cal R},\gamma)$ is a ribbon Hopf algebra
in a balanced category $\cal C$ then the category
${\cal C}_\cO{H,\overline H}$ is also balanced with
\begin{equation}
{}^{({\cal C}_\cO{H,\overline H})}\theta=
{}^{\cal C}\theta\circ(\triangleleft v)\,,
\qquad\qquad\hbox{where}\quad v:=(u\cdot\gamma)^{-1}\,.
\end{equation}
\end{proposition}

For quantum braided group cross product it is possible to formulate
the following analog of the lemma \ref{cross-sigma} without involving of
module.

\begin{lemma}
\label{cross-u}
Let $(A\ltimes B, \overline A\ltimes\overline B, {\cal R}_{A\ltimes B})$
be a cross product of quantum braided groups $(A,\overline A,{\cal R}_A)$
in $\cal C$ and
$(B,\overline B,{\cal R}_B)$ in ${\cal C}_{\cO{A,\overline A}}$.
Then
the 'element' $u_{A\ltimes B}:\,\underline 1\rightarrow A\ltimes B$
defined in \ref{element-u} is a product of corresponding
'elements' of $A$ and $B$:
\begin{equation}
u_{A\ltimes B}=u_A\cdot u_B=u_B\cdot u_A
\end{equation}
\end{lemma}

The next proposition is a corollary of the lemmas \ref{cross-square} and
\ref{cross-u}.

\begin{proposition}
Let $\cal C$ be a balanced category,
$(A,\overline A,{\cal R}_A,\gamma_A)$ be a ribbon braided group in $\cal C$
and $(B,\overline B,{\cal R}_B,\gamma_B)$ be a  ribbon braided group in
${\cal C}_\cO{A,\overline A}$.
Then the quantum group cross product
$(A\ltimes B,{\cal R}_{A\ltimes B},\gamma_{A\ltimes B})$
is a ribbon braided group in $\cal C$ with
\begin{equation}
\gamma_{A\ltimes B}=\gamma_A\cdot\gamma_B=\gamma_B\cdot\gamma_A
\end{equation}
\end{proposition}

\begin{lemma}
Let $(A,\overline A,{\cal R}_A)$,
$(H,\overline H,{\cal R}_H)$ be quantum braided groups in $\cal C$,
morphism $f:A\rightarrow H$ defines transmutation data and
$(\underline H,\underline{\overline H},{\cal R}_{\underline H})$
the transmutation.
Then the 'element' $u_{\underline H}$ defined in \ref{element-u}
for $\underline H$ is a quotient of corresponding elements of
$H$ and $A$:
\begin{equation}
u_{\underline H} ={u_A}^{-1}\cdot u_H=u_H\cdot{u_A}^{-1}
\end{equation}
\end{lemma}

\begin{proposition}
 Let $(A,\overline A,{\cal R}_A,\gamma_A)$,
$(H,\overline H,{\cal R}_H,\gamma_H)$ be ribbon braided groups in $\cal C$
and morphism $f:A\rightarrow H$ defines transmutation data.
Then the transmutation
$(\underline H,\underline{\overline H},{\cal R}_{\underline H})$
is a ribbon braided group with
\begin{equation}
\gamma_{\underline H} :=
{\gamma_A}^{-1}\cdot \gamma_H=\gamma_H\cdot{\gamma_A}^{-1}
\end{equation}
\end{proposition}

\subsec{}
More concrete applications of the theory that is developed here can
be find in \cite{B3} where we develop a fully braided analog of
Faddeev-Reshetikhin-Takhtajan construction of
quasitriangular bialgebra $A(X,R)$.
For given pairing $C$ factor-algebra $A(X,R;C)$ which is a dual quantum
braided group is built.
Corresponding inhomogeneous quantum group is obtained as a result of
generalized bosonization.
We define a fully  braided analog of Jur\v co construction of first order
bicovariant differential calculus.

\section* {Acknowledgments}
The author is pleased to thank the referees for the critical remarks which
do much toward improving the exposition.
The research described in this paper was made possible
by Grant No U4J200 from the International Science Foundation.

\begin{figure}
\begin{displaymath}
\matrix{\object{X}\step\object{Y}\Step\object{A}\Step\cr
\vvbox{\hbox{\id\step\id\step\cd}
       \hbox{\id\step\hx\step\cd}
       \hbox{\hx\step\hx\Step\id}
       \hbox{\id\step\ru\step\Ru}
       \hhbox{\krl\id\step\rd\step\Rd}
       \hbox{\hx\step\hx\Step\id}
       \hbox{\id\step\hx\step\cu}
       \hbox{\id\step\id\step\cu}}\cr
\object{X}\step\object{Y}\Step\object{A}\Step}
\quad ={}\quad
\vvbox{\hhbox{\krl\id\step\hstep\id\step\cd}
       \hbox{\id\step\hstep\hx\step\id}
       \hhbox{\krl\id\step\cd\hstep\id\step\id}
       \hbox{\hx\step\id\hstep\ru}
       \hbox{\id\step\ru\hstep\id}
       \hbox{\id\step\rd\hstep\id}
       \hbox{\hx\step\id\hstep\rd}
       \hhbox{\krl\id\step\cu\hstep\id\step\id}
       \hbox{\id\step\hstep\hx\step\id}
       \hhbox{\krl\id\step\hstep\id\step\cu}}
\quad ={}\quad
\vvbox{\hhbox{\krl\id\Step\id\step\cd}
       \hbox{\id\Step\hx\step\id}
       \hhbox{\krl\hrd\step\cd\hstep\ru}
       \hbox{\id\hstep\hx\step\id\hstep\id}
       \hhbox{\krl\hru\step\cu\hstep\rd}
       \hbox{\id\Step\hx\step\id}
       \hhbox{\krl\id\Step\id\step\cu}}
\quad ={}
\end{displaymath}
\begin{displaymath}
 ={}\quad
\vvbox{\hhbox{\krl\id\step\step\id\step\cd}
       \hbox{\id\step\step\hx\step\hdcd}
       \hbox{\id\step\step\id\step\hx\step\id}
       \hhbox{\krl\rd\step\id\step\id\step\ru}
       \hbox{\id\step\hx\step\id\step\id}
       \hhbox{\krl\ru\step\id\step\id\step\rd}
       \hbox{\id\step\step\id\step\hx\step\id}
       \hbox{\id\step\step\hx\step\hddcu}
       \hhbox{\krl\id\step\step\id\step\cu}}
\quad ={}\quad
\vvbox{\hbox{\id\step\step\id\step\cd}
       \hbox{\id\step\step\hx\step\hstep\hcd}
       \hhbox{\krl\rd\step\id\step\hrd\step\id\step\id}
       \hbox{\id\step\hx\step\id\hstep\hx\step\id}
       \hhbox{\krl\ru\step\id\step\hru\step\id\step\id}
       \hbox{\id\step\step\hx\step\hstep\hcu}
       \hbox{\id\step\step\id\step\cu}}
\quad ={}\quad
\vvbox{\hbox{\rd\step\rd\step\id}
       \hbox{\id\step\hx\step\id\hstep\hcd}
       \hhbox{\krl\id\step\id\step\cu\hstep\id\step\id}
       \hbox{\id\step\id\step\hstep\hx\step\id}
       \hhbox{\krl\id\step\id\step\cd\hstep\id\step\id}
       \hbox{\id\step\hx\step\id\hstep\hcu}
       \hbox{\ru\step\ru\step\id}}
\end{displaymath}
\caption{Proof of the crossed module axiom for tensor product of crossed
modules:
\protect\newline
We use coherence and associativty
(the first, the third, the fifth equalities),
the crossed module axiom for $X$ and $Y$ (the second and the fourth
equalities respectively).}
\label{Proof-DY-otimes}
\end{figure}

\begin{figure}
$$
{\rm L}_{{}_{A^{\rm op}}\DY{\overline{\cal C}}^{A^{\rm op}}}^{{}_-X}=\enspace
\vvbox{\hbox{\step\cd\step\id}
       \hbox{\cd\hstep\Obj{\enspace\mu^-_\ell}\hstep\lu}
       \hbox{\id\Step\Dash{{\rm R}_{{}_A\DY{\cal C}^A}^X}}
       \hbox{\Lu\step\Obj{\enspace\mu^-_\ell}\step\id}}
\quad\enspace=\enspace
\vvbox{\hbox{\step\cd\step\id}
       \hbox{\cd\hstep\Obj{\enspace\mu^-_\ell}\hstep\lu}
       \hbox{\id\Step\Dash{{\rm L}_{{}_A\DY{\cal C}^A}^X}}
       \hbox{\Lu\step\Obj{\enspace\mu^-_\ell}\step\id}}
\quad\enspace=
{\rm R}_{{}_{A^{\rm op}}\DY{\overline{\cal C}}^{A^{\rm op}}}^{{}_-X}
$$
{\scriptsize a) If $X\in{\rm Obj}({}_A\DY{\cal C}^A)$ then
${}_-X\in{\rm Obj}({}_{A^{\rm op}}\DY{\overline{\cal C}}^{A^{\rm op}})$.}
$$
\matrix{\object{X}\step\object{A}\Step\object{X^\vee}\cr
\vvbox{\hhbox{\krl\id\step\id\Step\id}
       \hbox{\id\step\dash{{\rm L}_{{}_A\DY{\cal C}^A}^{X^\vee}}}
       \hbox{\hev\Step\id}}\cr
\step\Step\object{A}}
\quad=\quad
\matrix{\object{X}\Step\object{A}\step\object{X^\vee}\cr
\vvbox{\hhbox{\krl\id\Step\id\step\id}
       \hbox{\dash{{\rm R}_{{}^A\DY{\cal C}_A}^X}\step\id}
       \hbox{\id\Step\hev}}\cr
\object{A}\Step\step}
\qquad\quad
\matrix{\object{X}\step\object{A}\Step\object{X^\vee}\cr
\vvbox{\hhbox{\krl\id\step\id\Step\id}
       \hbox{\id\step\dash{{\rm R}_{{}_A\DY{\cal C}^A}^{X^\vee}}}
       \hbox{\hev\Step\id}}\cr
\step\Step\object{A}}
\quad=\quad
\matrix{\object{X}\Step\object{A}\step\object{X^\vee}\cr
\vvbox{\hhbox{\krl\id\Step\id\step\id}
       \hbox{\dash{{\rm L}_{{}^A\DY{\cal C}_A}^X}\step\id}
       \hbox{\id\Step\hev}}\cr
\object{A}\Step\step}
$$
{\scriptsize b) If $X\in{\rm Obj}({}^A\DY{\cal C}_A)$ then
$X^\vee\in{\rm Obj}({}_A\DY{\cal C}^A)$.}
\caption{Constructions of new crossed modules.
         Proof of compatibility conditions.}
\label{Proof-new-DY}
\end{figure}

\begin{figure}
\begin{displaymath}
\matrix{\object{A\!\ltimes\! B}\Step\object{A\!\ltimes\! B}\cr
        \vvbox{\hbox{\cu}
               \hbox{\cd}}\cr
        \object{A\!\ltimes\! B}\Step\object{A\!\ltimes\! B}}
 =\enspace
\matrix{\hstep\object{A}\Step\step\object{B}
\step\hstep\object{A}\step\hstep\object{B}\hstep\cr
\vvbox{\hhbox{\hstep\id\Step\step\id\step\cd\step\id}
       \hbox{\hstep\id\Step\step\hx\step\id\step\id}
       \hhbox{\hstep\id\Step\hstep\dd\step\ru\step\id}
       \hhbox{\hstep\id\Step\dd\step\cd\step\cd}
       \hhbox{\cd\step\cd\step\id\step\id\step\hrd\hstep\id}
       \hbox{\id\step\hx\step\id\step\id\step\hx\hstep\id\hstep\id}
       \hhbox{\cu\step\cu\step\id\step\id\step\hru\hstep\id}
       \hhbox{\hstep\id\Step\d\step\cu\step\cu}
       \hhbox{\hstep\id\Step\hstep\d\step\rd\step\id}
       \hbox{\hstep\id\Step\step\hx\step\id\step\id}
       \hhbox{\hstep\id\Step\step\id\step\cu\step\id}}\cr
\hstep\object{A}\Step\step\object{B}
\step\hstep\object{A}\step\hstep\object{B}\hstep}
\enspace =\enspace
\vvbox{\hhbox{\cd\step\hstep\cd\step\cd\Step\hstep\cd}
       \hbox{\id\step\id\step\hstep\id\step\hx\step\d\step\hstep\id\step\d}
       \hhbox{\id\step\id\step\dd\hstep\cd\hstep\d\hstep\step\d\step\id\Step\d}
       \hbox{\id\step\id\step\hx\step\hdcd\d\step\id\step\rd\hstep\step\id}
       \hbox{\id\step\id\step\id\step\hx\step\id\step\id\step\hx\step
							     \hdcd\hstep\id}
       \hbox{\id\step\hx\step\id\step\k\step\hx\step\hx\step\id\hstep\id}
       \hbox{\id\step\id\step\id\step\hx\step\id\step\id\step\hx\step
							     \hddcu\hstep\id}
       \hbox{\id\step\id\step\hx\step\hddcu\dd\step\id\step\ru\hstep\step\id}
       \hhbox{\id\step\id\step\d\hstep\cu\hstep\dd\hstep\step\dd\step\id\Step
                                                                          \dd}
       \hbox{\id\step\id\step\hstep\id\step\hx\step\dd\step\hstep\id\step\dd}
       \hhbox{\cu\step\hstep\cu\step\cu\Step\hstep\cu}}
\quad =
\end{displaymath}
\begin{displaymath}
=\enspace
\vvbox{\hbox{\hcd\step\cd\step\hstep\hcd\step\cd}
       \hhbox{\id\step\id\step\rd\step\id\step\hstep\id\step\id\step\rd\step
									 \id}
       \hbox{\id\step\hx\step\id\step\id\step\hstep\id\step\hx\step\d\d}
       \hhbox{\id\step\id\step\cu\step\d\step\id\step\id\hstep\cd\step\cd
								 \hstep\id}
       \hbox{\id\step\id\step\hstep\d\step\hx\step\id\hstep\id\step\id\step\id
							\step\id\hstep\id}
       \hbox{\id\step\id\Step\hstep\hx\step\hx\hstep\id\step\hx\step\id
							       \hstep\id}
       \hbox{\id\step\id\step\hstep\dd\step\hx\step\id\hstep\id\step\id\step
                                                    \id\step\id\hstep\id}
       \hhbox{\id\step\id\step\cd\step\dd\step\id\step\id\hstep\cu\step\cu
								 \hstep\id}
       \hbox{\id\step\hx\step\id\step\id\step\hstep\id\step\hx\step\dd\dd}
       \hhbox{\id\step\id\step\ru\step\id\step\hstep\id\step\id\step\ru\step
									 \id}
       \hbox{\hcu\step\cu\step\hstep\hcu\step\cu}}
\enspace=\enspace
\matrix{\hstep\object{A\!\ltimes\! B}\Step\object{A\!\ltimes\! B}\hstep\cr
        \vvbox{\hhbox{\cd\step\cd}
               \hbox{\id\step\hx\step\id}
               \hhbox{\cu\step\cu}}\cr
         \hstep\object{A\!\ltimes\! B}\Step\object{A\!\ltimes\! B}\hstep}
\end{displaymath}
\caption{Proof of the bialgebra axiom for $A\ltimes B$:
\protect\newline
The first equality follows from the bialgebra axiom for $A$ and $B$.
The second uses the fact that $B$ is a coalgebra in ${\cal C}_A$
and an algebra in ${\cal C}^A$.
The third uses the crossed module axiom.
The last is the bialgebra axiom for $A$.}
\label{Proof-Bialg}
\end{figure}

\begin{figure}
\begin{displaymath}
\vvbox{\hbox{\ro{\overline{\cal R}_H}\step\ro{{\cal R}_H}}
       \hbox{\id\Step\hddcu\Step\id}
       \hbox{\x\Step\dd}
       \hbox{\id\Step\cu}}
\enspace=\enspace
\vvbox{\hbox{\Ro{\overline{\cal R}_B}\step\Ro{{\cal R}_B}}
      \hbox{\id\step\ro{\overline{\cal R}_A}\step\hcu\step\ro{{\cal R}_A}
                                               \step\id}
      \hhbox{\krl\cu\Step\d\step\id\step\cd\step\hstep\cu}
      \hbox{\hstep\d\Step\id\step\hxx\step\id\Step\id}
      \hbox{\step\hstep\d\step\hddcu\step\ru\step\dd}
      \hbox{\Step\hstep\hdcu\step\dd\step\dd}
      \hbox{\Step\step\hstep\hx\step\dd}
      \hhbox{\krl\Step\step\hstep\id\step\cu}}
\enspace=\enspace
\vvbox{\hbox{\ro{\overline{\cal R}_B}\step\Ro{{\cal R}_B}}
       \hbox{\d\step\hdcu\step\ro{{\cal R}_A}\step\id}
       \hbox{\step\d\step\ru\step\dd\dd}
       \hbox{\Step\hx\step\dd\dd}
       \hhbox{\krl\Step\id\step\cu\hstep\dd}
       \hhbox{\krl\Step\id\step\hstep\cu}}
\enspace=
\end{displaymath}
\begin{displaymath}
=\enspace
\vvbox{\hbox{\ro{\overline{\cal R}_B}\step\Ro{{\cal R}_B}}
       \hbox{\d\step\hdcu\step\ro{{\cal R}_A}\step\id}
       \hbox{\step\d\step\ru\step\dd\dd}
       \hbox{\Step\hx\step\dd\dd}
       \hbox{\step\dd\step\ru\dd}
       \hbox{\dd\ro{{\cal R}_A}\hdcu}
       \hbox{\ru\Step\hcu}}
\enspace=\qquad\quad
\vvbox{\hbox{\DDash{\overline{\cal R}_B^{\rm op}\cdot{\cal R}_B}}
       \hbox{\hhstep\hhstep\id\step\ro{{\cal R}_B}\step\id}
       \hbox{\hhstep\hhstep\ru\Step\hcu}}
\quad\enspace=\epsilon_H\otimes\epsilon_H
\end{displaymath}
\caption{Proof of the identity
$\overline{\cal R}^{\rm op}_{A\ltimes B}\cdot{\cal R}_{A\ltimes B}=
  \epsilon_H\otimes\epsilon_H$:
\protect\newline
The first and the third identities use the formula for multiplication in
$\overline A\ltimes\overline B$ and $A\ltimes B$ respectively.
The second and the last use the fact that
$\overline{\cal R}\cdot{\cal R}=\epsilon\otimes\epsilon$ for $A$ and for $B$
respectively.}
\label{randinverse}
\end{figure}

\begin{figure}
\begin{displaymath}
\matrix{
\vvbox{\hbox{\rh}
       \hbox{\id\step\hstep\hcd}}\cr
\object{H}\step\hstep\object{H}\step\object{H}}
\enspace=\enspace
\matrix{
\vvbox{\hbox{\Step\Rb}
       \hbox{\step\dd\Rb\d}
       \hbox{\dd\dd\step\ra\step\d\d}
       \hbox{\id\step\id\step\dd\ra\d\step\rd\d}
       \hbox{\hcu\step\hcu\Step\id\step\hx\step\id\step\id}
       \hbox{\hstep\cu\hstep\Step\id\step\id\step\hcu\step\id}}\cr
\hstep\step\object{H}\Step\step\hstep\object{A}\step\object{B}\step
\object{A}\step\object{B}}
\enspace=\enspace
\matrix{
\vvbox{\hbox{\Step\Rb}
       \hbox{\step\dd\Rb\d}
       \hbox{\dd\dd\step\ra\step\d\d}
       \hbox{\id\step\id\step\dd\ra\d\step\id\step\id}
       \hbox{\hcu\step\hcu\Step\id\step\id\step\id\step\id}
       \hbox{\hstep\cu\hstep\Step\id\step\id\step\id\step\id}}\cr
\hstep\step\object{H}\Step\step\hstep\object{A}\step\object{B}\step
\object{A}\step\object{B}}
\end{displaymath}
{\scriptsize a) ${\cal R}_H$ is an algebra-coalgebra copairing}
\begin{displaymath}
\matrix{
\vvbox{\hbox{\Rb}
       \hbox{\id\step\ra\step\id}
       \hbox{\id\step\hdcd\step\id\step\rd}
       \hbox{\hx\step\id\step\hx\step\id}
       \hbox{\id\step\ru\step\id\step\hcu}}
\enspace=\enspace
\vvbox{\hbox{\Ra}
       \hbox{\id\step\rb\step\id}
       \hbox{\id\step\id\Step\id\step\id}}
&\qquad\qquad&
\matrix{\hstep\object{A}\Step\hstep\object{B}\Step\hstep\cr
        \vvbox{\hbox{\Obj{{\cal R}_A\cdot\Delta_A}\hcd\step\cd
                                        \Step\Obj{{\cal R}_B\cdot\Delta_B}}
               \hbox{\id\step\id\step\id\hstep\ra\hhstep\d}
               \hhbox{\id\step\id\step\hru\step\hstep\dd\hstep\id}
               \hbox{\id\step\hx\step\dd\step\id}
               \hhbox{\id\step\id\step\cu\Step\id}}\cr
	\object{A}\step\object{B}\step\object{A}\Step\hstep\object{B}\step}
\cr
\hbox{\scriptsize b)The auxiliary identity}
&&
\hbox{\scriptsize c)
$\overline\Delta^{\rm op}\cdot{\cal R}={\cal R}\cdot\Delta$}
}
\end{displaymath}
\caption{Proof of the quantum braided group axioms:
\protect\newline
In the part a):
the first equality use the fact that ${\cal R}_A$ and ${\cal R}_B$ are
algebra-coalgebra pairings.
The second follows from the auxiliary identity from the  part b)
of this figure.
}
\label{proofpairing}
\end{figure}

\end{document}